\documentclass[11pt,reqno,letter]{amsart}

\usepackage{amssymb}
\usepackage[dvipdfmx]{graphicx}

\usepackage[usenames]{color}
\usepackage{caption}
\usepackage[dvips]{psfrag}
\usepackage{latexsym}
\usepackage{amsmath}
\usepackage{amsthm}
\usepackage{tabularx}
\usepackage{bm}
\usepackage{lscape}
\usepackage[numbers]{natbib}
\usepackage[usenames]{color}
\usepackage{url}
\usepackage{hyperref}
\usepackage{upgreek}
\usepackage{units}
\usepackage{fancyhdr}
\definecolor{OliveGreen}{rgb}{0,0.6,0}

\usepackage{enumitem}

\usepackage[normalem]{ulem}

\theoremstyle{plain}
\newtheorem{theorem}{Theorem}[section]
\newtheorem{lemma}[theorem]{Lemma}

\newtheorem{proposition}[theorem]{Proposition}

\theoremstyle{definition}
\newtheorem{definition}[theorem]{Definition}
\newtheorem{remark}[theorem]{Remark}
\newtheorem{example}[theorem]{Example}
\newtheorem{assumption}[theorem]{Assumption}

\numberwithin{equation}{section}
\numberwithin{theorem}{section}
\numberwithin{table}{section}
\numberwithin{figure}{section}

\DeclareMathOperator*{\argmax}{arg\,max}
\DeclareMathOperator*{\argmin}{arg\,min}

\title{Robustness Measures in Distributionally Robust Optimization}\thanks{An earlier version of this paper had the title ``Distributionally Robust Optimization is a Multi-Objective Problem".}

\date{\today}

\author[Gotoh]{Jun-ya Gotoh\textsuperscript{$\dagger$}}

\author[Kim]{Michael Jong Kim\textsuperscript{$\ddagger$}}

\author[Lim]{Andrew E.B. Lim\textsuperscript{$*$}}

\dedicatory{\textsuperscript{$\dagger$}Department of Data Science for Business Innovation, Chuo University, Tokyo, Japan. Email: jgoto@kc.chuo-u.ac.jp \\ \textsuperscript{$\ddagger$}Sauder School of Business, University of British Columbia, Vancouver, Canada. Email: mike.kim@sauder.ubc.ca \\ \textsuperscript{$*$}Department of Analytics and Operations, Department of Finance, and Institute of Operations Research and Analytics, National University of Singapore, Singapore. Email: andrewlim@nus.edu.sg}

\begin{document}

\begin{abstract}

Distributionally Robust Optimization (DRO) is  a worst-case approach to decision making when there is model uncertainty.   It is also well known that for certain uncertainty sets, DRO is approximated by a regularized nominal problem. We show that the regularizer  is not just a penalty function but  the worst-case sensitivity (WCS) of the expected cost with respect to deviations  from the nominal model,  giving it the interpretation of a {\it robustness measure}. This   has substantial consequences for robust modeling. It  shows that DRO is fundamentally a tradeoff between performance and robustness, where the robustness measure is determined by the uncertainty set. The  robustness measure reveals properties of a cost distribution that affect sensitivity to misspecification. This leads to  a systematic approach to selecting uncertainty sets. The family of DRO solutions obtained by varying the size of the uncertainty set traces a near Pareto-optimal performance--robustness frontier that can be used  to select its size. The frontier identifies problem instances where the price of robustness is high and provides insight into effective ways of redesigning the system to reduce this cost. We derive WCS for a collection of commonly used uncertainty sets, and illustrate these ideas in a number of applications.



\end{abstract}

\maketitle


\section{Introduction}

It is  a quirk of the literature on Distributionally Robust Optimization (DRO) that robustness, the quality being ``optimized,"  is not even measured. Contrast this with ``risk sensitive" optimization where quantitative risk measures  play a central role, from comparing the riskiness of decisions and  defining tradeoffs between performance and risk while inviting consideration of the properties a risk measure should possess. A quantitative robustness measure opens the door to similar questions in  robust optimization. Without one, even deciding  which of two decisions is ``more robust" is  challenging.

Though formulated as a worst-case problem, DRO can be  approximated for certain uncertainty sets  by a (non-worst-case) regularized nominal problem \cite{bertsimas2018,blanchet2019,Duchi2017variance,esfahani2018data,gao2023distributionally, gao2024wasserstein,kuhn2019wasserstein}. We show that this extends to a broader class of uncertainty sets, and that the induced regularizer is the sensitivity of expected cost to worst-case perturbations of the nominal model, and hence a measure of robustness. Our main contribution, however, is not the representation per se, but the modeling consequences of interpreting the regularizer (worst-case sensitivity (WCS)) as a robustness measure. 
This includes the selection of uncertainty-set families, calibration of size,  and system design for improving the tradeoff between performance and robustness. 


In the broader literature, expected cost has been the primary focus of DRO. Uncertainty sets are often tuned to minimize the out-of-sample expected cost and out-of-sample performance guarantees are provided only for the expected cost. DRO has been criticized for being ``too conservative" because it optimizes the worst case.  Without a robustness measure to quantify the benefits of being conservative this is difficult to defend.
WCS addresses this gap by explicitly quantifying robustness.  
This enables us to identify the errors in the model that have the highest impact on the performance of the nominal solution, and provides a principled approach to selecting uncertainty sets that target these vulnerabilities. It guards against DRO models with solutions  that  are even more fragile than the nominal. The family of DRO solutions obtained by varying the size of the uncertainty set defines a near Pareto-optimal cost--sensitivity frontier. This can be used not only to select the size of the uncertainty set but  to identify {\it conservative problem instances} where the price of robustness is high. Since it measures a quality (sensitivity) that is meaningful to  minimize, WCS provides insight into  redesigning the system to reduce this cost and improve the performance--sensitivity tradeoff, including the introduction of ``hedging strategies" that  target sensitivity reduction.

\subsection*{Summary of contributions}

We illustrate the consequences of interpreting the regularizer as a robustness measure and  show how this can guide system design, uncertainty set selection, and the calibration of its size.

\vspace{-0.25cm}
\begin{itemize}[leftmargin=0.2in]
\item We define the notion of worst-case sensitivity (WCS), a quantitative measure of robustness defined by the uncertainty set. 
We show that worst-case sensitivity is a {\it generalized measure of deviation} and hence a measure of the spread of the cost distribution. 
\item 
We derive explicit expressions of WCS for a number of  uncertainty sets including smooth $\phi$-divergence, total variation, hard constraints on the likelihood ratio, 
and the Wasserstein distance. 
This enables us to select uncertainty sets according to the sensitivities that need to be controlled. Cost-sensitivity frontier plots can be used to calibrate the size of the uncertainty set and identify conservative problem instances where the price of robustness is high.
\item Hedging for robustness. Explicit expressions for worst-case sensitivity 
can be used to guide system redesign when  the price of robustness is high. We illustrate this in an inventory problem where ``return contracts" are shown to directly target sensitivity and improve the cost--sensitivity tradeoff.
\item We extend these ideas to derive worst-case sensitivity for robust conditional-value-at-risk (RCVaR), a risk measure that can be sensitive to model error \cite{BLS}. WCS  measures the spread of the tail of the loss distribution, where its precise form depends on the choice of uncertainty set.
\end{itemize}

For clarity we focus on single-period problems: this is enough to convey the key ideas and is the setting for much of the DRO literature.  Generalizations 
(multi-period problems, online problems with learning, etc) build on the intuition developed in this paper and will be reported elsewhere.

\subsection*{Overview of paper}
We review the relevant literature in Section \ref{sec:LitRev}. 
In Section \ref{sec:DRO_sensitivity}, we introduce worst-case sensitivity as a robustness measure, show (under mild assumptions) that it is a measure of spread, and that solutions of DRO problems define a nearly-Pareto optimal tradeoff between mean cost and worst-case sensitivity.
We derive explicit expressions for worst-case sensitivity for a number of different uncertainty sets in Section \ref{sec:WCS}. This enables us to identify the measure of robustness being controlled by each choice of uncertainty set.
Using two applications (inventory control and CVaR minimization), we show in Section  \ref{sec:Examples}  how  the  sensitivity measures 
we derive can be used  in robust modeling and design. This includes selection of uncertainty sets and their size. 
When the price of robustness is high, the sensitivity measure guides the design of hedging instruments that improve the performance--robustness tradeoff.
It is a  good starting point for the first-time reader.




\section{Literature Review}
\label{sec:LitRev}

Several papers \cite{bertsimas2018,blanchet2019,Duchi2017variance,esfahani2018data,gao2023distributionally, gao2024wasserstein,kuhn2019wasserstein} show the connection between DRO and regularized empirical optimization.  For example,  \cite{Duchi2017variance} shows that DRO with $\phi$-divergence uncertainty sets is almost the same as variance regularization  and \cite{blanchet2019,gao2023distributionally, gao2024wasserstein} show that worst-case regression and classification with the Wasserstein uncertainty set are equivalent to regularized sample average approximation (SAA). A common perspective in all these papers  is that the worst-case problem and the regularized approximation are regarded as  single-objective problems; performance bounds are established for the expected cost, experiments compare expected cost. In this paper, we also start with the regularized SAA. However, there is an immediate diversion in that we interpret the regularizer not as a penalty but a measure of robustness (justified by showing it is equal to worst-case sensitivity) and hence an objective function equally important as the expected cost. We explore the consequences of interpreting the regularizer as a robustness measure for robust modeling, including the selection of uncertainty sets and system design.

The choice of uncertainty set is critical and can result in very different solutions. We consider uncertainty sets of the form ${\mathcal Q}(\varepsilon)    =  \{{\mathbb Q}:{\mathcal F} \rightarrow {\mathbb R} ~|~d({\mathbb Q}\,| \,{\mathbb P})\leq\varepsilon\}$ where ${\mathbb P}$ is the nominal model, $\mathbb Q$ an alternative model, and $d({\mathbb Q}\,| \,{\mathbb P})$ is a divergence between $\mathbb P$ and $\mathbb Q$. The robustness measure is the worst-case sensitivity of the expected cost when the size of the uncertainty set $\varepsilon$ vanishes. Both the ``family"  (i.e., divergence $d$) and ``size" $\varepsilon$ are important. In the literature, uncertainty sets, particularly the divergence $d$, are often  model primitives. However, there is little insight how to select between uncertainty set families.  

The regularized SAA approximation shows that DRO induces a tradeoff between expected cost and the robustness measure (WCS) induced by the uncertainty set. The uncertainty set family can then be chosen according to the sensitivity control most appropriate for the problem. 
To our knowledge, this approach to selecting the uncertainty set family  has not been discussed in the literature. 

With regard to its size $\varepsilon$,  several approaches have been proposed. It can be chosen so that the uncertainty set  contains the true generative model with high probability \cite{ben2013robust,BGK,BGK-SAA,del2010,duchi2021statistics,gup1,lam2016robust,LamZhou} where the confidence levels (e.g., $90\%$, $95\%$, or $99\%$) and the uncertainty set family are primitives of the model. It can be  selected so that the in-sample worst-case expected cost upper bounds the out-of-sample expected cost with high probability \cite{esfahani2018data,kuhn2019wasserstein} or that the worst-case  in-sample cost is below an acceptable user-specified target \cite{Brown2012satisfycing}. One concern with having confidence levels  as  primitives is that they are chosen independent of the data and the objective function is that there is no reason to expect that a $95\%$  confidence level, say, has anything to do with good out-of-sample performance. Indeed, examples in \cite{gotoh2021calibration} show that out-of-sample expected cost with confidence levels like $90\%$, $95\%$, or $99\%$ can be high and  that $42\%$ is the answer in many applications   \cite{adams1995hitch}.

Another approach is to optimize an estimate of the out-of-sample expected cost using cross-validation or the bootstrap, extending the methods used to calibrate regularized regression models in Machine Learning \cite{HTF}. Once again the key difference is that the expected cost is the only objective. With a robustness measure, the decision maker (DM) can consider a tradeoff between performance and robustness, and may accept a larger expected cost in return for substantial sensitivity reduction. Conversely, expected cost and WCS can be used to identify problem instances where the performance--robustness tradeoff is unfavorable, and the specific form of the sensitivity measure used to redesign the system and improve the tradeoff. This is only possible when the regularizer is interpreted as an objective function. Finally, while it is possible that the DRO solution has a smaller out-of-sample expected cost  than the nominal solution \cite{anderson2022improving,gotoh2021calibration,kim2015MAB}, there is no guarantee and any out-performance is likely to be small \cite{Lam2}.  It is always possible to reduce the SAA out-of-sample cost if optimistic  \cite{Chen2023optimistic,gotoh2023data,Xie2024optimistic,Nguyen2019Optimistic,Norton2017optimistic} as well as worst-case solutions are considered. However, being optimistic increases the sensitivity of the expected cost to model misspecification \cite{gotoh2023data}, which undermines our original robustness goals.

The papers \cite{gotoh2018robust, gotoh2021calibration,gotoh2023data} explore performance--robustness tradeoffs for DRO problems with smooth $\phi$-divergence uncertainty sets, showing that WCS is the variance and that the associated mean--variance frontier can be used to select the size of the uncertainty set. The scope of the present paper is larger. Specifically: (i) WCS is defined for a larger class of uncertainty sets, not just smooth $\phi$-divergence, and shown under relatively mild assumptions to be a {\it generalized measure of deviation};  (ii) the  performance--sensitivity tradeoff studied in  \cite{gotoh2018robust,gotoh2021calibration} is shown to hold for a much broader class of uncertainty sets; (iii) WCS is derived for a collection of uncertainty sets leading to a systematic approach to selecting both the family and size of the uncertainty set (\cite{gotoh2018robust, gotoh2021calibration}  only consider size); and (iv) WCS can be used to design systems that have favorable performance--robustness tradeoffs. 

Worst-case sensitivity is derived for smooth $\phi$-divergence in  \cite{gotoh2018robust,gotoh2021calibration,lam2016robust}, total variation in \cite{charalambous2013TV}, and for the Wasserstein metric in \cite{BDOW}. It is derived for likelihood (but not prior) misspecification in a Bayesian model in \cite{shapiro2023bayesian}. \cite{gotoh2018robust,gotoh2021calibration} use it to calibrate the size of smooth $\phi$-divergence uncertainty sets while \cite{BDOW,lam2016robust,shapiro2023bayesian} do not explore implications for robust modeling (selecting family and size of uncertainty sets, system design, etc.).

\section{Worst-case sensitivity}
\label{sec:DRO_sensitivity}


We introduce worst-case sensitivity (WCS) as a measure of robustness.
We generalize the well-known approximation of DRO as a regularized nominal problem to a larger class of uncertainty sets and show that the regularizer is WCS. This highlights the tradeoff between performance and robustness at the core of DRO. Our main contribution is understanding the consequences for robust design, including uncertainty set selection and calibration.

\subsection{Worst-case sensitivity}

\subsubsection*{Nominal model.} Consider a cost function $f(x, Y)$ where $x$ is the decision and $Y$ is a  random variable with (nominal) distribution ${\mathbb P}$. For example, $\mathbb P$ could be the empirical distribution associated with a historical sample of $Y$'s generated {\it iid} from some unknown distribution. We would like to find a decision that performs well out-of-sample. 

A natural candidate is the solution of the nominal/Sample Average Approximation (SAA) problem
\begin{align}
\min_x ~ {\mathbb E}_{\mathbb P}[f(x, Y)].
\label{eq:pop_opt}
\end{align}
However, out-of-sample performance of this decision can be poor if the  real-world distribution of $Y$ is different from the nominal $\mathbb P$ and the expected reward under the nominal is sensitive to these errors. We would like to modify \eqref{eq:pop_opt} to account for misspecification. 

\subsubsection*{\it Worst-case sensitivity} 

While it is customary at this point to introduce the DRO model, we adopt a different approach by first introducing a robustness measure. It quantifies the central object in DRO (``robustness") and will be valuable when selecting uncertainty sets and evaluating DRO solutions.

Let ${\mathcal Q}(\varepsilon)$ be an {\it uncertainty set} of size $\varepsilon$, specifically,    a set of probability distributions increasing in $\varepsilon$ that contains the nominal distribution for every $\varepsilon\geq 0$ and degenerates to the nominal when $\varepsilon=0$. The worst-case expected cost with respect to ${\mathcal Q}(\varepsilon)$ is
\begin{align}
V(\varepsilon; f(x,\cdot)) &:=
\max_{{\mathbb Q}\in{\mathcal Q}(\varepsilon)} \mathbb{E}_{\mathbb{Q}}[f(x, Y)].
\label{eq:V}
\end{align}
To ease the notation, we often write $V(\varepsilon, x) \equiv V(\varepsilon; f(x,\cdot))$.  

\begin{figure}[h]
\begin{center}
\includegraphics[scale=0.07]{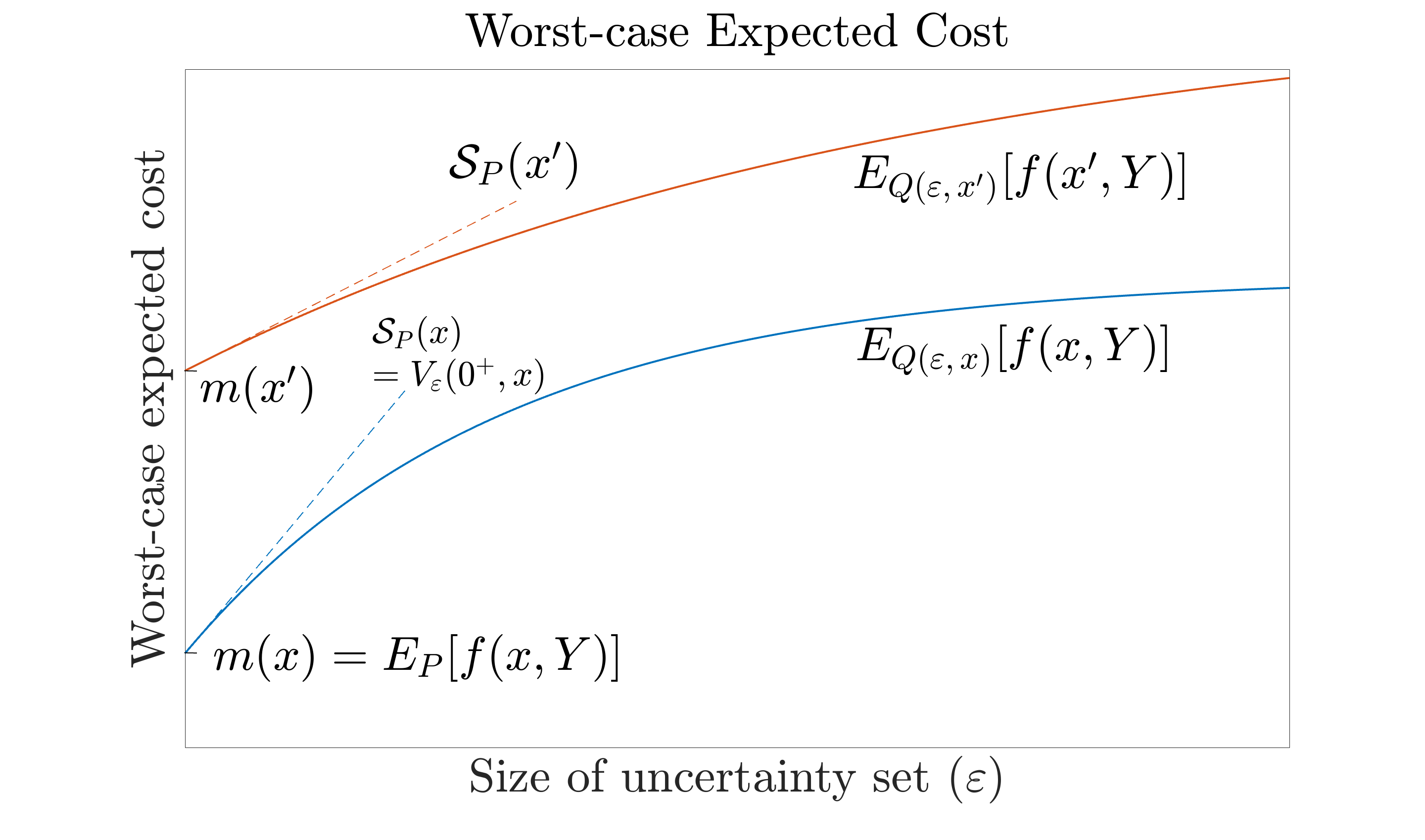} 
\includegraphics[scale=0.07]{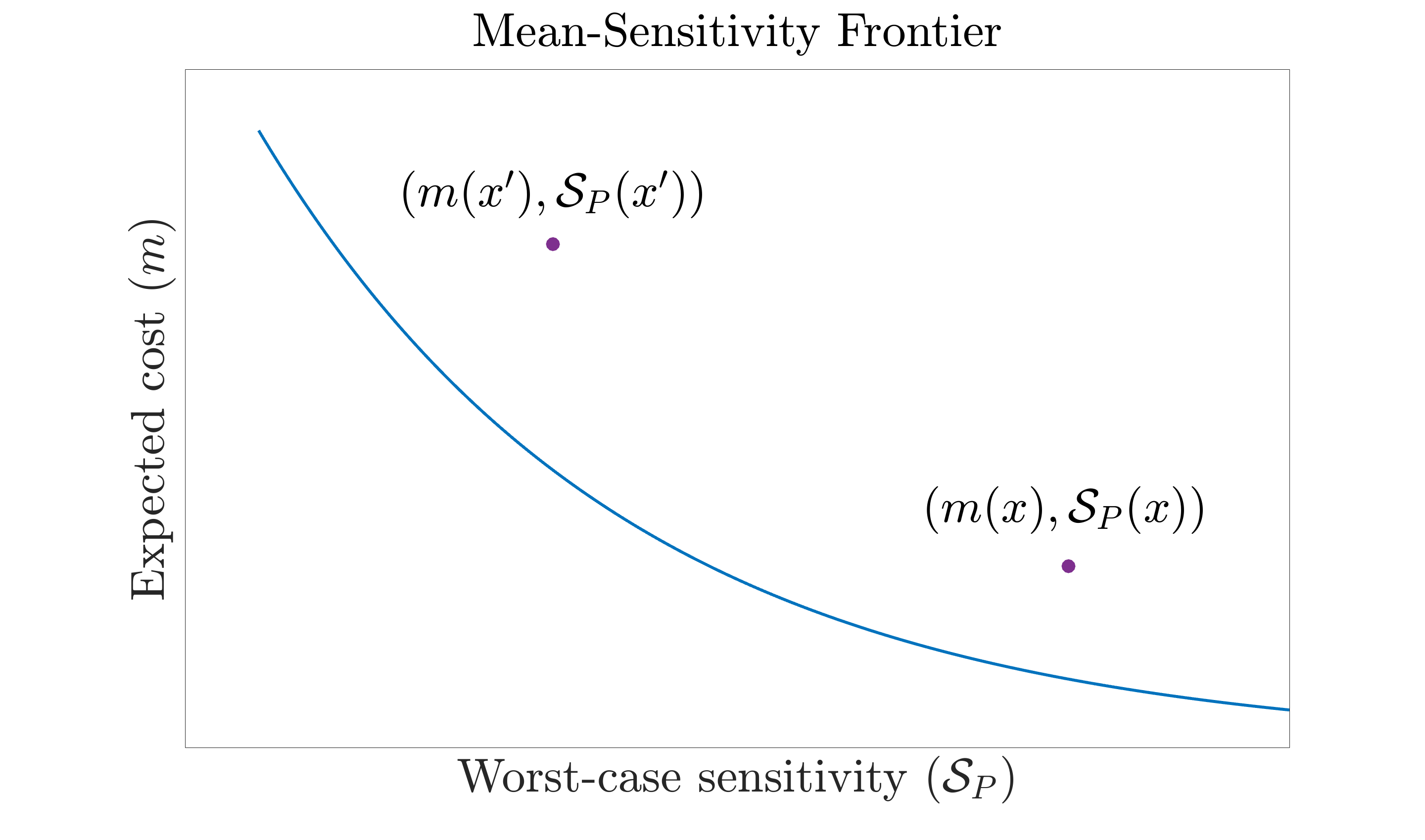} 
\end{center}
\caption{WCS and mean--sensitivity frontier}
\label{fig:illus}
{\raggedright\footnotesize
The figure on the left shows worst-case expected cost  as the size $(\varepsilon$) of the uncertainty set increases for decisions $x$ and $x'$. The intercept $(m)$ at $\varepsilon=0$ is the expected cost under the nominal $\mathbb P$; worst-case sensitivity $({\mathcal S}_{\mathbb P})$ is its slope. Every mean--sensitivity pair maps to a point  on the right. In this example, the expected cost under the nominal is smaller for $x$   ($m(x)< m(x')$) but is more sensitive to misspecification (${\mathcal S}_{\mathbb P}(x)>{\mathcal S}_{\mathbb P}(x')$). Pareto efficient decisions 
define a mean--sensitivity frontier shown on the right.
\par
}
\end{figure}

Suppose  there is a non-negative, concave and increasing function\footnote{The function $g(\varepsilon)$ depends on the uncertainty set. $V(\varepsilon; f(x, \cdot))$, and hence the difference ${\mathcal A}\big(\varepsilon; f(x, \cdot)\big) :=V(\varepsilon, x) - {\mathbb E}_{\mathbb P}[f(x, Y)]$, is concave in $\varepsilon$ if the set of alternative probability measures can be written in the form ${\mathcal Q}(\varepsilon) = \{{\mathbb Q}\,|\, d({\mathbb Q}|{\mathbb P})\leq \varepsilon\}$ for some convex function $d$ of ${\mathbb Q}$, which is true for all models considered in this paper.  $V(\varepsilon)-V(0)$ is $O(\sqrt \varepsilon)$  and $g(\varepsilon) = \sqrt{\varepsilon}$ for smooth $\phi$-divergence (Section \ref{sec:phi-div}), and $O(\varepsilon$) for budgeted uncertainty so $g(\varepsilon)=\varepsilon$ (Section \ref{sec:CVaR}). } $g(\varepsilon)$ such that 
\begin{eqnarray*}
V(\varepsilon, x) - {\mathbb E}_{\mathbb P}[f(x, Y)]  \sim O(g(\varepsilon))
\end{eqnarray*}
for every fixed  $x$.
Worst-case sensitivity 
\begin{eqnarray}
{\mathcal S}_{\mathbb P}[f(x, \cdot)] =  \begin{displaystyle} \lim_{\varepsilon\downarrow 0}\frac{V(\varepsilon, \, x)-{\mathbb E}_{\mathbb P}[f(x, Y)]}{g(\varepsilon)} \end{displaystyle}
\label{eq:sensitivity-general2}
\end{eqnarray}
measures the rate the expected cost increases under small perturbations of the nominal model. A higher sensitivity means the decision is more fragile to errors in the model. For decisions $x$ and $x'$,  worst-case sensitivity  is illustrated in the plot of the left in Figure \ref{fig:illus} being equal to the slope of the worst-case expected cost function at $\varepsilon=0$. The expected cost under $x$ is lower ($m<m')$ than under $x'$; sensitivity under $x$ is however higher (${\mathcal S}>{\mathcal S}'$).

WCS is determined by the  uncertainty set ${\mathcal Q}(\varepsilon)$ that defines the worst-case expected cost $V(\varepsilon, x)$. For example,  it is equal to the standard deviation of the cost when the uncertainty set defined in terms of smooth $\phi$-divergence (
Section \ref{sec:phi-div}), and the difference between the maximum and minimum values of the support of the reward distribution for total-variation deviation (Section \ref{sec:TV}).  Other uncertainty sets are also considered in Section \ref{sec:WCS}.

{\it Distributionally Robust Optimization (DRO)} accounts for misspecification in $\mathbb P$ by optimizing worst-case expected cost over a set of probability measures ${\mathcal Q}(\varepsilon)$ that includes $\mathbb P$
\begin{align}
\min_x V(\varepsilon; f(x,\cdot)) &\equiv \min_x \max_{{\mathbb Q}\in{\mathcal Q}(\varepsilon)} \mathbb{E}_{\mathbb{Q}}[f(x, Y)].
\label{eq:DRO}
\end{align}


Since   ${\mathcal Q}(\varepsilon)$ is increasing in $\varepsilon$ and contain the nominal $\mathbb P$, the worst-case expected cost \eqref{eq:V} is monotonically increasing in $\varepsilon$ so
\begin{eqnarray*}
V(\varepsilon, x) = {\mathbb E}_{\mathbb P}[f(x, Y)] + {\mathcal A}\big(\varepsilon; f(x, \cdot)\big)
\end{eqnarray*}
where the {\it ambiguity cost}
 ${\mathcal A}\big(\varepsilon; f(x, \cdot)\big)$ is a non-negative increasing function of $\varepsilon$ such that ${\mathcal A}\big(0; f(x, \cdot)\big)=0$. 
For conveniently chosen uncertainty sets,  ${\mathcal A}\big(\varepsilon; f(x, \cdot)\big)$ can be characterized through the dual problem of \eqref{eq:V}.  
By \eqref{eq:sensitivity-general2} we have 
\begin{eqnarray*}
{\mathcal A}\big(\varepsilon; f(x, \cdot)\big) =V(\varepsilon, x)  -  {\mathbb E}_{\mathbb P}[f(x, Y)] =  g(\varepsilon)\,{\mathcal S}_{\mathbb P}[f(x,\cdot)] + o(g(\varepsilon))
\end{eqnarray*}
so the worst-case problem \eqref{eq:DRO}  is locally equivalent to trading off expected cost and worst-case sensitivity:
\begin{eqnarray}
\label{eq:DRO-mean-sensitivity}
\min_x \max_{{\mathbb Q}\in{\mathcal Q}(\varepsilon)} \mathbb{E}_{\mathbb{Q}}[f(x, Y)] =  \min_x {\mathbb E}_{\mathbb P} [f(x,\,Y)] + g(\varepsilon)\,{\mathcal S}_{\mathbb P}[f(x,\cdot)] + o(g(\varepsilon)), \; \varepsilon>0.
\end{eqnarray}
While this extends existing ``regularized--nominal" approximations to a broader class of uncertainty sets, it also shows  that the regularizer is not just a penalty function but coincides with worst-case sensitivity and hence can be interpreted as a robustness measure.  This has consequences for robust modeling including uncertainty set selection and calibration and system design.


\begin{figure}[h]
\begin{center}
\includegraphics[scale=0.1]{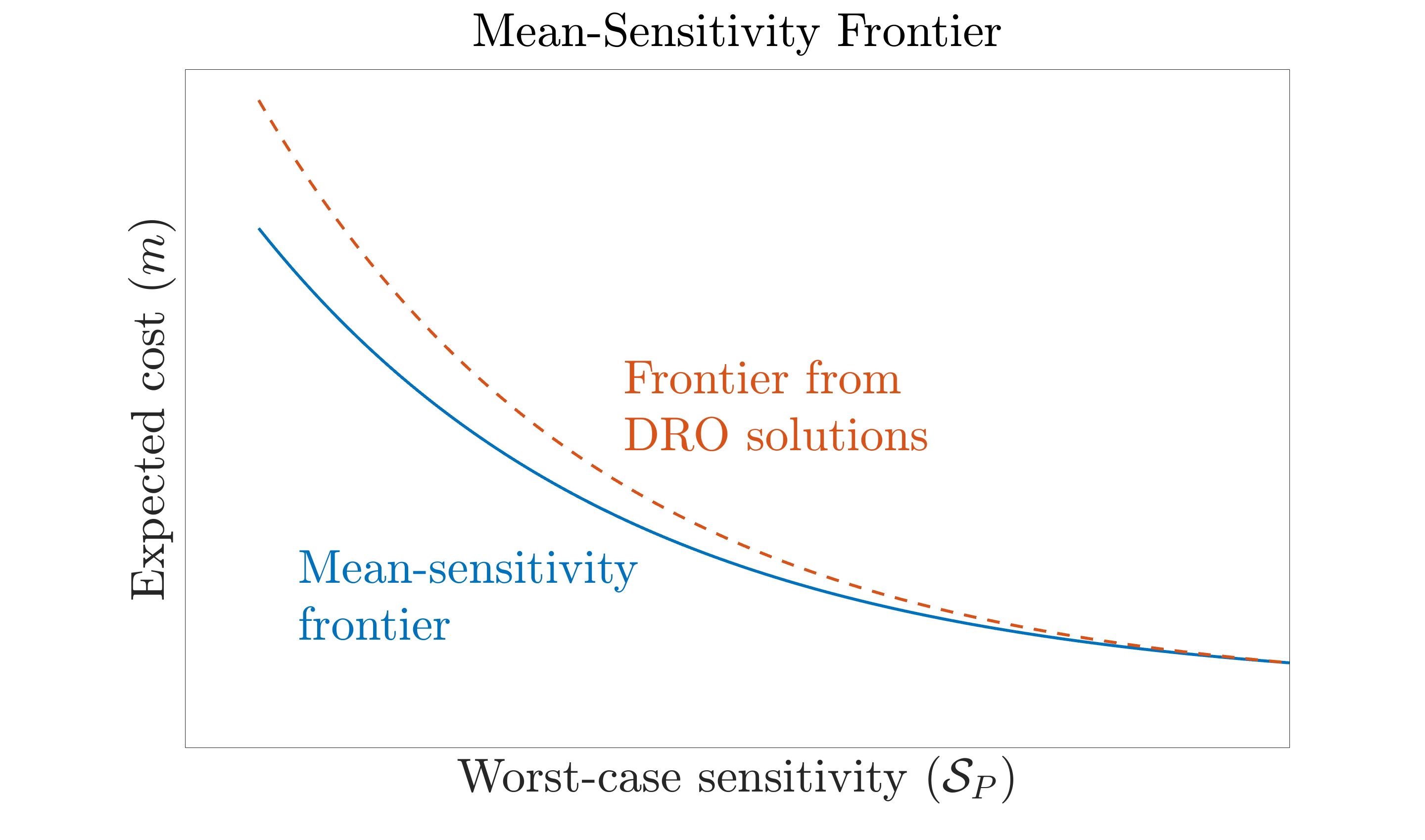} 
\end{center}
\caption{Solutions of the DRO problem \eqref{eq:DRO-mean-sensitivity} are nearly Pareto optimal for a tradeoff between expected cost and worst-case sensitivity.
}
\label{fig:illus2}
\end{figure}


\subsection{Sensitivity is a measure of deviation}


It follows from \eqref{eq:DRO-mean-sensitivity} that DRO makes a tradeoff between performance and robustness. For this reason, it is important to understand properties of the regularizer (WCS) as it is the robustness measure that the DRO solution is attempting to control \eqref{eq:DRO-mean-sensitivity}.

We begin with the  definition of a {\it Generalized Measure of Deviation} \cite{rockafellar2006generalized}. This  generalizes the notion of standard deviation and measures the spread of a random variable \cite{rockafellar2006generalized}. 
\begin{definition}\label{def:deviation}
Let $f$ be a  random variable.
 ${\mathcal H}[f]$ is a {\it Generalized Measure of Deviation} or a {\it Generalized Measure of Spread} of $f$ if
\begin{enumerate}
\item ${\mathcal H}[f] \geq 0$ with equality if and only if $f$ is constant;
\item ${\mathcal H}[\beta f] = \beta{\mathcal H}[f]$ for every constant $\beta\geq 0$;
\item ${\mathcal H}[\alpha + f]={\mathcal H}[f]$ for every constant $\alpha\in{\mathbb R}$.
\end{enumerate}
\end{definition}

It is well known that for some uncertainty sets ${\mathcal Q}(\varepsilon)$ the ambiguity cost 
\begin{eqnarray}
\label{eq:AC}
{\mathcal A}(\varepsilon; f)  =  \max_{{\mathbb Q}\in{\mathcal Q}(\varepsilon)}{\mathbb E}_{\mathbb Q} \big[f\big] - {\mathbb E}_{\mathbb P}\big[f\big]
\end{eqnarray}
is a generalized measure of deviation and a risk-measure (e.g., Theorem 1 of \cite{rockafellar2006generalized}).  The spread of the cost distribution is related to robustness because its expected value  is sensitive to small changes in its tail probabilities if it has a large spread, and hence is not robust.

The following result shows that worst-case sensitivity ${\mathcal S}_{\mathbb P}[f(x, \cdot)]$ is a {\it generalized measure of deviation}, and hence, measures the spread of the cost distribution under the nominal.  
 The proof is in the Appendix \ref{sec:proof_WCS_is_GMD}.
\begin{proposition}
\label{prop:deviation}
Let ${f}:\Omega\rightarrow{\mathbb R}$ be a random variable on $(\Omega, {\mathcal F})$ and let ${\mathbb P}$ a probability measure on this space. Suppose that $d({\mathbb Q}\,|\,{\mathbb P})$ is convex and continuous over the set of probability measures $\mathbb Q$ which are absolutely continuous with respect to $\mathbb P$ and that $d({\mathbb Q}|{\mathbb P})=0$ if and only if ${\mathbb Q}={\mathbb P}$. 
Suppose too that there is a constant $k>0$ such that for every probability measure ${\mathbb Q}^{(\delta)}$ such that
\begin{eqnarray*}
\frac{\rm{d} {\mathbb Q}}{\rm{d} {\mathbb P}}^{(\delta)} = 1 + \delta {\Delta} > 0
\end{eqnarray*}
for some random variable $\Delta: \Omega \rightarrow {\mathbb R}$ satisfying ${\mathbb E}_{\mathbb P}(\Delta)=0$, that
\begin{align}
d({\mathbb Q}^{(\delta)}\,|\,{\mathbb P}) & \sim O(\delta^k) ~\mbox{when}  ~\delta \rightarrow 0.
\label{eq:continuity d}
\end{align}
Let ${\mathcal A}(\varepsilon; f) $ be the ambiguity cost \eqref{eq:AC} for an uncertainty set
\begin{eqnarray*}
{\mathcal Q}(\varepsilon)    =  \Big\{{\mathbb Q}:{\mathcal F} \rightarrow {\mathbb R} ~\Big|~d({\mathbb Q}\,| \,{\mathbb P})\leq\varepsilon\Big\}.
\end{eqnarray*}
 For every $\varepsilon>0$, ${\mathcal A}(\varepsilon;f)$ is a  generalized measure of deviation of $f$,
 ${\mathcal A}(\varepsilon;f)\sim O(\varepsilon^\frac{1}{k})$, and  worst-case sensitivity ${\mathcal S}_{\mathbb P}[f(x, \cdot)]$ defined by \eqref{eq:sensitivity-general2}  is a generalized measure of deviation with $g(\varepsilon)=\varepsilon^\frac{1}{k}$.
\end{proposition}
When the uncertainty set is a constraint on smooth $\phi$-divergence, $g(\varepsilon)=\sqrt{\varepsilon}$ in the definition \eqref{eq:sensitivity-general2} of worst-case sensitivity and \eqref{eq:continuity d} holds with $k=2$. It is linear in $\varepsilon$ in all other cases considered in this paper. More generally, $g(\varepsilon)$ is determined by the continuity property \eqref{eq:continuity d} of the divergence.

Proposition \ref{prop:deviation} tells us that the expected cost is sensitive to misspecification if its distribution has a large spread. 
However, different uncertainty sets induce different worst-case perturbations and have different measures of spread; when we select an uncertainty set, we are choosing a robustness measure. In this regard, Proposition \ref{prop:deviation} is the starting point for selecting uncertainty sets according to the vulnerabilities one wants to control.

\section{Worst-case sensitivity: Explicit formulas}
\label{sec:WCS}

\begin{figure}[h]
\centering
\begin{tabular}{cc}
\includegraphics[width=.45\linewidth]{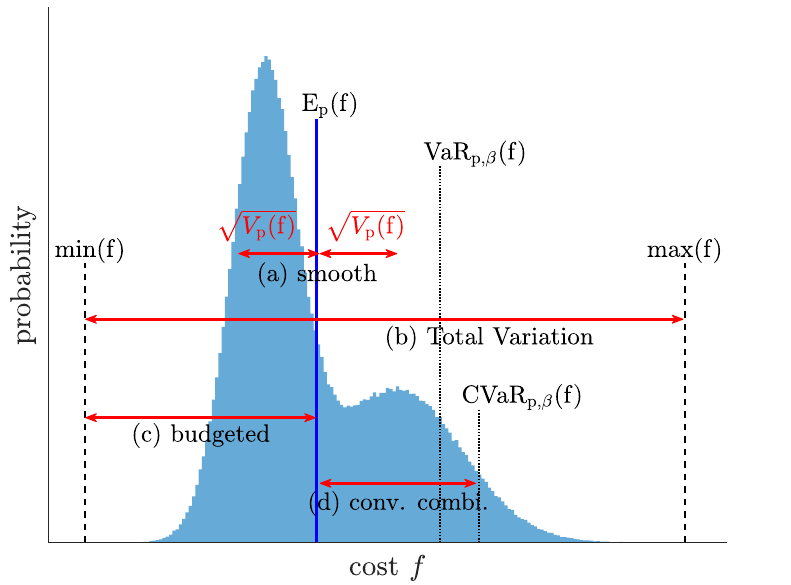}&
\includegraphics[width=.45\linewidth]{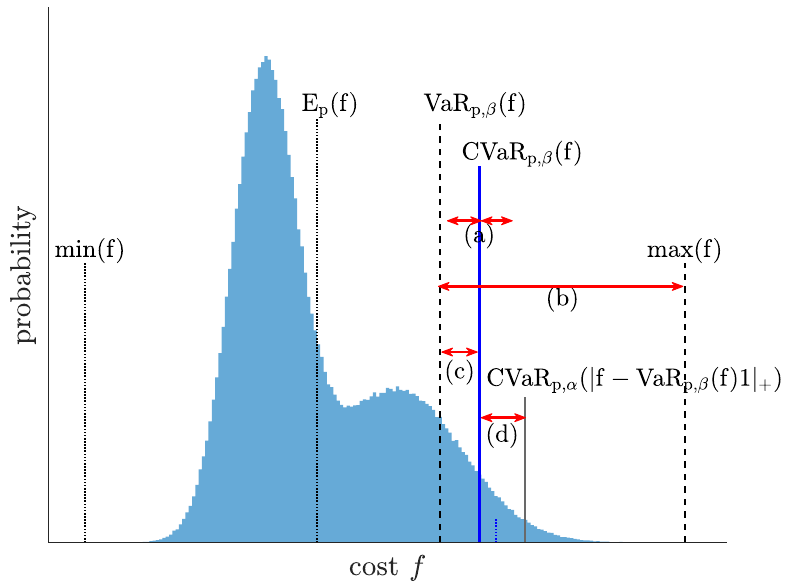}\\
$\langle$i$\rangle$~ ~WCSs for the mean & $\langle$ii$\rangle$~ ~WCSs for $\beta$-CVaR
\end{tabular}
\caption{Different worst-case sensitivities represent stress tests against spreads in different parts}
\label{fig:sensitivity}
{\raggedright \footnotesize See the Appendix \ref{sec:comparingWCS} for an analytical comparison of the deviation measures (a), (b), (c), and (d).\par}
\end{figure}

The regularized--nominal approximation \eqref{eq:DRO-mean-sensitivity} shows that DRO makes a tradeoff between performance and robustness and that DRO solutions depend critically on the sensitivity measure induced by the uncertainty set.
In this section we characterize the sensitivity measure induced by a selection of commonly used uncertainty sets (smooth $\phi$-divergence, total variation, budgeted uncertainty sets, uncertainty sets corresponding to a convex combination of the nominal distribution and a CVaR-type uncertainty set, and the Wasserstein metric). The resulting sensitivity measures can be very different and have a substantial impact on the DRO solution, highlighting the importance of sensitivity when selecting uncertainty sets. Table \ref{table:summary} and Figure \ref{fig:sensitivity} show the sensitivity measures induced by different uncertainty sets.

\begin{table}
\begin{center}
\caption{
Worst-case sensitivities considered in Section \ref{sec:WCS}}
\label{table:summary}
\begin{tabular}{ll|c}
\multicolumn{3}{c}{$\langle$i$\rangle$~ $\max_{\mathsf{q}}\{\mathbb{E}_{\mathsf{q}}(\mathsf{f})\vert d(\mathsf{q}\vert\mathsf{p})\leq \varepsilon\}$}\\
\hline
\multicolumn{2}{l|}{Type of divergence/worst-case objective} & Worst-case sensitivity \eqref{eq:sensitivity-general2} \\
\hline
\multicolumn{2}{l|}{(a) smooth $\phi$-divergence 
\qquad
$\phi''(1)>0=\phi(1)=\phi'(1)$} & $\sqrt{\frac{2{\mathbb V}_{\mathsf{p}}(\mathsf{f})}{\phi''(1)}}$ \\ [2pt]
(b) total variation &
$\phi(z)=|z-1|$ & $\frac{1}{2}\big(\max(\mathsf{f})-\min(\mathsf{f})\big)$\\ [2pt]
(c) ``budgeted'' 
& $\phi=\delta_{[0,1+\varepsilon]}$ & $\mathbb{E}_{\mathsf{p}}(\mathsf{f})-\min(\mathsf{f})$ \\[2pt]
(d) $(1-\varepsilon)$``mean
''$+\varepsilon$``$\alpha$-CVaR'' & 
$\phi=\delta_{[1-\varepsilon,\frac{1-(1-\varepsilon)\alpha}{1-\alpha}]}$
& $\mathrm{CVaR}_{\mathsf{p},
\alpha}(\mathsf{f})-\mathbb{E}_{\mathsf{p}}
(\mathsf{f})$\\ [5pt]
(d') $(1-\varepsilon)$``mean
''$+\varepsilon$``max
'' 
& $\phi=\delta_{[1-\varepsilon,
\infty)}$
 & $\max(\mathsf{f})-\mathbb{E}_{\mathsf{p}}(\mathsf{f})$\\[2pt]
(d'') symmetric ($\alpha=\frac{1}{2-\varepsilon}$) & 
$\phi=\delta_{[1-\varepsilon,\frac{1}{1-\varepsilon}]}$ & $\mathrm{CVaR}_{\mathsf{p},\frac{1}{2}}(\mathsf{f})-\mathbb{E}_{\mathsf{p}}(\mathsf{f})$ \\[5pt]
\multicolumn{2}{l|}{(e) type-1 Wasserstein distance
~w/ Lipschitz smooth $f$ 
} & $\max\limits_{i=1,\cdots,n} \max\limits_{z_i}\frac{f(z_i)-f(Y_i)}{\|z_i-Y_i\|}$ \\
\hline
\multicolumn{3}{c}{}\\
\multicolumn{3}{c}{$\langle$ii$\rangle$~ $\max_{\mathsf{q}}\{ \mathrm{CVaR}_{\mathsf{q},\beta}(\mathsf{f})\vert d(\mathsf{q}\vert\mathsf{p})\leq \varepsilon\}$}\\
\hline
\multicolumn{2}{l|}{Type of divergence/worst-case objective} & Worst-case sensitivity \eqref{eq:wsc_for_cvar} \\
\hline
\multicolumn{2}{l|}{(a) smooth $\phi$-divergence 
\qquad
$\phi''(1)>0=\phi(1)=\phi'(1)$} & $
\frac{1}{1-\beta}\sqrt{\frac{2}{\phi''(1)}\mathbb{V}_{\mathsf{p}}\big(|\mathsf{f}-\mathrm{VaR}_{\mathsf{p},\beta}(\mathsf{f})\mathsf{1}|_+\big)}$ \\
(b) total variation & $\phi(z)=|z-1|$ & $\frac{1}{2(1-\beta)}\big(\max(\mathsf{f})-\mathrm{VaR}_{\mathsf{p},\beta}(\mathsf{f})\big)$\\
(c) ``budgeted" & $\phi=\delta_{[0,1+\varepsilon]}$ & $\mathrm{CVaR}_{\mathsf{p},\beta}(\mathsf{f})-\mathrm{VaR}_{\mathsf{p},\beta}(\mathsf{f})$\\
(d) $(1-\varepsilon)$``mean
''$+\varepsilon$``$\alpha$-CVaR'' & 
$\phi=\delta_{[1-\varepsilon,\frac{1-(1-\varepsilon)\alpha}{1-\alpha}]}$ & $\frac{1}{1-\beta}
\mathrm{CVaR}^{\Delta}_{\mathsf{p},\alpha}(\left|\mathsf{f}-\mathrm{VaR}_{\mathsf{p},\beta}(\mathsf{f})\mathsf{1}\right|_+)$\\ 
\hline
\end{tabular}\\
\footnotesize $\mathrm{CVaR}^{\Delta}_{\mathsf{p},\alpha}\equiv
\mathrm{CVaR}_{\mathsf{p},\alpha}-\mathbb{E}_{\mathsf{p}}$ denotes the operator of $\alpha$-CVaR Deviation.
\end{center}
\end{table}

We focus in this section on discrete random variables and treat random variables as an $n$-vector.
Let $\mathsf{f}:=(f_1,...,f_n)^\top$ 
and $\mathsf{p}:=(p_1,...,p_n)^\top$, where 
$p_i$ is the probability mass on $f_i$.
Without loss of generality, we assume that $p_i>0$ for all $i=1,...,n$. 
Likewise, we reserve $\mathsf{q}:=(q_1,...,q_n)^\top$ for an alternative probability distribution $\mathbb{Q}$, and $\mathbb{E}_{\mathsf{p}}(\mathsf{f}):=\mathsf{p}^\top\mathsf{f}$, $\mathbb{V}_{\mathsf{p}}(\mathsf{f})$, and $\mathcal{S}_{\mathsf{p}}(\mathsf{f})$ denote, respectively, expectation $\mathbb{E}_{\mathbb{P}}[f]$, variance $\mathbb{V}_{\mathbb{P}}[f]$, and worst-case sensitivity $\mathcal{S}_{\mathbb{P}}[f]$ of $\mathsf{f}\equiv f$ under $\mathsf{p}\equiv \mathbb{P}$. In addition, $\mathsf{1}$ and $\mathsf{0}$ denote column vectors with all entries being $1$'s and $0$'s, respectively. 
Many of our results have generalizations to more complex settings. The restriction to discrete random variables enables us to communicate our message with minimal technical fuss.  For the purposes of readability, all proofs can be found in Appendix \ref{App:WCSproofs}.

{

\subsection{Smooth $\phi$-divergences}
\label{sec:phi-div}
Consider the worst-case objective
\begin{align}
V_{\phi}(\varepsilon;\mathsf{f}) &:= \max_{\mathsf{q}\in{\mathcal Q}_{\phi}(\varepsilon)}\sum_{i=1}^{n}q_i 
f_i,
\label{eq:phi-theta}
\end{align}
where
\begin{align}
{\mathcal Q}_{\phi}(\varepsilon) &:= \Big\{\mathsf{q}:=(q_1,...,q_n)^\top\in{\mathbb R}^n \, \Big| \, \sum_{i=1}^{n}p_i \phi\Big(\frac{q_i}{p_i}\Big) \leq \varepsilon,~\mathsf{1}^\top\mathsf{q}=1,~\mathsf{q}\geq\mathsf{0}\Big\}.
\label{eq:phi-div-US}
\end{align}
We assume the following.
\begin{assumption}
\label{ass:phi}
$\phi(z)$ is strictly convex, twice continuously differentiable in $z$, with $\phi(1)=0$, $\phi'(1)=0$ and $\phi''(1)>0$.
\end{assumption}

Proposition \ref{prop:phi-div} in the Appendix shows that the  worst-case distribution (the maximizer in \eqref{eq:phi-theta}) is
\begin{align*}
q_i(\varepsilon) & = p_i\Big\{1 + \sqrt\frac{2\varepsilon}{{{\mathbb V}_{\mathsf{p}}(\mathsf{f})}} \left(f_i-\mathbb{E}_{\mathsf{p}}(\mathsf{f})\right)\Big\}+ o(\sqrt\varepsilon).
\end{align*}
It follows that the worst-case expected cost is 
\begin{align*}
V_{\phi}(\varepsilon) & = \mathbb{E}_{\mathsf{q}(\varepsilon)}(\mathsf{f})\nonumber \\
& = 
\mathbb{E}_{\mathsf{p}}(\mathsf{f}) +\sqrt{\varepsilon}\sqrt{\frac{2\mathbb{V}_{\mathsf{p}}(\mathsf{f})}{\phi''(1)}} + o(\sqrt{\varepsilon})
\end{align*}
and that $V_{\phi}(\varepsilon)-V_{\phi}(0)$ is $O(\sqrt \varepsilon)$. The following result shows that WCS for smooth $\phi$-divergence is the standard deviation of the cost distribution.
\begin{proposition}
Suppose Assumption \ref{ass:phi} is satisfied. Then
\begin{align}
{\mathcal S}_{\mathsf{p}}(\mathsf{f})
 &=\lim_{\varepsilon\downarrow 0}\frac{V_{\phi}(\varepsilon)-V_{\phi}(0)}{\sqrt{\varepsilon}}
  =\sqrt{\frac{2{\mathbb V}_{\mathsf{p}}(\mathsf{f})}{\phi''(1)}}.
\label{eq:phi-sensitivity}
\end{align}
\end{proposition}
A closely related result for the case $\phi(z)$ is relative entropy was derived in \cite{lam2016robust}, while \cite{gotoh2018robust} derives worst-case sensitivity for the ``penalty formulation" of the DRO model. Extensions to the penalty case can be found in Appendix \ref{App:Penalty}.

\begin{example}
When $\phi$-divergence is  modified $\chi^2$-deviation, $\phi(z) = \frac{1}{2}(z-1)^2$, the expressions in Proposition \ref{prop:phi-div} are exact at $O(\sqrt{\varepsilon})$ and ${\mathcal S}_{\mathsf{p}}(\mathsf{f}) = \sqrt{2 {\mathbb V}_{\mathsf{p}}(\mathsf{f})}$.
\end{example}

\subsection{Nonsmooth $\phi$-divergences}\label{sec:nonsmooth_phi}
Next consider uncertainty sets \eqref{eq:phi-div-US} defined by $\phi$-divergences 
with $\phi(z)=|z-1|$ or $\phi(z)=\delta_{I}(z)$ for a closed interval $I\subset[0,\infty)
$ such that $1\in I$, where $\delta_{I}(z):=0$ for $z\in I;~\infty$,~otherwise. Neither satisfies Assumption \ref{ass:phi} since the former is not differentiable at $z=1$ and $\phi''(1)=0$ for the latter. 

\subsubsection{Total variation}\label{sec:TV}
For $\phi(z)=|z-1|$ and $\varepsilon\geq 0$, the uncertainty set is given as  
\[
\mathcal{Q}_\mathrm{TV}(\varepsilon):=\Big\{\mathsf{q}\in\mathbb{R}^{n}\,\Big|\,
\mathsf{1}^\top|\mathsf{q}-\mathsf{p}|\leq\varepsilon,~
\mathsf{1}^\top\mathsf{q}=1,~\mathsf{q}\geq\mathsf{0}\Big\},
\]
where $|\mathsf{z}|:=(|z_1|,...,|z_n|)^\top$, and denote the worst-case expected cost
\begin{align}
V_{\rm TV}(\varepsilon;\mathsf{f}) := \max_{\mathsf{q}\in\mathcal{Q}_\mathrm{TV}(\varepsilon)}~\mathbb{E}_{\mathsf{q}}(\mathsf{f}).
\label{def:wcobj:TV}
\end{align}
We can focus on $\varepsilon\leq 2$ since 
$\mathcal{Q}_\mathrm{TV}(\varepsilon)$ coincides with the unit simplex, i.e., 
entire space of 
all distributions,  
 for $\varepsilon>2$. 
Consider an ordering of the components of the cost vector $\mathsf{f}=(f_1,...,f_n)^\top\in\mathbb{R}^n$ from largest to smallest and denote the  $i^{th}$ largest component by $f_{(i)}$, i.e., $f_{(1)}\geq\cdots\geq f_{(n)}$,
and let $p_{(i)}$ denote the probability mass corresponding to $f_{(i)}$, whereas $q_{(i)}$ denotes the worst-case probability mass on $f_{(i)}$. For $\varepsilon \in (0, \min(\mathsf{p}))$, 
it can be shown (Proposition \ref{lemma:sen_tv1} in the Appendix) that a worst-case probability distribution for \eqref{def:wcobj:TV} is
\begin{equation*}
(q_{(1)},q_{(2)},...,q_{(n-1)},q_{(n)})=
\big(p_{(1)}+\frac{\varepsilon}{2},~p_{(2)},~...,~p_{(n-1)},~p_{(n)}-\frac{\varepsilon}{2}\big)
\end{equation*}
and the worst-case objective is
\[
V_{\rm TV}(\varepsilon;\mathsf{f})=
\mathbb{E}_\mathsf{p}(\mathsf{f})+\varepsilon \cdot \frac{\max(\mathsf{f})-\min(\mathsf{f})}{2}
\]
(this is also derived in  \cite{charalambous2013TV}). Accordingly,  $g(\varepsilon)=\varepsilon$, and the worst-case sensitivity of \eqref{def:wcobj:TV} is as follows.
\begin{align}
{\mathcal S}_{\mathsf{p}}(\mathsf{f})
=\lim_{\varepsilon\downarrow 0}\frac{V_{\rm TV}(\varepsilon)-V_{\rm TV}(0)}{\varepsilon}
= \frac{\max(\mathsf{f})-\min(\mathsf{f})}{2}
\equiv\frac{1}{2}\times\mbox{``Range of }\mathsf{f}.\mbox{''}
\label{eq:sen_tv}
\end{align}

\subsubsection{Budgeted uncertainty}
\label{sec:CVaR}
With $\phi(z)=\delta_{[0,1+\varepsilon]}(z)$ for $\varepsilon\geq 0$, the uncertainty set 
becomes
\begin{align}
{\mathcal Q}_{\mathrm{b}}(\varepsilon) & = \Big\{\mathsf{q}\in\mathbb{R}^n\,\Big|\,\mathsf{1}^\top\mathsf{q}=1,\mathsf{0}\leq\mathsf{q}\leq(1+\varepsilon)\mathsf{p}\Big\}
\label{eq:CVaR_uncertainty}
\end{align}
and define the associated worst-case expected cost as
\begin{equation}
V_{\rm b}(\varepsilon;\mathsf{f}):=\max_{\mathsf{q}\in{\mathcal Q}_{\rm b}(\varepsilon)} \mathbb{E}_{\mathsf{q}}(\mathsf{f}).
\label{eq:dual_cvar_eps}
\end{equation}
Clearly, ${\mathcal Q}_{\mathrm{b}}(\varepsilon)=\{\mathsf{p}\}$ if and only if $\varepsilon=0$.  
We refer to the set \eqref{eq:CVaR_uncertainty} as the  \emph{budgeted uncertainty set} since the worst-case distribution for \eqref{eq:dual_cvar_eps} is obtained, as below, by applying a greedy algorithm for (a scaled version of) continuous knapsack problem (see the Appendix \ref{sec:proof:wcs_budgeted} for the proof).
\begin{proposition}\label{prop:cvar_explicit}
Let $\varepsilon>0$. 
The worst-case distribution for \eqref{eq:dual_cvar_eps} is given by 
\[
(q_{(1)},..., q_{(k)},q_{(k+1)},q_{(k+2)},...,q_{(n)})=
\big((1+\varepsilon)p_{(1)},...,(1+\varepsilon)p_{(k)},1-(1+\varepsilon)\sum_{i=1}^kp_{(i)},0,...,0\big)
\]
and the worst-case expected cost 
is
\[
V_{\rm b}(\varepsilon;\mathsf{f})=
(1+\varepsilon)\sum_{i=1}^{k}p_{(i)}f_{(i)}+\big\{1-(1+\varepsilon)\sum_{i=1}^{k}p_{(i)}\big\}f_{(k+1)},
\]
where $k\in\{0,...,n-1\}$ is such that $\sum\limits_{i=1}^{k}p_{(i)}\leq \frac{1}{1+\varepsilon}<\sum\limits_{i=1}^{k+1}p_{(i)}$, where we denote $\sum\limits_{i=1}^0p_{(i)}=0$. 
\end{proposition}
The worst-case expected cost \eqref{eq:dual_cvar_eps} 
is related to {\it Conditional Value-at-Risk} (CVaR, for short) \cite{rockafellar2002conditional}. For 
$\alpha\in[0, 1)$, it is referred to as $\alpha$-CVaR and defined by
\begin{align}
\mathrm{CVaR}_{\mathsf{p},\alpha}(\mathsf{f})&:=\min_\gamma\Big\{\gamma+\frac{1}{1-\alpha}\mathbb{E}_{\mathsf{p}}(\max\{\mathsf{f}-\gamma\mathsf{1},\mathsf{0}\})\Big\}
\label{eq:primal_cvar}
\\
&=\max_{\mathsf{q}}\Big\{
\mathsf{f}^\top\mathsf{q}
\,\Big|\,\mathsf{1}^\top\mathsf{q}=1,\mathsf{0}\leq\mathsf{q}\leq\frac{1}{1-\alpha}\mathsf{p}\Big\},
\label{eq:dual_cvar}
\end{align}
where the last equality follows from the duality. 
Obviously, $V_{\rm b}(\varepsilon)=\mathrm{CVaR}_{\mathsf{p},\frac{\varepsilon}{1+\varepsilon}}(\mathsf{f})$.
\begin{proposition}
\label{cor:sen_cvar} 
For  uncertainty set \eqref{eq:CVaR_uncertainty} and $g(\varepsilon)=\varepsilon$, worst-case sensitivity of \eqref{eq:dual_cvar_eps} is given by
\begin{align}
{\mathcal S}_{\mathsf{p}}(\mathsf{f}) 
=\lim_{\varepsilon\downarrow 0}\frac{V_{\rm b}(\varepsilon)-V_{\rm b}(0)}{\varepsilon}
= \mathbb{E}_{\mathsf{p}}(\mathsf{f})-\min(\mathsf{f}).
\label{eq:sen_cvar}
\end{align}
\end{proposition}
While worst-case sensitivity \eqref{eq:sen_cvar} is a measure of spread, it only depends on the  ``good" side of the cost distribution. In contrast, those for smooth $\phi$-divergence \eqref{eq:phi-sensitivity} and total variation \eqref{eq:sen_tv} depend on the entire distribution.

\subsubsection{Convex combination of expected cost and CVaR}
\label{sec:ExpLoss_CVaR}

Let $\alpha\in[0,1)$ be a fixed parameter and
\begin{align*}
\mathcal{Q}_{\rm CVaR}(\alpha):=
\Big\{\mathsf{q}\in\mathbb{R}^{n}\,\Big|\,
\mathsf{1}^\top\mathsf{q}=1,~\mathsf{0}\leq\mathsf{q}\leq\frac{1}{1-\alpha}\mathsf{p}\Big\}
\end{align*}
(the feasible set of \eqref{eq:dual_cvar}).
Consider the uncertainty set
\begin{align}
\mathcal{Q}_{\rm c}(\varepsilon)&:=
(1-\varepsilon)\{\mathsf{p}\}+\varepsilon\mathcal{Q}_{\rm CVaR}(\alpha)
\label{eq:CVaR_uncertainty_set}
\end{align}
parameterized by $\varepsilon\in[0, 1]$, and the associated worst-case expected cost 
\begin{align}
V_{\rm c}(\varepsilon;\mathsf{f})&:=\max_{\mathsf{q}\in\mathcal{Q}_\mathrm{c}(\varepsilon)}~\mathbb{E}_{\mathsf{q}}(\mathsf{f}).
\label{eq:V_c}
\end{align}
Observe that $\mathcal{Q}_{\rm c}(0)=\{\mathsf{p}\}$, so there is no robustness if $\varepsilon=0$. 
The uncertainty set \eqref{eq:CVaR_uncertainty_set}  was considered in \cite{anderson2022improving} (see also \cite{Tsang2025TRO}) and is equivalent to
\begin{align*}
\mathcal{Q}_{\rm c}(\varepsilon) &= \Big\{\mathsf{q}\in{\mathbb R}^n\,\Big|\, \mathsf{1}^\top\mathsf{q}=1,~\mathsf{p}(1-\varepsilon)\leq \mathsf{q}\leq \mathsf{p}(1-\varepsilon)+\frac{\varepsilon}{1-\alpha}\mathsf{p}\Big\},
\end{align*}
which is equivalent to the case with $\phi(z)=\delta_{[1-\varepsilon,\frac{1-(1-\varepsilon)\alpha}{1-\alpha}]}(z)$.
\begin{proposition}
\label{cor:wcs:comb}
The worst-case probability distribution for \eqref{eq:V_c} is given by
\begin{align*}
\lefteqn{\big(q_{(1)},~...,~q_{(k)},~q_{(k+1)},~q_{(k+2)},~...,~q_{(n)}\big) }
 \\
&=\Big((1+\frac{\alpha\varepsilon}{1-\alpha})p_{(1)},...,(1+\frac{\alpha\varepsilon}{1-\alpha})p_{(k)},(1-\varepsilon)p_{(k+1)}+\varepsilon\big(1-\frac{1}{1-\alpha}\sum_{i=1}^kp_{(i)}\big),\\
&\qquad\qquad\qquad\qquad(1-\varepsilon)p_{(k+2)},...,(1-\varepsilon)p_{(n)}\Big),
\end{align*}
where $k$ is such that $\sum\limits_{i=1}^kp_{(i)}\leq 1-\alpha<\sum\limits_{i=1}^{k+1}p_{(i)}$,  
and the worst-case 
expected cost \eqref{eq:V_c} is
\begin{equation}
V_{{\rm c}}(\varepsilon;\mathsf{f})
=(1-\varepsilon)\mathbb{E}_\mathsf{p}(\mathsf{f})+\varepsilon\mathrm{CVaR}_{\mathsf{p},\alpha}(\mathsf{f})
=\mathbb{E}_\mathsf{p}(\mathsf{f})+\varepsilon\big(\mathrm{CVaR}_{\mathsf{p},\alpha}(\mathsf{f})-\mathbb{E}_\mathsf{p}(\mathsf{f})\big),
\label{eq:sen_cvxcmb_expectation+cvar}
\end{equation}
for $\varepsilon\in[0,1]$. Accordingly, for uncertainty set \eqref{eq:CVaR_uncertainty_set} and $g(\varepsilon)=\varepsilon$, worst-case sensitivity is 
\begin{align}
{\mathcal S}_{\mathsf{p}}(\mathsf{f}) 
=\lim_{\varepsilon\downarrow 0}\frac{V_{\rm c}(\varepsilon)-V_{\rm c}(0)}{\varepsilon}
= \mathrm{CVaR}_{\mathsf{p},\alpha}(\mathsf{f})-\mathbb{E}_\mathsf{p}(\mathsf{f})
~\equiv~``\alpha\mbox{-CVaR Deviation of }\mathsf{f}.\mbox{''}
\label{eq:sen-CVaRdev}
\end{align}
\end{proposition}
See the Appendix \ref{sec:proof:convcomb} for the proof.   

Expression \eqref{eq:sen_cvxcmb_expectation+cvar} suggests that $V_{\rm c}(\varepsilon)$ is linearly increasing at a rate of its \emph{CVaR Deviation} \cite{rockafellar2006generalized}  for any non-uniform vector $\mathsf{f}$ and $\varepsilon\in(0,1]$. 
The worst-case sensitivity \eqref{eq:sen-CVaRdev} is the spread of the  ``bad" part of the cost distribution, which contrasts with \eqref{eq:sen_cvar}. 
In particular, for $\alpha\in[1-p_{(1)},1)$, $\mathrm{CVaR}_{\mathsf{p},\alpha}(\mathsf{f})=\max(\mathsf{f}):=\max\{f_1,...,f_n\}$ and $\mathcal{S}_{\mathsf{p}}(\mathsf{f})=\max(\mathsf{f})-\mathbb{E}_\mathsf{p}(\mathsf{f})$. 

From \eqref{eq:sen_cvxcmb_expectation+cvar} and \eqref{eq:sen-CVaRdev}, we see that the WCS is equal to the average sensitivity, or that 
$V_{\rm c}(\varepsilon;\mathsf{f})=\mathbb{E}_{\mathsf{p}}(\mathsf{f})+\varepsilon \mathcal{S}_{\mathsf{p}}(\mathsf{f}),
$ 
so, unlike the case depicted in Figure \ref{fig:illus2}, the frontier from the DRO solution coincides with the mean--sensitivity frontier. 
This feature differentiates $V_{\rm c}$ from the other $\phi$-divergence-based objectives, $V_{\phi},V_{\rm TV}$, and $V_{\rm b}$. 
Also, we can associate the frontier of $V_{\rm c}$ with 
risk-averse decision making. Specifically, 
any point on the mean--CVaR Deviation frontier for $\varepsilon\in[0,1]$ is 
consistent with risk-averse expcted utilities 
if the point $(\mathbb{E}_{\mathsf{p}}(\mathsf{f}),\mathcal{S}_{\mathsf{p}}(\mathsf{f}))$ corresponds to a unique $\mathsf{f}$ (see \cite{ogryczak2002dual} for the details).

\begin{remark}
More generally, for any $L,U$ such that $0 \leq L\leq 1 \leq U$, the uncertainty set
\begin{align*}
\mathcal{Q}_{\rm h}(L,U) & := \Big\{\mathsf{q}\in{\mathbb R}^n\,\Big|\,\mathsf{1}^\top\mathsf{q}=1,~L\mathsf{p}\leq\mathsf{q}\leq U\mathsf{p}
\Big\}
\end{align*}
corresponds to $\phi(z)=\delta_{[L,U]}(z)$ and the worst-case expected cost 
\begin{align*}
V_{\rm h}(L,U;\mathsf{f})&:=
L\cdot\mathbb{E}_\mathsf{p}(\mathsf{f})+(1-L)\cdot\mathrm{CVaR}_{\mathsf{p},\frac{U-1}{U-L}}(\mathsf{f}).
\end{align*}
If $(L,U)=(0,1+\varepsilon)$, $V_{\rm h}(L,U)$ is $V_{\rm b}(\varepsilon)$. 
If $(L,U)=(1-\varepsilon,\frac{1-(1-\varepsilon)\alpha}{1-\alpha})$, $V_{\rm h}(L,U)$ is $V_{\rm c}(\varepsilon)$. 
In particular for the latter, 
while the parameter $\alpha$ for CVaR is usually fixed (e.g., at $0.95$ or $0.99$), it can be viewed as another hyperparameter, in addition to $\varepsilon$, that defines the uncertainty set. A reasonable option to merge the two parameters $L,U$ into a single one is to set as $U=\frac{1}{L}=\frac{1}{1-\varepsilon}$ for $\varepsilon\in[0,1]$, with which the uncertainty set becomes symmetric between $\mathsf{q}$ and $\mathsf{p}$ since $\delta_{[L,U]}(\frac{q_i}{p_i})=\delta_{[L,U]}(\frac{p_i}{q_i})$ or 
\begin{align*}
\mathcal{Q}_{\rm s}(\nu)&:=\Big\{\mathsf{q}
\,\Big|\, \mathsf{1}^\top\mathsf{q}=1,
(1-\varepsilon)\mathsf{p}\leq\mathsf{q}\leq
\frac{1}{1-\varepsilon}
\mathsf{p}
\Big\}=\Big\{\mathsf{q}
\,\Big|\, \mathsf{1}^\top\mathsf{q}=1,
(1-\varepsilon)
\mathsf{q}\leq\mathsf{p}\leq
\frac{1}{1-\varepsilon}
\mathsf{q}
\Big\}.
\end{align*}
For the symmetric case $\mathcal{Q}_{\rm s}(\varepsilon)$, the worst-case expected cost is then equivalent to the convex combination of the mean cost and $\frac{1}{2-\varepsilon}$-CVaR, and the 
 worst-case sensitivity becomes 
$\mathcal{S}_{\mathsf{p}}(\mathsf{f})=\mathrm{CVaR}_{\mathsf{p},\frac{1}{2}}(\mathsf{f})-\mathbb{E}_\mathsf{p}(\mathsf{f})$, i.e., $\frac{1}{2}$-CVaR Deviation of $\mathsf{f}$.
\end{remark}


\subsection{Extensions to Robust CVaR}
Consider distributionally robust $\beta$-CVaR 
\begin{align*}
\mathrm{RCVaR}_{\mathsf{p},\beta}^{\varepsilon}(\mathsf{f})
:=\max_{\mathsf{q}\in\mathcal{Q}(\varepsilon)} \mathrm{CVaR}_{\mathsf{q},\beta}(\mathsf{f})
\end{align*}
where 
$\mathcal{Q}(\varepsilon)$ denotes an uncertain set of probability distribution $\mathsf{q}$. 
The worst-case sensitivity can also be defined with a growth rate $g(\cdot)$ as
\begin{align}
\mathcal{S}_{\mathsf{p}}(\mathsf{f})=\lim_{\varepsilon\downarrow 0}\frac{\mathrm{RCVaR}_{\mathsf{p},\beta}^{\varepsilon}(\mathsf{f})
-\mathrm{CVaR}_{\mathsf{p},\beta}(\mathsf{f})}{g(\varepsilon)}.
\label{eq:wsc_for_cvar}
\end{align}
To ease complexity, we assume the following nondegeneracy condition. 
\begin{assumption}\label{ass:nondegeneracy_for_cvar}
There is an integer 
$k$ such that $\sum\limits_{i=1}^{k}p_{(i)}<1-\beta<\sum\limits_{i=1}^{k+1}p_{(i)}$. 
\end{assumption}
Note that Assumption \ref{ass:nondegeneracy_for_cvar} implies that $\beta$-VaR of $\mathsf{f}$ is given as $\mathrm{VaR}_{\mathsf{p},\beta}(\mathsf{f})=f_{(k+1)}$. 
\begin{proposition}
\label{prop:rcvar_sensitivities}
Suppose that distribution of $\mathsf{f}$ satisfies Assumption \ref{ass:nondegeneracy_for_cvar}. 
\begin{list}{}{
\topsep 0em
\parskip 0em
\partopsep 0em
\leftmargin 0em
\itemindent 0em
}
\item[(a)] With 
$\mathcal{Q}(\varepsilon)=\mathcal{Q}_{\phi}(\varepsilon)$ and $g(\varepsilon)=\sqrt{\varepsilon}$, worst-case sensitivity \eqref{eq:wsc_for_cvar} is  
\begin{align*}
\mathcal{S}_{\mathsf{p}}(\mathsf{f})
&=
\frac{1}{1-\beta}\sqrt{\frac{2}{\phi''(1)}\mathbb{V}_{\mathsf{p}}\big(|\mathsf{f}-\mathrm{VaR}_{\mathsf{p},\beta}(\mathsf{f})\mathsf{1}|_+\big)},
\end{align*}
where $|\mathsf{f}|_+:=\max(\mathsf{f},\mathsf{0})\equiv(\max\{f_1,0\},...,\max\{f_n,0\})^\top$ for $\mathsf{f}\in\mathbb{R}^n$. 
%
\item[(b)] With 
$\mathcal{Q}(\varepsilon)=\mathcal{Q}_{\rm TV}(\varepsilon)$ and  $g(\varepsilon)=\varepsilon$, worst-case sensitivity \eqref{eq:wsc_for_cvar} is  
\begin{align*}
\mathcal{S}_{\mathsf{p}}(\mathsf{f})
&=\frac{1}{2(1-\beta)}\big(\max(\mathsf{f})-\mathrm{VaR}_{\mathsf{p},\beta}(\mathsf{f})\big)
.
\end{align*}
\item[(c)] 
With 
$\mathcal{Q}(\varepsilon)=\mathcal{Q}_{\rm b}(\varepsilon)$ and  $g(\varepsilon)=\varepsilon$, worst-case sensitivity \eqref{eq:wsc_for_cvar} is  
\begin{align*}
\mathcal{S}_{\mathsf{p}}(\mathsf{f})
&=\mathrm{CVaR}_{\mathsf{p},\beta}(\mathsf{f})-\mathrm{VaR}_{\mathsf{p},\beta}(\mathsf{f}).
\end{align*}
\item[(d)] With 
$\mathcal{Q}(\varepsilon)=\mathcal{Q}_{\rm c}(\varepsilon)$ and  $g(\varepsilon)=\varepsilon$, worst-case sensitivity \eqref{eq:wsc_for_cvar} is
\begin{align*}
\mathcal{S}_{\mathsf{p}}(\mathsf{f})
&=\frac{1}{1-\beta}
\mbox{``}\alpha\mbox{-CVaR Deviation of }\left|\mathsf{f}-\mathrm{VaR}_{\mathsf{p},\beta}(\mathsf{f})\mathsf{1}\right|_+.\mbox{''}
\end{align*}
\end{list}
\end{proposition}
See the Appendix \ref{sec:proof:wcs_cvar} for the proof. 

Recalling that CVaR is approximately equal to the conditional expected cost over the quantile (VaR), i.e., the mean of $f_{(1)}-\mathrm{VaR}_{\mathsf{p},\beta}(\mathsf{f}),...,f_{(k)}-\mathrm{VaR}_{\mathsf{p},\beta}(\mathsf{f})$, 
Proposition \ref{prop:rcvar_sensitivities} shows that worst-case sensitivity of the robust CVaR objective with different uncertainty set can be captured as different measures of spread 
 of the tail distribution above VaR 
(Figure \ref{fig:sensitivity}$\langle$ii$\rangle$).


\subsection{Wasserstein metric}
\label{sec:Wasserstein}
Consider the worst-case expected cost with a constraint on the type-$1$ Wasserstein metric \cite{BDOW,blanchet2019,esfahani2018data,gao2024wasserstein}:
\begin{align}
V_{\rm w}(\varepsilon) &:= \max_{\gamma \in {\mathcal X}} \Big\{ \int_z f(x,\,z) \Big(\sum_{i=1}^{n} \gamma_i({\rm d}z) \Big) \,\Big|\,
\sum_{i=1}^{n}\int_z \|z-Y_i\|\gamma_i({\rm d}z) \leq \varepsilon \Big\},
\label{eq:W2}
\end{align}
where $\|\cdot\|$ is any norm on ${\mathbb R}^n$ and 
\begin{align*}
{\mathcal X} = \Big\{\gamma\;\Big|\; 
\int_z\gamma_i({\rm d}z) = p_i,\; 
\gamma_i({\rm d}z) \geq 0,\; i=1,\cdots,n
\Big\}.
\end{align*}
When $\varepsilon$ is small it can be shown (see the Appendix \ref{sec:proof:wcs_wasserstein}) that
\begin{align*}
V_{\rm w}(\varepsilon) & = \sum_{i=1}^np_i f(Y_i) + \varepsilon \Big(\max_{i=1,\cdots,\,n} \max_{z_i}\frac{f(z_i)-f(Y_i)}{\|z_i-Y_i\|}\Big) + o(\varepsilon).
\end{align*}

The proof of the following result can be found in the Appendix \ref{sec:proof:wcs_wasserstein}.
\begin{proposition}
\label{prop:Wass_sensitivity}
 For \eqref{eq:W2}, we have
\begin{align}
{\mathcal S}_{\mathbb{P}}[f] &= 
\lim_{\varepsilon\downarrow 0}\frac{V_{\rm w}(\varepsilon)-V_{\rm w}(0)}{\varepsilon}
 =  \max_{i=1,\cdots,\,n} \max_{z_i}\frac{f(z_i)-f(Y_i)}{\|z_i-Y_i\|}.
\label{eq:Wass-sensitivity}
\end{align}
\end{proposition}

\medskip

The Wasserstein uncertainty set is quite different from the others in that sensitivity is measured by ``local perturbations" corresponding to moving probability mass from the support of the nominal distribution to a worst-case point. 


\begin{example}
\label{ex:Wass-inv}
Consider the cost function
\begin{align}
f(x, Y)= -r \min\{x,  Y\} + q\max\{x-Y, 0\}+ s\max\{Y-x, 0\} + c x
\label{eq:inv2}
\end{align}
where $0\leq q<c<r$ and $s \geq 0$. The negative of this cost function is the reward function for an inventory problem, so minimizing ${\mathbb E}_{\mathbb P}[f(x, Y)]$ is equivalent to maximizing expected reward. If $x \in (\min_i Y_i, \max_i Y_i)$, it can be shown that for a Wasserstein metric with $p=1$
\begin{align*}
{\mathcal S}_{\mathbb P} [f(x, \cdot)] = \max\{r-q, s\},
\end{align*}
so  the SAA optimizer is also the solution of the DRO problem for a large range of $\varepsilon$, beyond which, the order quantity is either smaller than $\min_i Y_i$ or larger than $\max_i Y_i$, which is not sensible. This suggests that the Wasserstein uncertainty set with $p=1$ may not be a good choice for the robust inventory problem.
\end{example}

\subsection{Summary}
Worst-case sensitivity is a measure of robustness. Each uncertainty set defines a robustness measure. Table \ref{table:summary} and  Figure \ref{fig:sensitivity} give WCS for a number of commonly used uncertainty sets, which we have derived in this section.
The collection of sensitivity measures can be very different; for example, ``budgeted" sensitivity measures the spread of the left tail of the cost distribution whereas the ``convex combination" is the spread of the right. A decision can be robust under one robustness measure but not another. 

The  uncertainty set affects the DRO decision through  the robustness measure (regularizer) induced by the uncertainty set  \eqref{eq:DRO-mean-sensitivity} and should be chosen according to the sensitivity measure that needs to be reduced. For example,  budgeted uncertainty should be chosen if the DM wants to reduce the left  tail of the cost distribution whereas a convex-combination should be chosen to reduce the right. While reducing  its own sensitivity measure, a DRO solution can  increase other  sensitivity measures, making it less robust than the nominal solution from the perspective of other uncertainty sets (we illustrate this in the following section).

From \eqref{eq:DRO-mean-sensitivity} the family of robust solutions  maps out a nearly Pareto-optimal tradeoff between performance and robustness which  can be used not only to select the size of the uncertainty set but to identify conservative problem instances where the price of robustness is high. In such cases, the form of the sensitivity measure provides guidance on how to redesign the problem to improve the performance--robustness tradeoff, which we illustrate in Section \ref{sec:hedging}.

\section{Examples}
\label{sec:Examples}

We consider two examples, inventory control and CVaR minimization.  Both illustrate how worst-case sensitivity can be used to select the family and size of uncertainty sets in a DRO model. In particular sensitivity measures can be used to identify vulnerabilities of the nominal solution and the uncertainty set that most efficiently controls these sensitivities. Mean--sensitivity frontiers are used to choose the size of the uncertainty set. 

Both examples show that the   sensitivity measure of the uncertainty set can have  a substantial impact on the DRO solution; a decision is ``robust" with respect to the stress test defined by an uncertainty set and DRO solutions reduce sensitivity with respect to its uncertainty set. For the inventory problem, we even have conflicting recommendations on how to be ``more robust." (There are also examples where robust solutions under different uncertainty sets are almost the same (Appendix \ref{sec:LR}).) When the price of robustness is high, which is the case in the inventory application, the sensitivity measure provides guidance on redesigning the system, in this case with ``return contracts," to improve the tradeoff.

\subsection{Inventory control}

Consider once again the inventory cost function \eqref{eq:inv2}. In this example, we show how a robustness measure can be used to select both the family and size of the uncertainty set by (i) looking at the distribution of the cost under the SAA solution to select the uncertainty set according to the sensitivity control that is needed, and (ii) using the mean--sensitivity frontier to select its size. We also show that the choice of uncertainty set can have a substantial impact on the DRO solution and that an injudicious choice of uncertainty set can make the decision ``less robust." 
We generated $n=100$ demand realizations $\{Y_1,\cdots,\,Y_n\}$ by sampling from a mixture of two exponential distributions with means $\mu_L=10$ and $\mu_H=100$, where the probability of a sample from population $L$ is 0.9. We assume $r=10$, $c=2$, $q=0$ and $s=4$. 

\begin{figure}[h]
\centering
\includegraphics[scale=0.1]{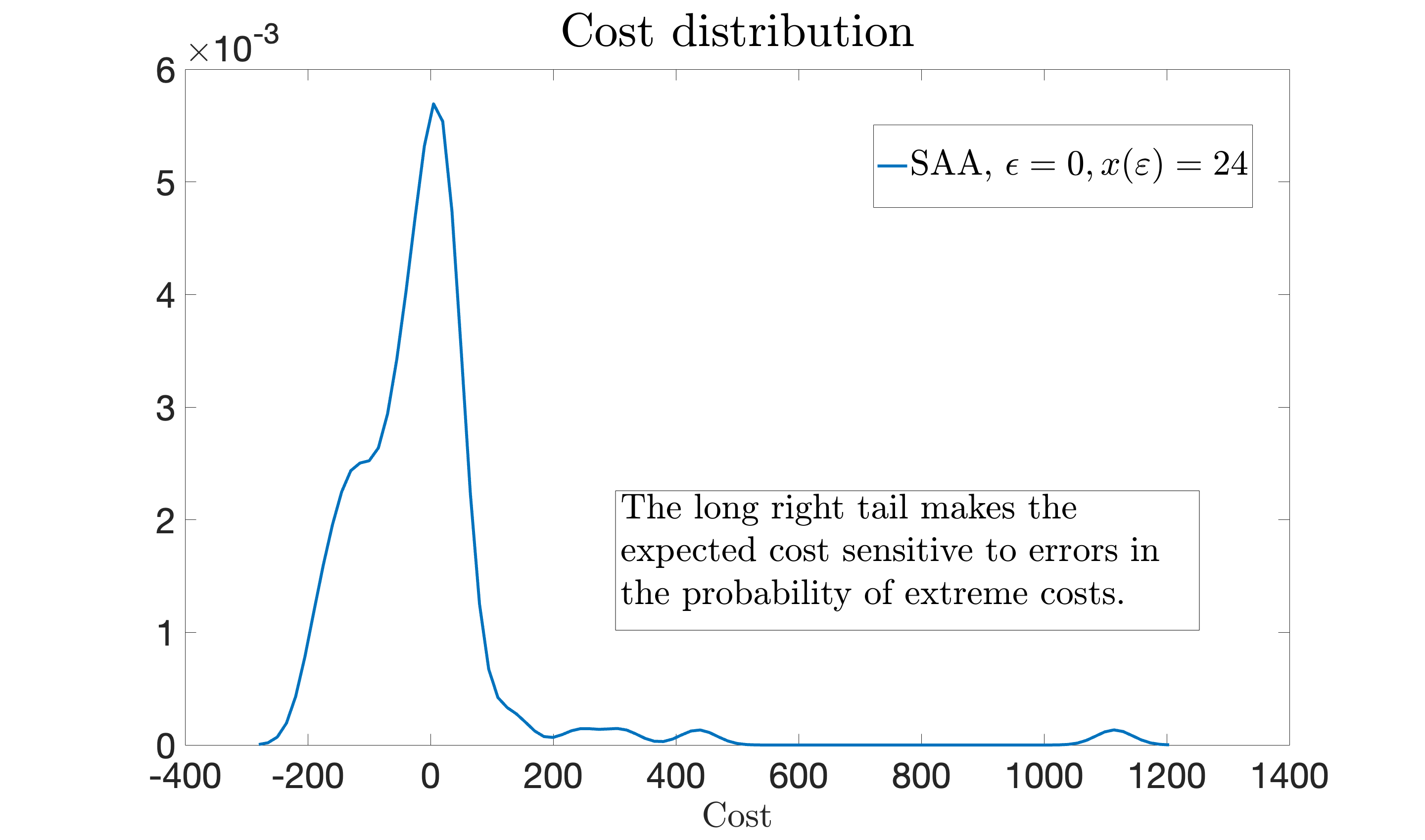}
\caption{Distribution of the cost under the optimal order quantity for the SAA problem.}
\label{fig:hist_cost}
{\raggedright\footnotesize
The cost has a long right tail, which suggests we should choose an uncertainty set with a sensitivity measure that shrinks this tail.
\par}
\end{figure}

We begin by looking at the distribution of the cost under the SAA solution to determine whether there is a need for sensitivity control.   Figure \ref{fig:hist_cost} shows  that the resulting cost distribution has a long right tail, meaning that the expected cost is sensitive to errors in tail probabilities. It follows that we should choose an uncertainty set in the DRO problem that has a WCS that depends on the right tail. From Table \ref{table:summary}, all uncertainty sets seem reasonable except for the ``budgeted" uncertainty set which has a sensitivity that penalizes the left tail but does not depend on the right. For this experiment, we select the modified $\chi^2$-uncertainty set.

\begin{figure}[h]
\centering
\includegraphics[scale=0.1]{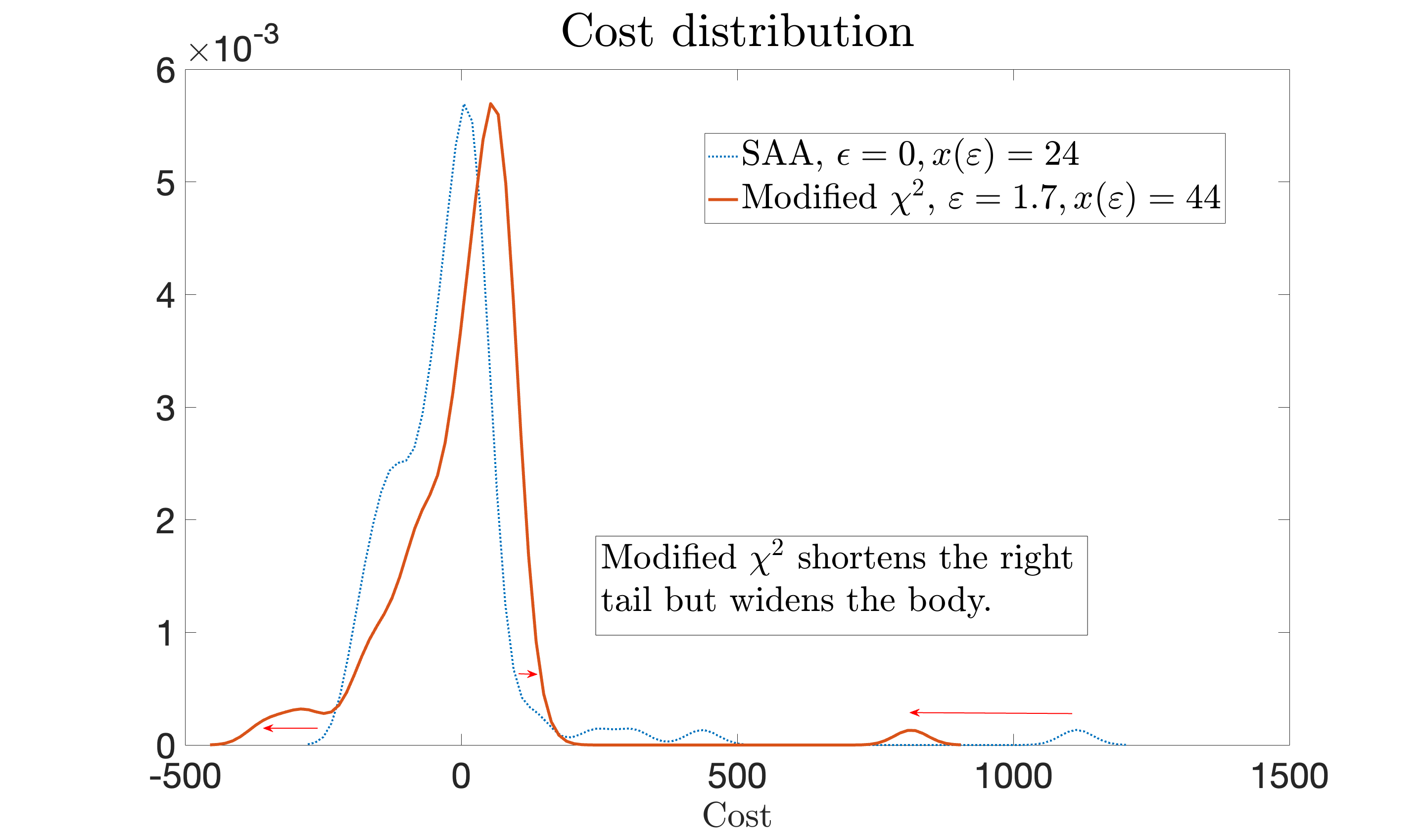} 
\caption{Distribution of the cost with the modified $\chi^2$-uncertainty set.}
\label{fig:hist_chi2}
{\raggedright\footnotesize
 Worst-case sensitivity is the standard deviation of the cost distribution. The DRO solution reduces the standard deviation (sensitivity) by shrinking the right tail. 
\par}
\end{figure}

\begin{figure}
\includegraphics[scale=0.1]{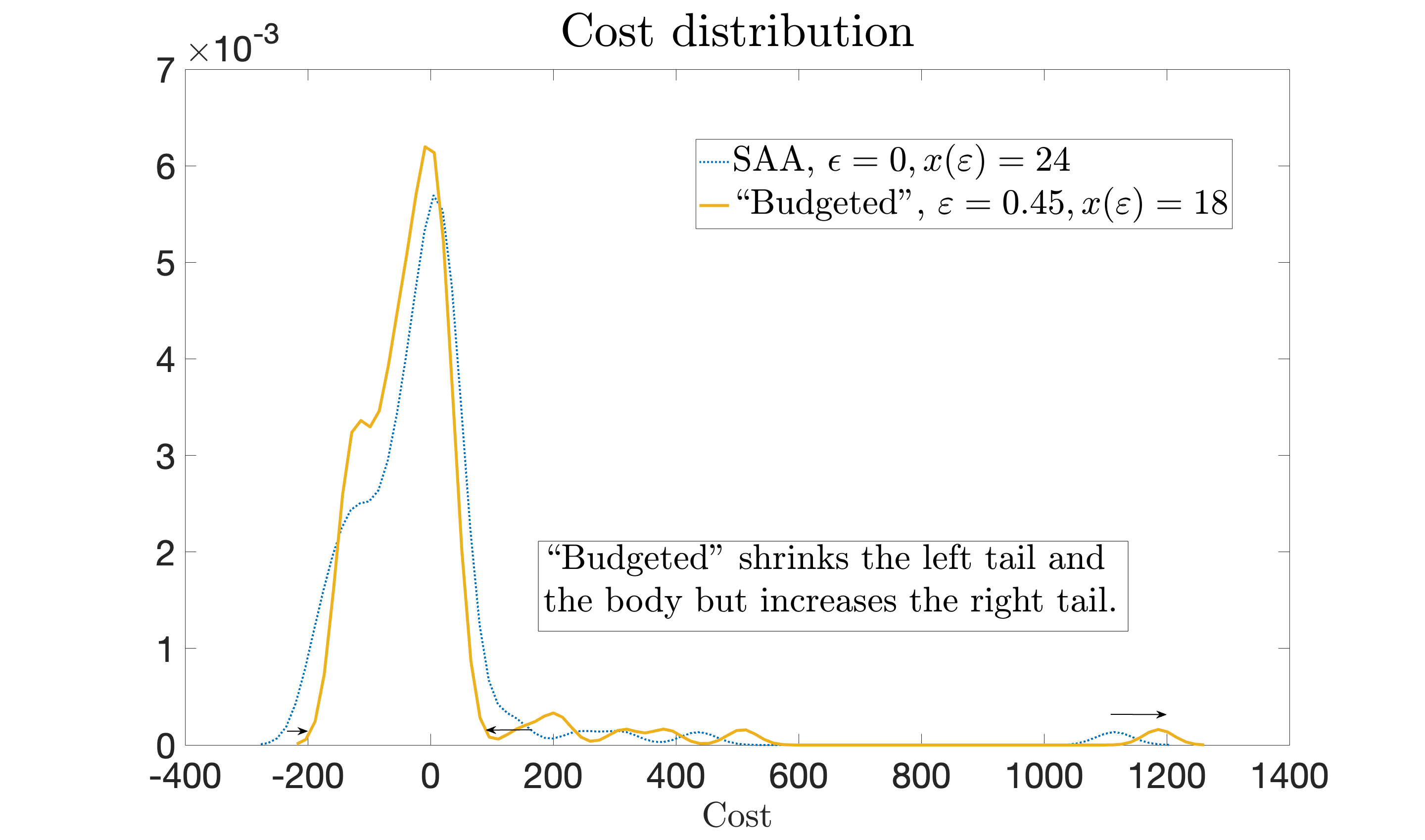}
\caption{Distribution of the cost with the ``budgeted" uncertainty set.}
\label{fig:hist_budgeted}
{\raggedright\footnotesize
 Under this uncertainty set, worst-case sensitivity is the size of the left tail. The solution of the DRO problem shrinks the left tail (sensitivity), which increases the right tail of the distribution. 
\par}
\end{figure}

Figure \ref{fig:hist_chi2} shows the cost distribution under the solution of a DRO problem with a modified $\chi^2$-uncertainty set with $\varepsilon=1.7$. As expected, the DRO solution shrinks the right tail, which comes at the cost of a ``wider body" and a larger expected cost. Nevertheless, it reduces the sensitivity (standard deviation) as desired. The optimal order quantity for the modified $\chi^2$-uncertainty set ($x(\varepsilon)=44$) is larger than for SAA ($x(0)=24$). 

To highlight the importance of the uncertainty set on the DRO solution, Figure \ref{fig:hist_budgeted} shows the cost distribution for the DRO solution with a budgeted uncertainty set ($\varepsilon=0.45$). In this case, the DRO solution shrinks the left-tail of the cost distribution which corresponds to reducing the sensitivity, $\mathbb{E}_{\mathsf{p}}(\mathsf{f})-\min(\mathsf{f})$, of the budgeted uncertainty set. While this reduces sensitivity to perturbations in left-tail probabilities, it increases the length of the right tail which magnifies the sensitivity to errors in right tail probabilities. From the perspective of the modified $\chi^2$ uncertainty set, a ``robust" solution obtained using budgeted uncertainty is less robust than the SAA solution. The optimal order quantity with budgeted uncertainty ($x(\varepsilon)= 18$) is smaller than the SAA solution ($x(0)=24$). It was larger ($x(\varepsilon)=44$) with the modified $\chi^2$ uncertainty set.  This example shows that a decision is only robust {\it with respect to an uncertainty set}, and that the DRO solution with respect to one uncertainty set can be less robust than the nominal solution under another uncertainty set. The choice of uncertainty set should align with the sensitivity measure the DM would like to reduce. 


\begin{figure}
\includegraphics[scale=0.1]{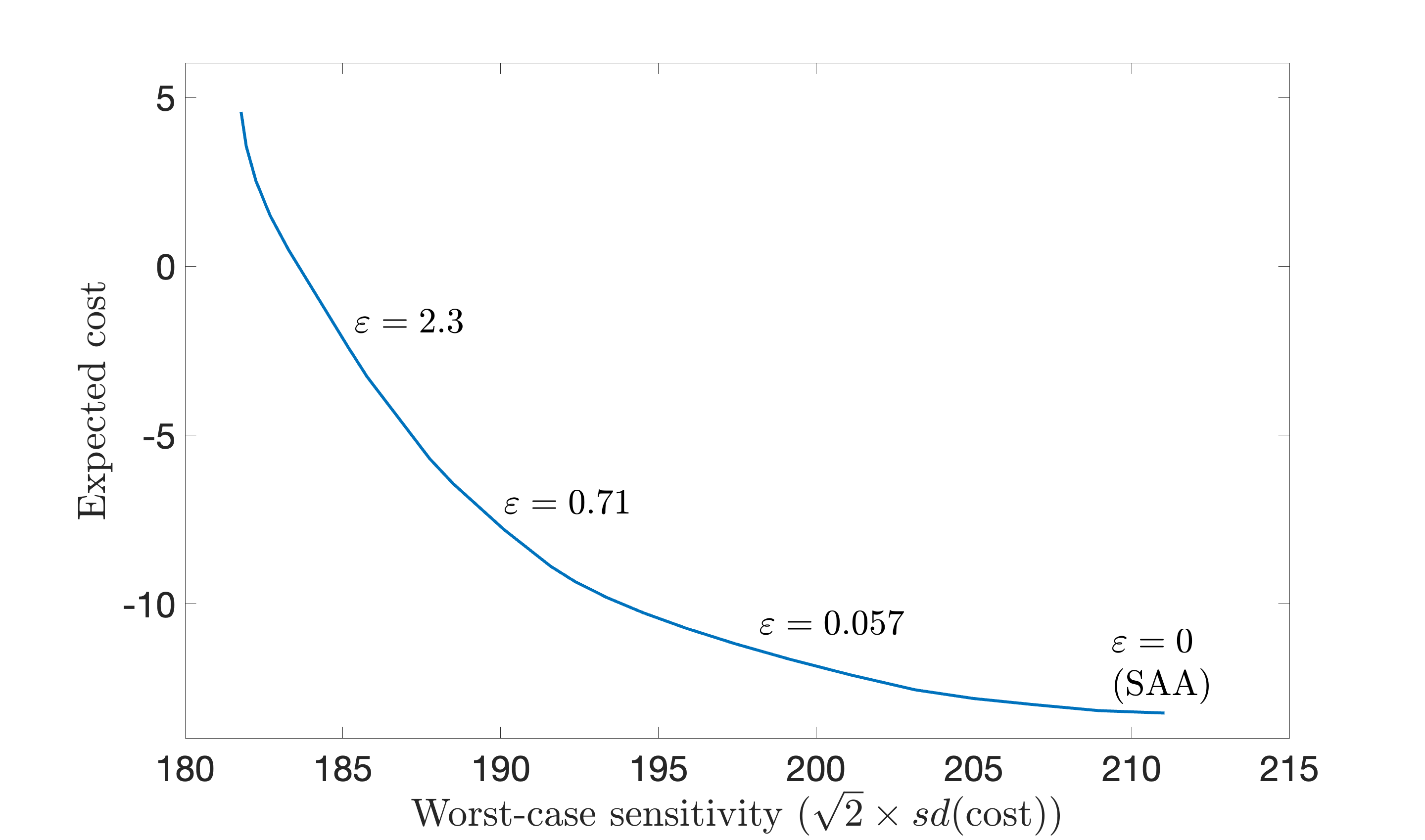}
\caption{Mean--sensitivity frontier for DRO solutions for a modified $\chi^2$ uncertainty set.}
\label{fig:frontier-mean_chi2}
\end{figure}

Figure \ref{fig:frontier-mean_chi2} shows the mean--sensitivity frontier generated by solutions of the DRO problem with a modified $\chi^2$ uncertainty set.  If $\varepsilon=0.114$ in the DRO problem, sensitivity is $190$, a reduction of $7.7\%$ from its value of $211$ under the SAA solution; expected cost is $-10.3$, an increase of $22\%$ from its minimum of $-13.2$. As with any tradeoff, there is no ``correct" choice of $\varepsilon$; it depends on the decision maker and the context of the application.

For this example, KL-divergence, total variation, and the convex-combination ($\alpha=0.9$) generate the same set of robust solutions because the decision variable is scalar. It will be seen in Section \ref{sec:CVaRexample} that this is not generally the case.

\subsection{Hedging for robustness}
\label{sec:hedging}

We see from Figure \ref{fig:frontier-mean_chi2} that  the price of robustness is high; a relatively modest decrease in sensitivity comes with a large increase in expected cost. This is a property of the problem instance: different parameters ($r$, $q$, $s$ and $c$) may result in a less conservative (flatter) frontier. 

We can improve the cost--sensitivity tradeoff  by directly targeting WCS. For modified $\chi^2$ divergence, where the uncertainty set is the standard deviation, an additional income stream with a small expected value that is negatively correlated with the cost, like the option to return leftover items to the supplier for an upfront fee, will do the job
\begin{eqnarray*}
f(x, Y) -  \underbrace{\alpha \big\{\max(x-Y, 0)- {\mathbb E}_{\mathbb P}\left[\max(x-Y, 0) \right]\big\}}_{\tiny \mbox{Hedge that is negatively correlated with $f(x, Y)$}}.
\end{eqnarray*}
A large cost $f(x, Y)$ due to low sales ($Y<x$) is offset by income $\alpha \max(x-Y, 0)$ from  returning unsold items for $\alpha$ per item; the price of this option is $\alpha{\mathbb E}_{\mathbb P}\left[\max(x-Y, 0) \right]$.  Figure \ref{fig:frontier-mean_chi2_hedging} shows a flatter (less conservative) frontier plot after hedging ($\alpha=0.75$).  In reality the price of the option is likely to exceed its expected value and this additional cost will lift the frontier; we can also optimize over $\alpha$. We ignore these embellishments to highlight the key idea which is the benefit of directly targeting the robustness (sensitivity) measure.  The  cost function \eqref{eq:inv2} with positive salvage value $q$ (other parameters the same) is a special case of a return option where the option  is free. This simultaneously reduces expected cost while flattening the frontier.

This example illustrates how the sensitivity measure  provides guidance  on how to redesign the system when the DRO solution is conservative. Other sensitivity measures may require a different hedge. Figure \ref{fig:histchi2-hedged} compares the distribution of the cost for the hedged and unhedged DMs when the order quantity is $44$.

\begin{figure}
\includegraphics[scale=0.27]{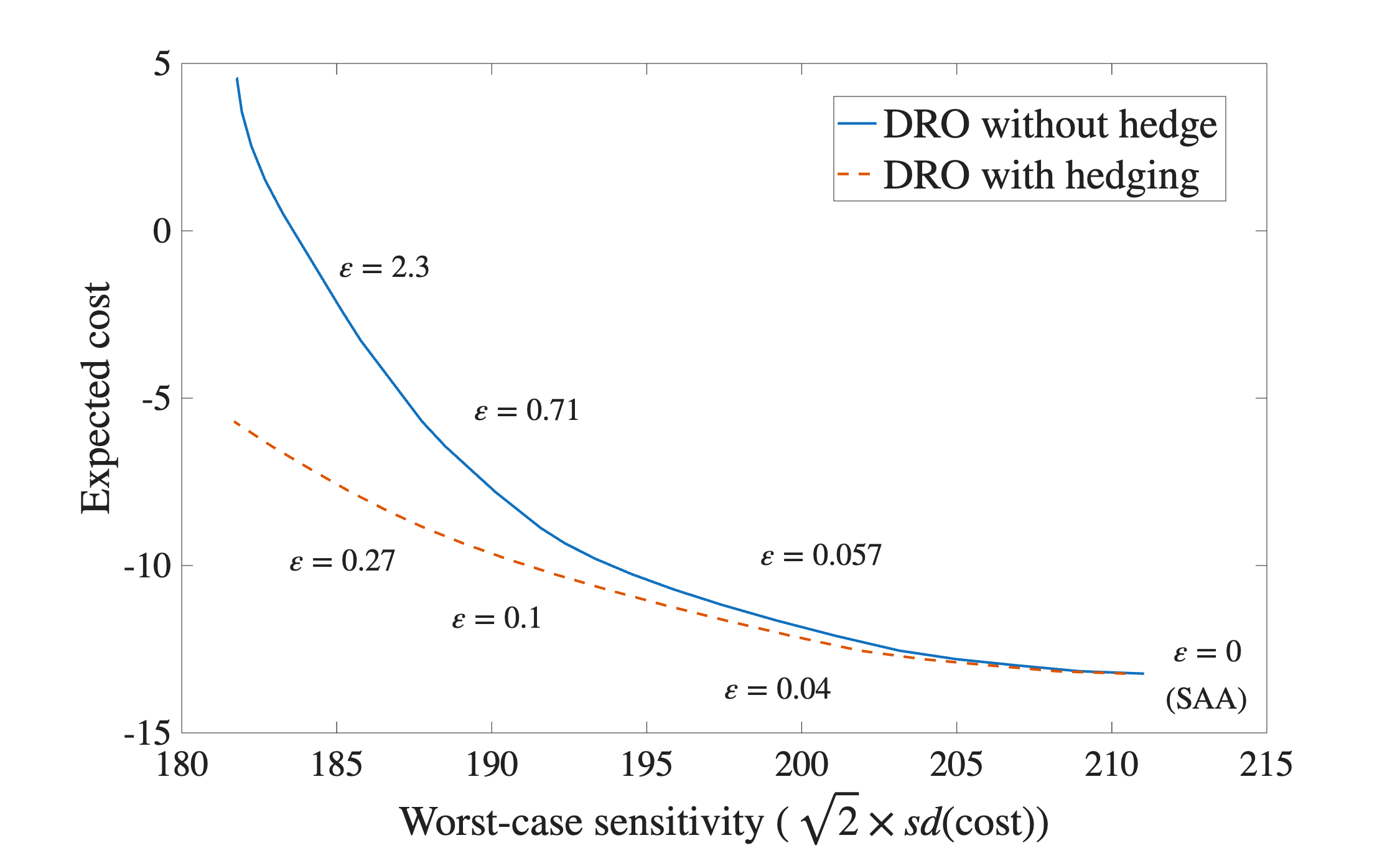}
\caption{Improving mean--sensitivity frontier through hedging}
\label{fig:frontier-mean_chi2_hedging}
{\raggedright\footnotesize
WCS for modified $\chi^2$ divergence can be reduced by providing the DM with an income stream that is negatively correlated with the cost. One example of  a hedge is the option to return items at a reduced price. This reduces the price of robustness, as seen by the flatter frontier for the hedged DM. 
\par}
\end{figure}

\begin{figure}
\includegraphics[scale=0.29]{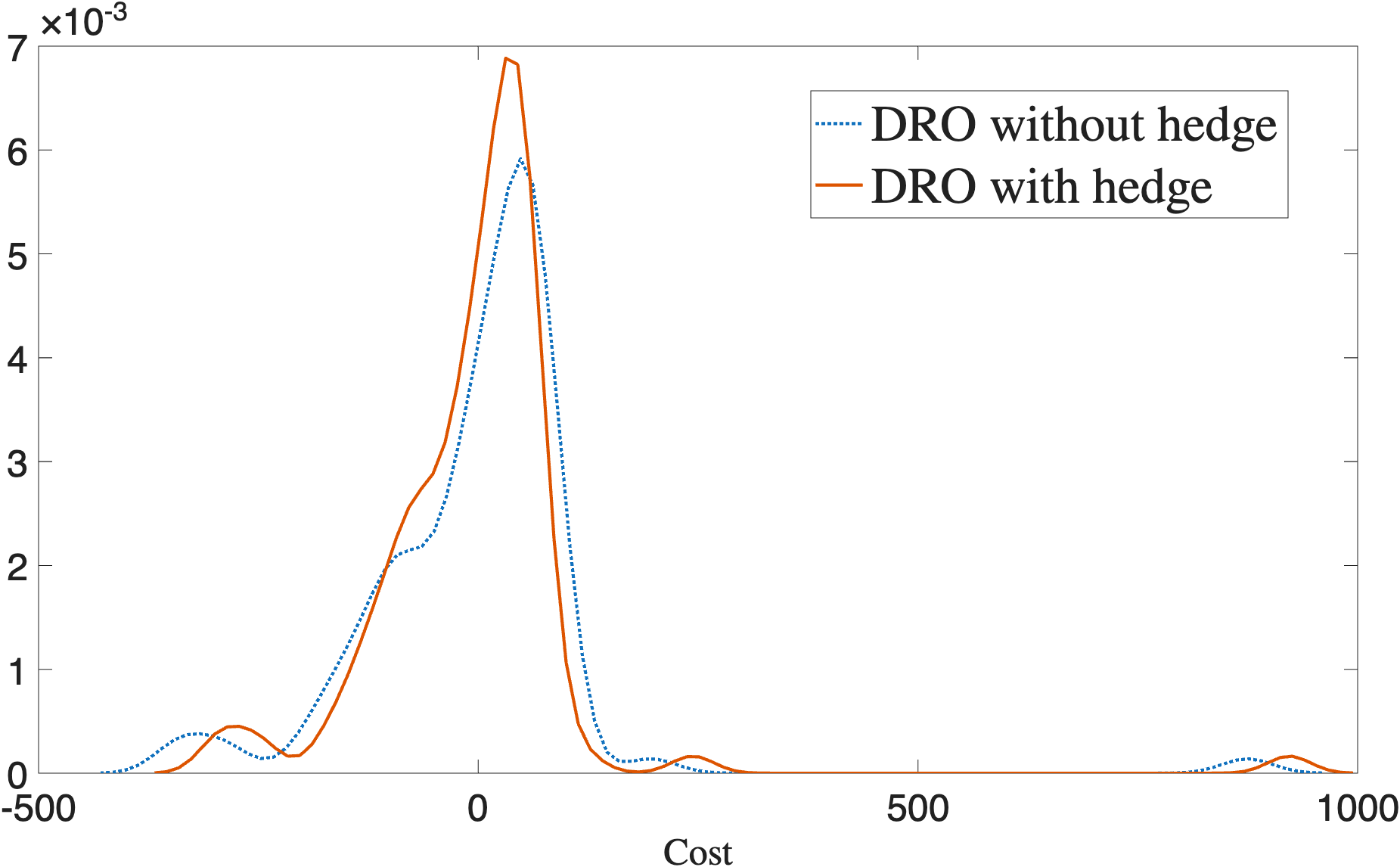}
\caption{Histogram of unhedged and hedged cost.}
\label{fig:histchi2-hedged}
{\raggedright\footnotesize
 The hedge enables the DM to reduce WCS (standard deviation) with minimal impact on the expected cost. In both cases, the order quantity is $44$ and the expected cost is $-3.3$. WCS is $186$ (unhedged) and $178$ (hedged).
\par}
\end{figure}

\subsection{Minimum CVaR portfolio}
\label{sec:CVaRexample}
Consider the robust minimum 90\%-CVaR portfolio selection: 
\begin{align*}
\min_{\mathsf{x}\in X}\max_{\mathsf{q}\in\mathcal{Q}(\varepsilon)} \mathrm{CVaR}_{\mathsf{q},0.9}(-\mathsf{R}\mathsf{x}),
\end{align*}
where $\mathsf{R}\in\mathbb{R}^{n\times d}$ is a rate of return data set and $X=\{\mathsf{x}\in\mathbb{R}^d\vert\mathsf{1}^\top\mathsf{x}=1\}$ is the feasible set of portfolio $\mathsf{x}$. We use the 30 industry portfolio data set \cite{french}, which 
consists of $d=30$ industry stock indices and $n=1082$ historical monthly rate of returns. Minimum CVaR optimization for the ambiguity neutral  problem ($\varepsilon=0$) is sensitive to model misspecification \cite{Ban2016CVaR,BLS}.

\begin{figure}[h]
\centering
\includegraphics[width=4in]{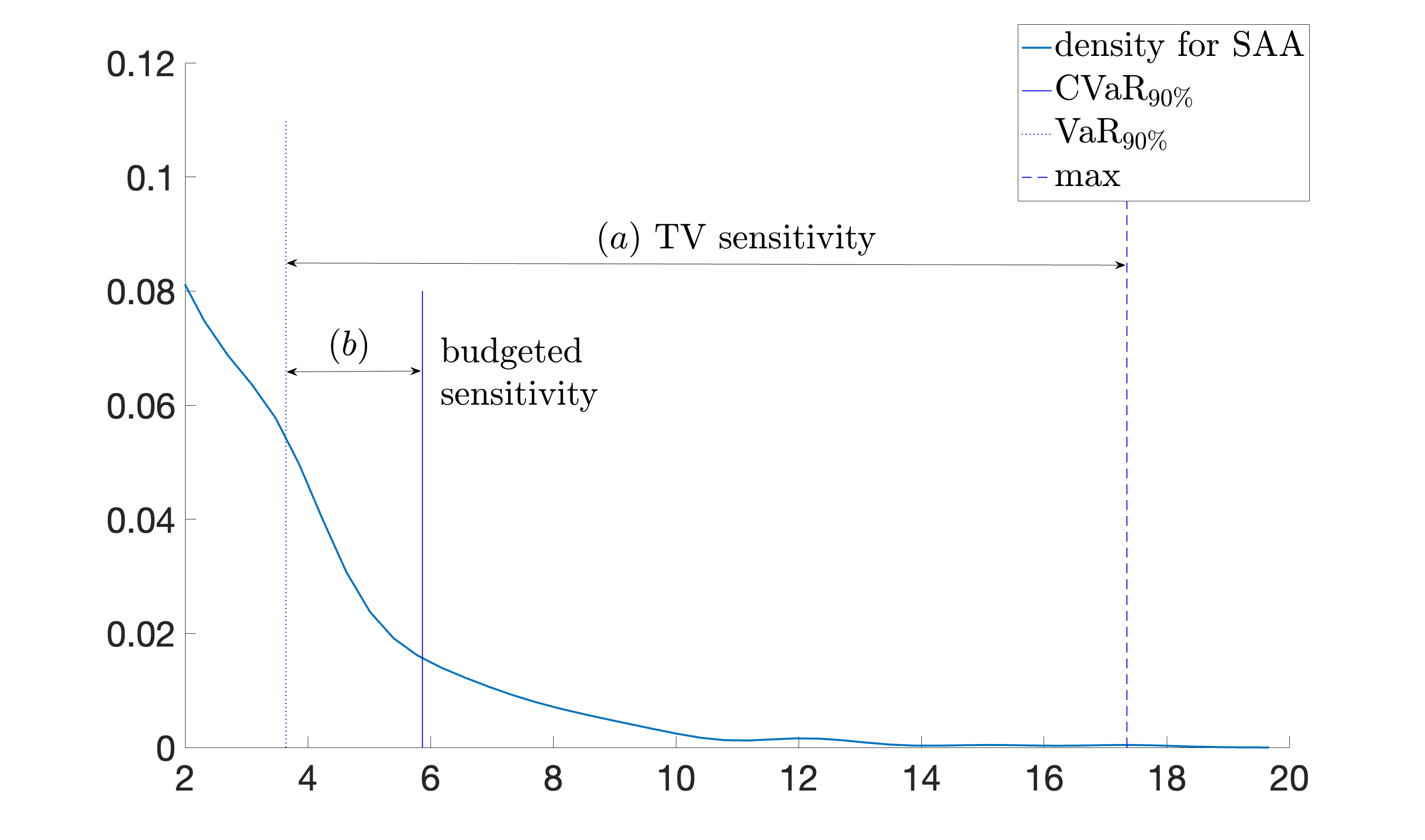}
\caption{Tail of the cost distribution of the  SAA solution}
\label{fig:SAA-hist}
{\raggedright\footnotesize
 $(a)$ is the WCS for the TV uncertainty set and $(b)$ the WCS for the budgeted uncertainty set.
\par}
\end{figure}

We first solve the (ambiguity-neutral) minimum CVaR problem (i.e., $\varepsilon=0$). Figure \ref{fig:SAA-hist} shows the kernel density estimate of the right tail of the associated loss distribution. The long right tail suggests that for this data set, optimal ambiguity-neutral CVaR is sensitive to worst-case deviations defined by a total-variation uncertainty set (recall from Figure \ref{fig:sensitivity}$\langle$ii$\rangle$ that TV-sensitivity is proportional to ``$\max(\mathsf{f})-\mathrm{VaR}_{\mathsf{p},\beta}(\mathsf{f})$,''  indicated as $(a)$ in the plot). Also shown is the worst-case sensitivity for  the budgeted uncertainty set $(b)$.

The mean--sensitivity tradeoff suggests that TV- (budgeted-) sensitivity can be controlled most efficiently by solving a robust CVaR problem with the TV- (budgeted-) uncertainty set. To illustrate this, we compare the loss distributions of empirical CVaR minimization and   the robust CVaR problem  with  TV-,  budgeted- and modified $\chi^2$-uncertainty sets.

\begin{figure}[h]
\centering
\includegraphics[width=5in]{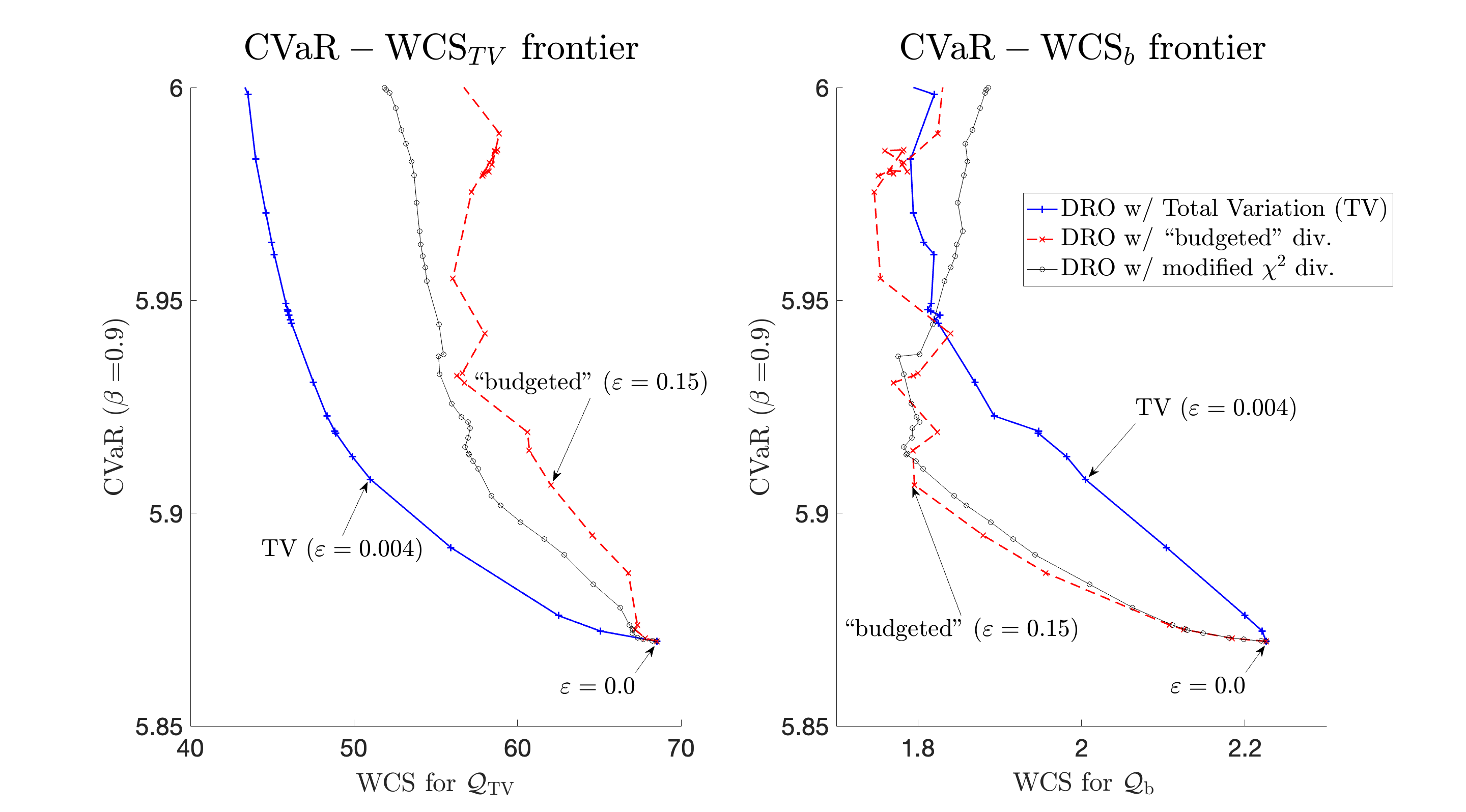}
\caption{Frontiers with different sensitivities}
\label{fig:cvar_wcs_frontiers_comp09_two}
{\raggedright\footnotesize
The plot on the left shows the CVaR--TV sensitivity frontiers  generated by solutions of the robust $\mbox{CVaR}_{90\%}$ problem with TV, budgeted, and modified $\chi^2$ uncertainty sets. All frontiers are different because the robust problems have different solutions.  Not surprisingly, the most efficient frontier is generated with the TV uncertainty set, though robust solutions with the budgeted and modified $\chi^2$ uncertainty sets also reduce TV-sensitivity.  The plot on the right shows CVaR--budgeted sensitivity frontiers for the same robust solutions. The most efficient frontier is now generated by the budgeted uncertainty set. In both cases, robust solutions with the modified $\chi^2$ uncertainty set generate a frontier that is intermediate between the other two. 
\par
}
\end{figure}

The left plot of Figure \ref{fig:cvar_wcs_frontiers_comp09_two} 
 shows the CVaR--TV sensitivity frontier generated by solutions of the robust CVaR problem with the TV-, budgeted- and modified $\chi^2$-uncertainty sets (i.e., we solve the robust CVaR problem with each uncertainty set for a collection of $\varepsilon$, and  compute the CVaR--TV frontier for the collection of decisions). CVaR is $5.91$ when $\varepsilon=0.004$ for DRO with a TV uncertainty set and $\varepsilon=0.15$ when a budgeted uncertainty set is adopted. All frontiers are different because  the uncertainty sets are different. Of course, the most efficient CVaR--TV sensitivity frontier is generated by solving the robust CVaR problem with a TV uncertainty set, though unlike the inventory problem, robust solutions for the other two uncertainty sets also reduce TV-sensitivity. The right plot of Figure \ref{fig:cvar_wcs_frontiers_comp09_two} 
shows the mean--budgeted sensitivity frontiers for the same DRO problems/solutions; budgeted solutions are now the most efficient. This example shows again there is no ``universal" notion of robustness; the choice of uncertainty set matters and affects the DRO solution and all sensitivities differently. The budgeted uncertainty set reduces budgeted sensitivity most efficiently but is not a good choice if a large nominal TV-sensitivity is the main concern. It makes sense to monitor multiple frontiers (Table \ref{table:summary}) when selecting the size of an uncertainty set.


\begin{figure}[h]
\centering
\includegraphics[width=4in]{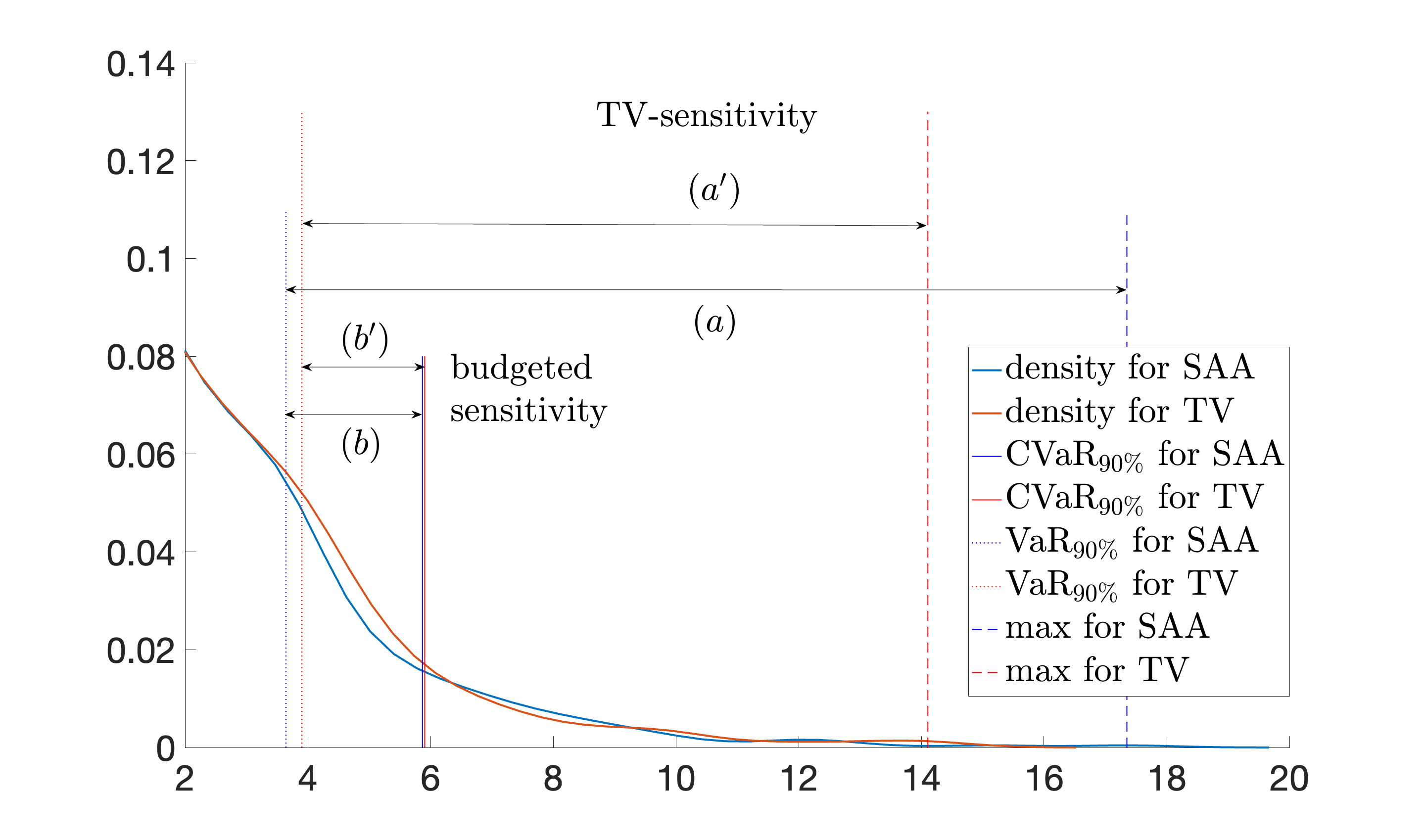}
\caption{Tails of cost distributions of empirical $\mbox{CVaR}_{90\%}$ and DRO $\mbox{CVaR}_{90\%}$ with the TV uncertainty set ($\varepsilon=0.004$)}
\label{fig:SAA-DROTV-hist}
{\raggedright\footnotesize
The size of the uncertainty set is chosen so that CVaR is approximately $5.91$. $(a)$ is TV-sensitivity and $(b)$ budgeted sensitivity for the SAA and robust CVaR solutions, respectively;  $(b)$ and $(b')$ are the respective budgeted sensitivities. The robust solution reduces TV-sensitivity ($(a)$ versus $(a')$) and budgeted-sensitivity ($(b)$ versus $(b')$) of the (ambiguity neutral) CVaR minimizer.
\par}
\end{figure}

\begin{figure}[h]
\centering
\includegraphics[width=4in]{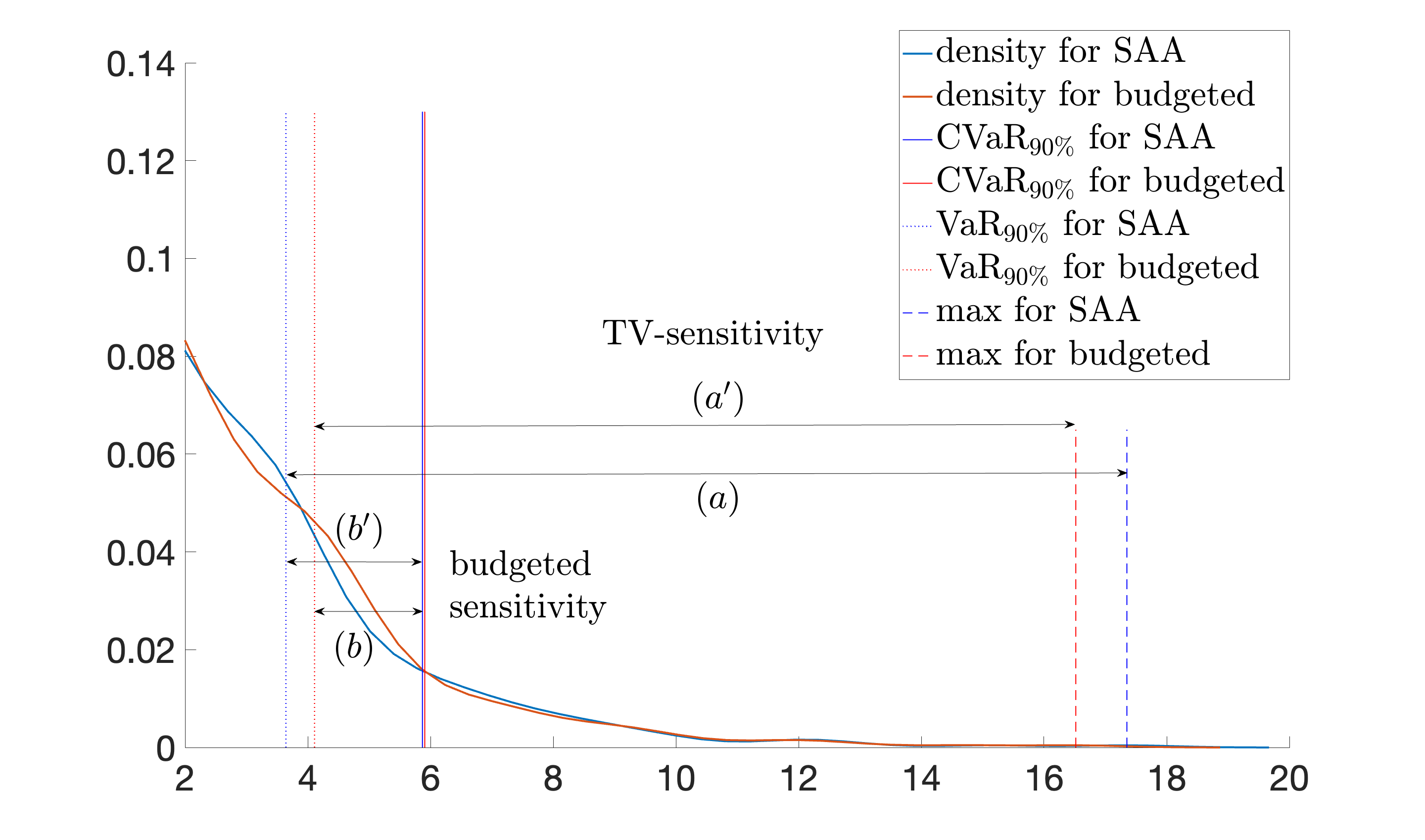}
\caption{Tails of cost distributions of empirical $\mbox{CVaR}_{90\%}$ and DRO $\mbox{CVaR}_{90\%}$ with budgeted uncertainty set ($\varepsilon=0.15$)}
\label{fig:SAA-DRObudgeted-hist}
{\raggedright\footnotesize
 The size of the uncertainty set is chosen so that CVaR is approximately $5.91$.  $(a)$ is TV-sensitivity amd $(b)$ budgeted sensitivity for the SAA and robust CVaR solutions, respectively;  $(b)$ and $(b')$ are the respective budgeted sensitivities. The robust solution reduces TV-sensitivity ($(a)$ versus $(a')$) and budgeted-sensitivity ($(b)$ versus $(b')$) of the (ambiguity neutral) CVaR minimizer. 
\par}
\end{figure}

Figure \ref{fig:SAA-DROTV-hist} superimposes the tail of the loss distribution of the solution of the
robust CVaR problem with a TV uncertainty set ($\varepsilon=0.004$) on the ambiguity neutral distribution from Figure \ref{fig:SAA-hist}, while Figure \ref{fig:SAA-DRObudgeted-hist} does the same but now with the  solution of the DRO-budgeted ($\varepsilon = 0.15$) problem. $\varepsilon$ is chosen so that CVaR is approximately the same ($5.91$) for both robust solutions. Consistent with the frontier plots from Figure \ref{fig:cvar_wcs_frontiers_comp09_two}, the solution of the DRO-TV (DRO-budgeted) problem reduces both TV-  and budgeted-sensitivity, but the reduction in TV-sensitivity is larger.

\section{Conclusion}

Though formulated as a single-objective worst-case problem, DRO can be understood  as a tradeoff between performance and  robustness. The robustness objective (worst-case sensitivity) is determined by the choice of uncertainty set; when we choose an uncertainty set, we select a robustness measure. DRO  maps out a (nearly) Pareto-optimal performance--robustness frontier as the size of the uncertainty set is varied. 

This interpretation of DRO resolves many of the misgivings about the use of  worst-case models when there is model uncertainty: That it is too conservative with no systematic approach for selecting the family or size of the uncertainty set or clear insight on why  the resulting decision is even ``robust" (particularly when different uncertainty sets give contradicting recommendations).

Worst-case sensitivity quantifies errors in the nominal model that have large sensitivities and the uncertainty sets that reduce these sensitivities most effectively. The size of the uncertainty set can be determined from the underlying performance--robustness  frontier.
It shows that DRO as a framework is not inherently ``conservative" but a computationally tractable tool for mapping out the performance--robustness tradeoff. Under standard assumptions, worst-case sensitivity is a generalized measure of deviation. Intuitively, the expected cost under the nominal distribution can be sensitive to changes in the probability assigned to extreme values so a cost distribution with a large spread is not robust. 

We derive sensitivity measures for expected cost and for a risk measure for a collection of commonly used uncertainty sets. This makes explicit the robustness measure that is controlled for each choice of uncertainty set. The performance--sensitivity frontier identifies {\it conservative problem instances}, ones where the price of robustness is high. In such cases,  ``hedging instruments" that directly target sensitivity reduction can be used to reduce the price of robustness and improve the tradeoff.

By quantifying robustness and putting the  tradeoff between performance and robustness on center stage, we provide a new way of understanding  Distributionally Robust Optimization and building robust optimization models. We focus on the single-period problem as it is sufficient to convey the key ideas and is the setting for much of the DRO literature.  The main insights, however, go well beyond the setting of this paper and generalizations will be reported elsewhere.

\paragraph{\bf Acknowledgments}
Jun-ya Gotoh is supported in part by the MEXT Grant-in-Aid 24K01113.
Michael Kim is supported in part by the Natural Sciences and Engineering Research Council (NSERC) Discovery Grant RGPIN-2015-04019.
Andrew Lim is supported by the  Ministry of Education, Singapore, under its 2021 and 2024 Academic Research Fund Tier 2 grant calls (Award ref: MOE-T2EP20121-0014 and MOE-T2EP20224-0018).

\bibliographystyle{plainnat}
\bibliography{references}

\newpage

\appendix

\section{Proofs from Section \ref{sec:DRO_sensitivity}}

\subsection{Proof of Proposition \ref{prop:deviation}}\label{sec:proof_WCS_is_GMD}
Without loss of generality, we assume ${\mathbb E}_{\mathbb P}[f(x,Y)]=0$.
\paragraph{${\mathcal A}(\varepsilon;f)$ is a generalized deviation measure}
Clearly, ${\mathcal A}(\varepsilon;f)$ satisfies properties (2) and (3) of Definition \ref{def:deviation}. As for property (1), ${\mathcal A}(\varepsilon;f)\geq 0$ follows from the observation that the nominal measure ${\mathbb P}\in{\mathcal Q}(\varepsilon)$ is feasible so
\begin{align*}
{\mathcal A}(\varepsilon;f) \geq {\mathbb E}_{\mathbb P}[f(Y)]= 0.
\end{align*}
It is also clear that ${\mathcal A}(\varepsilon;f)=0$ if $f\equiv 0$. We now show the converse.

Suppose that ${\mathcal A}(\varepsilon;f)=0$, ${\mathbb E}_{\mathbb P}[f]=0$ but that $f$ is not equal to $0$. 
Clearly, we can select $\delta>0$ such that 
\begin{align*}
\frac{{\rm d}{\mathbb Q}}{{\rm d}{\mathbb P}}^{(\delta)} = 1 + \delta f(Y)
\end{align*}
is strictly positive and defines a probability measure ${\mathbb Q}^{(\delta)}$ such that $d({\mathbb Q}^{(\delta)}\,|\, {\mathbb P}) < \varepsilon$.
It follows that
\begin{eqnarray*}
{\mathcal A}(\varepsilon;f) & \geq  &{\mathbb E}_{{\mathbb Q}^{(\delta)}}[f(Y)] \\
& = & {\mathbb E}_{\mathbb P}\Big\{\frac{{\rm d}{\mathbb Q}}{{\rm d}{\mathbb P}}^{(\delta)} f(Y)\Big\} \\
& = & {\mathbb E}_{\mathbb P}[f(Y)] + \delta {\mathbb V}_{\mathbb P}(f(Y)) \\
& > & 0
\end{eqnarray*}
which contradicts our initial assumption that ${\mathcal A}(\varepsilon;f)=0$. It follows that $f$ must equal $0$.
Since ${\mathcal A}(\varepsilon;f)$ is a measure of deviation, so too is ${\mathcal S}(\varepsilon;f)$.

\paragraph{${\mathcal S}_{\mathsf p}(\mathsf{f})$ is a generalized measure of deviation}
We first establish the growth condition on ${\mathcal A}(\varepsilon;f)$. Suppose there is a constant $k>0$ such that the continuity condition \eqref{eq:continuity d} holds for every ${\Delta}:\Omega \rightarrow {\mathbb R}$ satisfying ${\mathbb E}_{\mathbb P}[\Delta]=0$. Given $f$ and $\varepsilon>0$, let ${\mathbb Q}^{(\varepsilon)}$ denote the solution of the worst-case problem \eqref{eq:V}. Since $d({\mathbb Q}^{(\varepsilon)}\, | \,{\mathbb P}) = \varepsilon$, it follows that 
\begin{eqnarray*}
\frac{{\rm d}{\mathbb Q}}{{\rm d}{\mathbb P}}^{(\varepsilon)} = 1 + \varepsilon^\frac{1}{k}\xi \geq 0
\end{eqnarray*}
for some $\xi$ such that ${\mathbb E}_{\mathbb P}[\xi]=0$, so
\begin{align*}
{\mathcal A}(\varepsilon;f) = {\mathbb E}_{{\mathbb Q}^{(\varepsilon)}}[f(Y)] = {\varepsilon}^\frac{1}{k} ~ {\mathbb E}_{\mathbb P}[\xi f(Y)]
\end{align*}
We now show that ${\mathcal S}_{\mathbb P}(\mathsf{f})$ is a generalized measure of deviation. 
Clearly, ${\mathcal S}_{\mathbb P}(\mathsf{f})$ satisfies conditions (2) and (3) of Definition \ref{def:deviation}.
As for Condition (1), it is easy to show that ${\mathcal S}_{\mathbb P}({f(Y)})\geq 0$ and that ${\mathcal S}_{\mathbb P}({f})=0$ if $f$ is equal to $0$; these properties follow from those of ${\mathcal A}(\varepsilon;f)$. To prove the converse,  suppose that
\begin{align}
{\mathcal S}_{\mathsf p}(\mathsf{f}) := \lim_{\varepsilon \downarrow 0}\frac{{\mathcal A}(\varepsilon;f)}{\varepsilon^{\frac{1}{k}}} = 0
\label{eq:contradiction1}
\end{align}
but that $f
$ is not equal to $0$.
Since $d({\mathbb Q}|{\mathbb P})$ is continuous in ${\mathbb Q}$,  there is a constant $\delta>0$ such that
\begin{align*}
\frac{{\rm d}{\mathbb Q}}{{\rm d}{\mathbb P}}^{(\varepsilon)} = 1 + \varepsilon^\frac{1}{k}\delta f(Y) > 0
\end{align*}
and from the continuity assumption \eqref{eq:continuity d}
\begin{align*}
d\big({\mathbb Q}^{(\varepsilon)} \,\big|\, {\mathbb P}\big) \leq \varepsilon
\end{align*}
for all $\varepsilon>0$ sufficiently small.
Since ${\mathbb Q}^{(\varepsilon)}\in{\mathcal Q}(\varepsilon)$, it follows that
\begin{align*}
0 < {\mathbb E}_{{\mathbb Q}^{(\varepsilon)}}[f(Y)] = {\mathbb E}_{\mathbb P}[f(Y)] +  \varepsilon^\frac{1}{k} \delta {\mathbb V}_{\mathbb P}(f(Y)) \leq V(\varepsilon)
\end{align*}
so
\begin{align*}
\frac{{\mathcal A}(\varepsilon)}{ \varepsilon^\frac{1}{k}} \geq \delta {\mathbb V}_{\mathbb P}(f(Y)) >0.
\end{align*}
This contradicts \eqref{eq:contradiction1}.\hfill$\Box$

\section{Penalty version} \label{App:Penalty}
\subsection{Penalty version for smooth $\phi$-divergence}
In \cite{gotoh2018robust}, the penalty version of the worst-case problem is used to define worst-case sensitivity. Specifically, a family of worst-case distributions $\{\tilde{\mathsf q}(\delta)\,|\,\delta\geq 0\}$ is given by the solutions of the worst-case problem
\begin{align}
\tilde{\mathsf{q}}(\delta) & := \left\{
\begin{array}{cl}
\begin{displaystyle}\argmax_{\mathsf q}\Big\{  \sum_{i=1}^nq_i f_i -  \frac{1}{\delta} \sum_{i=1}^n {p}_i \phi\Big(\frac{q_i}{{p}_i}\Big)\Big\}, \end{displaystyle}& \delta>0,\\
\mathsf{p}, & \delta=0,
\end{array}\right.
\label{eq:worst-case-Q}
\end{align}
where the parameter $\delta$ determines the penalty on deviations from the nominal. In particular, $\delta=0$ gives the nominal and increasing $\delta$ is analogous to increasing the size of the uncertainty set in \eqref{eq:phi-theta}. When the penalty version is used to define the set of worst-case measures, the worst-case expected cost under 
$\tilde{\mathsf{q}}(\delta)$ is linear in the ambiguity parameter
\begin{align*}
V(\delta) & \equiv \mathbb{E}_{\tilde{\mathsf{q}}(\delta)}(\mathsf{f}) = \mathbb{E}_{\mathsf{p}}(\mathsf{f})
                +\frac{\delta}{\phi''(1)}\mathbb{V}_{\mathsf{p}}(\mathsf{f}) + o(\delta),
\end{align*}
so the standard definition of sensitivity \eqref{eq:sensitivity-general2} 
can be used:
\begin{align}
{\mathcal S}_{\mathsf{p}}(\mathsf{f})
 &=\lim_{\delta\downarrow 0}\frac{\mathbb{E}_{\tilde{\mathsf{q}}(\delta)}(\mathsf{f}) -\mathbb{E}_{\mathsf{p}}(\mathsf{f})}{\delta}
  =\frac{1}{\phi''(1)}\mathbb{V}_{\mathsf{p}}(\mathsf{f}).
\label{eq:phi-sensitivity-var}
\end{align}
While this leads to a different sensitivity measure, its qualitative nature is the same as \eqref{eq:phi-sensitivity}.

\section{Proofs from Section \ref{sec:WCS}}
\label{App:WCSproofs}

\subsection{WCS for $\phi$-divergence uncertainty sets}
\label{App:phi-div}

Since $p_i>0$ for all $i$, Assumption \ref{ass:phi} and
convex duality imply that
for small $\varepsilon>0$,
\begin{align*}
V_{\phi}(\varepsilon) &= \min_{\delta>0,\, c} \max_{\mathsf{q>0}} \sum_{i=1}^{n}q_i f
 + \frac{1}{\delta}\Big(\varepsilon-\sum_{i=1}^{n}p_i \phi\Big(\frac{q_i}{p_i}\Big)\Big)  + c\Big(\sum_{i=1}^{n}q_i - 1\Big).
\end{align*}
The following result characterizes the solution $(c(\varepsilon), \delta(\varepsilon))$ of the dual problem and the worst-case distribution.
\begin{proposition}
\label{prop:phi-div}
Suppose that $\phi$ satisfies Assumption \ref{ass:phi}. Then
\begin{align}
c(\varepsilon)  & =  -{\mathbb E}_{\mathsf{p}}(\mathsf{f}) + O(\sqrt{\varepsilon}), \label{eq:c_phi}
\\
\delta(\varepsilon) & =  \sqrt{\varepsilon}\sqrt{\frac{2 \phi''(1)}{{\mathbb V}_{\mathsf{p}}(\mathsf{f})}}+ o(\sqrt{\varepsilon}).
\label{eq:delta_phi}
\end{align}
The family of worst-case distributions $\{{\mathsf q}(\varepsilon)\,|\,\varepsilon\geq 0\}$ satisfies
\begin{align*}
q_i(\varepsilon) & = p_i\Big\{1 + \sqrt\frac{2\varepsilon}{{{\mathbb V}_{\mathsf{p}}(\mathsf{f})}} \left(f_i-\mathbb{E}_{\mathsf{p}}(\mathsf{f})\right)\Big\}+ o(\sqrt\varepsilon).
\end{align*}
\end{proposition}

\begin{proof}
\begin{align*}
V_{\phi}(\varepsilon)  = & \min_{\delta\geq 0,\, c} \; \frac{1}{\delta} \sum_{i=1}^{n} p_i \max_{q_i} \Big\{ \frac{q_i}{p_i}\delta \big(f_i+c\big) - \phi\Big(\frac{q_i}{p_i}\Big)\Big\} + \frac{\varepsilon}{\delta}-c \\
= &  \sum_{i=1}^n p_i f_i + \min_{\delta\geq 0,\, c} \frac{1}{\delta} \Big\{\sum_{i=1}^{n} p_i\Big[ \phi^*\Big(\delta \big(f_i+c\big)\Big) - \delta(f_i+c)\Big]+ {\varepsilon}\Big\}
\end{align*}
where
\begin{align*}
\phi^*(\zeta) = \max_z \big\{\zeta z - \phi(z)\big\}
\end{align*}
is the convex conjugate of $\phi(z)$ and
\begin{align}
q_i = p_i[\phi']^{-1}\Big(\delta(f_i+c)\Big) = \argmax_{q_i} \Big\{ \frac{q_i}{p_i}\delta \big(f_i+c\big) - \phi\Big(\frac{q_i}{p_i}\Big)\Big\}
\label{eq:wcq_phi_temp}
\end{align}
is the optimizer in the first equality. The worst-case measure ${\mathsf q}(\varepsilon) = \big(q_1(\varepsilon),\cdots,\,q_n(\varepsilon)\big)^\top$ is obtained by substituting the optimizers over $\delta$ and $c$.

Under the assumptions about $\phi(z)$ from \cite{gotoh2018robust}, $\phi^*(\zeta)$ is convex and twice continuously differentiable with
\begin{align}
\phi^*(\zeta)=\zeta + \frac{\zeta^2}{2\phi''(1)}+o(\zeta^2)
\label{eq:phi*-expansion}
\end{align}
Differentiating with respect to $c$ and $\delta$, the first order conditions are
\begin{align}
\sum_{i=1}^{n}p_i[\phi^*]'\Big(\delta[f_i+c]\Big)   &=   1 \label{eq:foc-c} \\
\sum_{i=1}^{n}p_i\Big\{\phi^*\Big(\delta[f_i+c]\Big) -  [\phi^*]'\Big(\delta[f_i+c]\Big)\delta[f_i+c]\Big\}+\varepsilon  &=  0. \label{eq:foc-delta}
\end{align}
Clearly $c(0) = -{\mathbb E}_{{\mathsf p}}(\mathsf{f})$. We can apply the Implicit Function Theorem to show that $c(\delta)$ is continuously differentiable in $\delta$ in the neighborhood of $\delta=0$ so it follows from \eqref{eq:phi*-expansion} and \eqref{eq:foc-c} that
\begin{align*}
c(\delta) = -{\mathbb E}_{\mathsf{p}}(\mathsf{f}) + O(\delta).
\end{align*}
Together with the expansion of $\phi^*(\zeta)$, \eqref{eq:foc-delta} becomes
\begin{align*}
\frac{\delta^2}{2 \phi^{''}(1)}\sum_{i=1}^{n} p_i\big(f_i-{\mathbb E}_{\mathsf{p}}(\mathsf{f})\big)^2 + o(\delta^2) = \varepsilon.
\end{align*}
The Implicit Function Theorem can again be used to show that $\delta(\varepsilon)$ is continuously differentiable on some open interval\footnote{Consider \begin{align*}\frac{y}{2 \phi^{''}(1)}\sum_{i=1}^{n} p_i\big(f_i-{\mathbb E}_{\mathsf{p}}(\mathsf{f})\big)^2 + o(y) = \varepsilon
\end{align*} Observe that $y=0$ when $\varepsilon=0$. The Implicit Function Theorem implies that $y(\varepsilon)$ is continuously differentiable on an open interval containing $\varepsilon=0$. Continuous differentiability of $\delta(\varepsilon)$ on some open interval $(0,\,b)$ follows from the observation that $\delta(\varepsilon)=\sqrt{y(\varepsilon)}$ for $\varepsilon>0$.} $(0,\,b)$ and given by \eqref{eq:delta_phi}.
It follows that  $c(\varepsilon)$ is given by \eqref{eq:c_phi}.
%
Finally, the worst-case distribution is obtained by substituting \eqref{eq:c_phi} and \eqref{eq:delta_phi} into \eqref{eq:wcq_phi_temp}.
\end{proof}

\subsection{WCS for total variation uncertainty set}
\label{App:TV}

\begin{proposition}\label{lemma:sen_tv1}
For $\varepsilon \in (0, \min(\mathsf{p}))$, 
a worst-case probability distribution of \eqref{def:wcobj:TV} is
\begin{equation}
(q_{(1)},q_{(2)},...,q_{(n-1)},q_{(n)})=
\big(p_{(1)}+\frac{\varepsilon}{2},~p_{(2)},~...,~p_{(n-1)},~p_{(n)}-\frac{\varepsilon}{2}\big)
\label{eq:worst-case_solution_to_tv}
\end{equation}
and the worst-case objective \eqref{def:wcobj:TV} is
\[
V_{\rm TV}(\varepsilon;\mathsf{f})=
\mathbb{E}_\mathsf{p}(\mathsf{f})+\varepsilon \cdot \frac{\max(\mathsf{f})-\min(\mathsf{f})}{2}.
\]
\end{proposition}


%
\begin{proof}
The worst-case objective \eqref{def:wcobj:TV} 
is written by the following optimization problem:
\[
V_{\rm TV}(\varepsilon;\mathsf{f})\equiv
\begin{array}[t]{|ll}
\underset{\mathsf{q}}{\mbox{maximize}}&\mathsf{f}^\top\mathsf{q}\\
\mbox{subject to}&\mathsf{1}^\top\mathsf{q}=1,\\
                 &\mathsf{1}^\top|\mathsf{p}-\mathsf{q}|\leq\varepsilon,
\end{array}
\]
where the nonnegativity condition ``$\mathsf{q}\geq\mathsf{0}$'' is omitted as we suppose that $p_i>0$ for all $i$ and  $\varepsilon>0$ is sufficiently small. 
In addition, let us introduce nonnegative vectors $\mathsf{u},\mathsf{v}\geq{0}$ such that $\mathsf{u}-\mathsf{v}=\mathsf{q}-\mathsf{p}$.
For 
small $\varepsilon$, the worst-case objective can be rewritten by
\[
\begin{array}{|ll}
\underset{\mathsf{u},\mathsf{v}}{\mbox{maximize}}&\mathsf{f}^\top\mathsf{u}-\mathsf{f}^\top\mathsf{v}+\mathsf{p}^\top\mathsf{f}\\
\mbox{subject to}&\mathsf{1}^\top\mathsf{u}-\mathsf{1}^\top\mathsf{v}=0,\\
                 &\mathsf{1}^\top\mathsf{u}+\mathsf{1}^\top\mathsf{v}\leq\varepsilon,\\
                 &\mathsf{u},\mathsf{v}\geq\mathsf{0}.
\end{array}
\]
The dual LP problem is derived as
\[
\begin{array}{|ll}
\underset{\lambda\geq 0,\theta}{\mbox{minimize}}&\varepsilon\lambda+\mathsf{p}^\top\mathsf{f}\\
\mbox{subject to}&\mathsf{1}\theta+\mathsf{1}\lambda\geq\mathsf{f},\\
                 &-\mathsf{1}\theta+\mathsf{1}\lambda\geq-\mathsf{f},\\
                 &\lambda\geq 0,
\end{array}
\]
which can be reduced to
\[
\begin{array}{ll}
\underset{\theta}{\mbox{minimize}}&\varepsilon\max\{|f_1-\theta|,...,|f_n-\theta|\}+\mathsf{p}^\top\mathsf{f}.
\end{array}
\]
It is easy to see that its optimality attained at $\theta=\frac{f_{(1)}+f_{(n)}}{2}$ with the optimal value being $\frac{\varepsilon(f_{(1)}-f_{(n)})}{2}+\mathsf{p}^\top\mathsf{f}$, and that the solution \eqref{eq:worst-case_solution_to_tv} attains this optimal value and is in $\mathcal{Q}_{\rm TV}(\varepsilon)$ for $\varepsilon\in(0,\min(\mathsf{p}))$. 
%
\end{proof}

\subsection{Proof of Propositions \ref{prop:cvar_explicit} and \ref{cor:sen_cvar} (WCS for budgeted uncertainty set)}\label{sec:proof:wcs_budgeted}
\paragraph{Proof of Proposition \ref{prop:cvar_explicit}}
Solving 
LP \eqref{eq:dual_cvar_eps}, which is equivalent to \eqref{eq:dual_cvar} under the change of parameter $\varepsilon=\alpha/(1-\alpha)$, by a greedy algorithm, 
an optimal solution to 
\eqref{eq:dual_cvar_eps} or \eqref{eq:dual_cvar} is derived as 
\begin{equation}
q_{(i)}=
\left\{
\begin{array}{ll}
\displaystyle (1+\varepsilon)p_{(i)}\equiv\frac{1}{1-\alpha}p_{(i)},&i=1,...,k,\\
\displaystyle 1-(1+\varepsilon)\sum_{i=1}^{k}p_{(i)}\equiv 1-\frac{1}{1-\alpha}\sum_{i=1}^{k}p_{(i)},&i=k+1,\\
\displaystyle 0,&i=k+2,...,n,
\end{array}
\right.
\label{eq:dual_cvar_sol}
\end{equation}
by the  strong duality of LP, the expression of $V_{\rm b}(\varepsilon;\mathsf{f})$ in the statement of Proposition \ref{prop:cvar_explicit} is obtained as the optimal objective function value, and 
CVaR for \eqref{eq:dual_cvar} is expressed 
as 
\begin{align*}
\mathrm{CVaR}_{\mathsf{p},\alpha}(\mathsf{f})=\frac{1}{1-\alpha}\sum_{i=1}^{k}p_{(i)}f_{(i)}+\big(1-\frac{1}{1-\alpha}\sum_{i=1}^{k}p_{(i)}\big)f_{(k+1)},
\end{align*}
which is equivalent to $V_{\rm b}(\alpha/(1-\alpha);\mathsf{f})$.\hfill (End of Proof of Proposition \ref{prop:cvar_explicit})

It is easy to see that the worst-case expected cost $V_{\rm b}(\varepsilon)$ is piecewise linear, concave, and increasing in $\varepsilon$. 
The following result characterizes the slope of $V_{\rm b}(\varepsilon)$ for all values of $\varepsilon$.
\begin{lemma}
\label{prop:sen_cvar_env}
Let $\varepsilon>0$. Suppose $k\in \{0,...,n-1\}$ is an integer such that
\begin{equation}
\varepsilon \in 
\Big(\frac{\sum_{i=k+2}^{n}p_{(i)}}{\sum_{i=1}^{k+1}p_{(i)}},\frac{\sum_{i=k+1}^{n}p_{(i)}}{\sum_{i=1}^{k}p_{(i)}}\Big),
\label{eq:intervals_cvar}
\end{equation}
where 
$\sum\limits_{i=n+1}^np_{(i)}=\sum\limits_{i=1}^0p_{(i)}=0$ and $\frac{1}{0}=\infty$.
For $\Delta>0$ satisfying $\varepsilon+\Delta<\frac{\sum_{i=k+1}^{n}p_{(i)}}{\sum_{i=1}^{k}p_{(i)}}$, we have
\begin{align}
S_{\rm b}(\varepsilon):=
\frac{V_{\rm b}(\varepsilon+\Delta)-V_{\rm b}(\varepsilon)}{\Delta}
&=\sum_{i=1}^kp_{(i)}\big(f_{(i)}-f_{(k+1)}\big)\label{eq:constant_slope_cvar}
\\
&=\frac{1}{1+\varepsilon}\Big(
\mathrm{CVaR}_{\mathsf{p},\frac{\varepsilon}{1+\varepsilon}}(\mathsf{f})-\mathrm{VaR}_{\mathsf{p},\frac{\varepsilon}{1+\varepsilon}}(\mathsf{f})
\Big).
\label{eq:sen_cvar_env}
\end{align}
\end{lemma}
Expression \eqref{eq:constant_slope_cvar} shows that the slope in the intetaval \eqref{eq:intervals_cvar} is constant (i.e., independent of $\varepsilon$), but expression \eqref{eq:sen_cvar_env} provides another insightful interpretation of the slope.
\paragraph{Proof of Lemma \ref{prop:sen_cvar_env}}
With the formula of $V_{\rm b}(\varepsilon)$ in the statement of Proposition \ref{prop:cvar_explicit}, we have
\begin{align}
\lefteqn{V_{\rm b}(\varepsilon+\Delta)-V_{\rm b}(\varepsilon)} \nonumber \\
& =  \mathrm{CVaR}_{\mathsf{p},\frac{\varepsilon+\Delta}{1+\varepsilon+\Delta}}(\mathsf{f})-\mathrm{CVaR}_{\mathsf{p},\frac{\varepsilon}{1+\varepsilon}}(\mathsf{f})\nonumber\\
&= (1+\varepsilon+\Delta)\sum_{i=1}^{k}p_{(i)}f_{(i)}+\Big\{1-(1+\varepsilon+\Delta)\sum_{i=1}^{k}p_{(i)}\Big\}f_{(k+1)}\nonumber\\
&\qquad -(1+\varepsilon)\sum_{i=1}^{k}p_{(i)}f_{(i)}-\Big\{1-(1+\varepsilon)\sum_{i=1}^{k}p_{(i)}\Big\}f_{(k+1)}\nonumber\\
&=\frac{\Delta}{1+\varepsilon}\Big[
\underbrace{
(1+\varepsilon)\sum_{i=1}^{k}p_{(i)}f_{(i)}+\Big\{1-(1+\varepsilon)\sum_{i=1}^{k}p_{(i)}\Big\}f_{(k+1)}
}_{V_{\rm b}(\varepsilon)\equiv\mathrm{CVaR}_{\mathsf{p},\frac{\varepsilon}{1+\varepsilon}}(\mathsf{f})}-\underbrace{f_{(k+1)}}_{\mathrm{VaR}_{\mathsf{p},\frac{\varepsilon}{1+\varepsilon}}(\mathsf{f})}\Big]\label{eq:cvar-cvar}\\
&=\frac{\Delta}{1+\varepsilon}\Big(
\mathrm{CVaR}_{\mathsf{p},\frac{\varepsilon}{1+\varepsilon}}(\mathsf{f})-\mathrm{VaR}_{\mathsf{p},\frac{\varepsilon}{1+\varepsilon}}(\mathsf{f})
\Big),\nonumber
\end{align}
which proves \eqref{eq:sen_cvar_env}. Note that \eqref{eq:constant_slope_cvar} follows from \eqref{eq:cvar-cvar}.
(End of Proof of Lemma \ref{prop:sen_cvar_env})
\paragraph{Proof of Proposition  \ref{cor:sen_cvar}}
\eqref{eq:constant_slope_cvar} defines the constant slope of the linear piece over the interval \eqref{eq:intervals_cvar}. 
Since $V_{\rm b}(\varepsilon)$ is continuous as well as concave and increasing for $\varepsilon\geq 0$, its slope is the largest over the left-most interval $[0,\frac{p_{(n)}}{1-p_{(n)}}]$. For $\varepsilon\in(0,\frac{p_{(n)}}{1-p_{(n)}})$ and $\Delta$ sufficiently small, \eqref{eq:sen_cvar_env} becomes
\begin{align*}
\mathcal{S}_{\mathsf{p}}(\mathsf{f})
&=\lim_{\varepsilon'\downarrow 0}\frac{V_{\rm b}(\varepsilon')-V_{\rm b}(0)}{\varepsilon'}\\
&=
\frac{V_{\rm b}(0+\Delta)-V_{\rm b}(0)}{\Delta} \\
&=
\frac{V_{\rm b}(\varepsilon+\Delta)-V_{\rm b}(\varepsilon)}{\Delta} \\
&\underset{\eqref{eq:constant_slope_cvar}}{=}
\sum_{i=1}^{n-1}p_{(i)}(f_{(i)}-f_{(n)})+p_{(n)}(f_{(n)}-f_{(n)})\\
&=\sum_{i=1}^{n}p_{(i)}f_{(i)}-f_{(n)}
\\
&= {\mathbb E}_{\mathsf{p}}(\mathsf{f}) - \min(\mathsf{f}),
\end{align*}
which is the desired result. \hfill$\Box$

More specifically, by using \eqref{eq:sen_cvar_env}, we can obtain an explicit piecewise linear decomposition of $V_{\rm b}(\varepsilon)$. 
Let
\begin{eqnarray*}\varepsilon_{(h)}=\frac{\sum_{i=n-h+1}^{n}p_{(i)}}{\sum_{i=1}^{n-h}p_{(i)}},~h=0,1,2,...,n-1.
\end{eqnarray*}
 For $\varepsilon\in(\varepsilon_{(j)},\varepsilon_{(j+1)}]$, the $(j+1)^{th}$ interval, 
\begin{align*}
V_{\rm b}(\varepsilon)
&=V_{\rm b}(0)+\sum_{h=0}^{j-1}S_{\rm b}(\varepsilon_h)(\varepsilon_{(h+1)}-\varepsilon_{(h)})+S_{\rm b}(\varepsilon_j)(\varepsilon-\varepsilon_{(j)})\\
&=V_{\rm b}(0)+\sum_{h=0}^{j-1}\frac{
1}{1+\varepsilon_{(h)}}\big(\mathrm{CVaR}_{\mathsf{p},\frac{\varepsilon_{(h)}}{1+\varepsilon_{(h)}}}(\mathsf{f})-\mathrm{VaR}_{\mathsf{p},\frac{\varepsilon_{(h)}}{1+\varepsilon_{(h)}}}(\mathsf{f})\big)(\varepsilon_{(h+1)}-\varepsilon_{(h)})\\
&\qquad\qquad\qquad +\frac{
1}{1+\varepsilon_{(j)}}\big(\mathrm{CVaR}_{\mathsf{p},\frac{\varepsilon_{(j)}}{1+\varepsilon_{(j)}}}(\mathsf{f})-\mathrm{VaR}_{\mathsf{p},\frac{\varepsilon_{(j)}}{1+\varepsilon_{(j)}}}(\mathsf{f})\big)(\varepsilon-\varepsilon_{(j)}).
\end{align*}
%


\subsection{Proof of Proposition \ref{cor:wcs:comb} (WCS for convex combination of expected cost and $\alpha$-CVaR)}\label{sec:proof:convcomb}
With the change of variables $\mathsf{Q}=\frac{1}{\varepsilon}\{\mathsf{q}-(1-\varepsilon)\mathsf{p}\}$, we have 
 \begin{align*}
 V_{\rm c}(\varepsilon;\mathsf{f})&=
 \max_{\mathsf{q}} \Big\{ \mathsf{f}^\top\mathsf{q} \Big\vert 
 \bm{1}^\top\mathsf{q}=1, (1-\varepsilon)\mathsf{p}\leq\mathsf{q}\leq(1-\varepsilon)\mathsf{p}+\frac{\varepsilon}{1-\alpha}\mathsf{p} \Big\} \\
 &= (1-\varepsilon)\mathsf{f}^\top\mathsf{p}+\varepsilon \underbrace{\max_{\mathsf{Q}} \Big\{ \mathsf{f}^\top\mathsf{Q} \Big\vert 
 \bm{1}^\top\mathsf{Q}=1, \mathsf{0}\leq \mathsf{Q}\leq\frac{1}{1-\alpha}\mathsf{p} \Big\}}_{\mathrm{CVaR}_{\mathsf{p},\alpha}(\mathsf{f})},
 \end{align*}
which equals the left-hand side of \eqref{eq:sen_cvxcmb_expectation+cvar}. 
Similar to \eqref{eq:dual_cvar_sol}, an optimal solution, $\mathsf{Q}=\mathsf{Q}^*$, of the last maximization can be obtained, 
and the worst-case probability is obtained via the formula $\mathsf{q}^*=\varepsilon\mathsf{Q}^*+(1-\varepsilon)\mathsf{p}$.
\hfill$\Box$

\subsection{Proof of Proposition \ref{prop:rcvar_sensitivities} (WCSs for CVaR)}\label{sec:proof:wcs_cvar}
Note first that under Assumption \ref{ass:nondegeneracy_for_cvar}, optimal solutions 
for the CVaR formula are uniquely determined and $\gamma=f_{(k+1)}=\mathrm{VaR}_{\mathsf{p},\alpha}(\mathsf{f})$.\\ 
\noindent
(a) With $\mathcal{Q}(\varepsilon)=\mathcal{Q}_{\phi}(\varepsilon)$, the robust CVaR objective turns out to be
\begin{align}
\mathrm{RCVaR}_{\mathsf{p},\beta}^{\phi,\varepsilon}(\mathsf{f})\equiv
\begin{array}[t]{|lll}
\underset{\mathsf{q},\mathsf{Q}}{\mbox{maximize}}&\mathsf{f}^\top\mathsf{Q}\\
\mbox{subject to}
&\mathsf{1}^\top\mathsf{Q}=1,&\gets\gamma\\
&\mathsf{1}^\top\mathsf{q}=1,&\gets\nu\\
&\mathsf{Q}-\frac{1}{1-\beta}\mathsf{q}\leq\mathsf{0},&\gets\mathsf{\lambda}\geq\mathsf{0}\\
&\sum\limits_{i=1}^{n}p_i\phi(\frac{q_i}{p_i})\leq\varepsilon,&\gets\kappa\geq 0\\
&\mathsf{Q}\geq\mathsf{0},&\gets\mathsf{\theta}\geq\mathsf{0},
\end{array}
\label{eq:rcvar_w/phi_dual}
\end{align}
where the non-negativity condition ``$\mathsf{q}\geq\mathsf{0}$'' is omitted since it is implied by the first and third inequalities. 
The Lagrangian is 
\begin{align*}
L(\mathsf{Q},\mathsf{q},\gamma,\nu,\mathsf{\lambda},\kappa,\mathsf{\theta})
=(\mathsf{f}-\gamma\mathsf{1}-\mathsf{\lambda}+\mathsf{\theta})^\top\mathsf{Q}
+\kappa\sum_{i=1}^np_i\Big[\frac{q_i}{p_i}\Big(\frac{\lambda_i}{1-\beta}-\nu\Big)\cdot\frac{1}{\kappa}-\phi\big(\frac{q_i}{p_i}\big)\Big]
+\gamma+\nu+\varepsilon\kappa,
\end{align*}
and the Lagrangian dual becomes
\begin{align*}
\begin{array}{|ll}
\underset{\gamma,\nu,\mathsf{\lambda},\kappa,\mathsf{\theta}}{\mbox{minimize}}~ & \displaystyle
\gamma+\nu+\varepsilon\kappa+
\kappa\sum_{i=1}^np_i\phi^*\Big(\frac{1}{\kappa}\big(\frac{\lambda_i}{1-\beta}-\nu\big)\Big)
,\\
\mbox{subject to}
\quad&\mathsf{\lambda}\geq\mathsf{f}-\gamma\mathsf{1},~\mathsf{\lambda}\geq\mathsf{0},~\kappa>0 
\end{array}
\end{align*}
With the formula $\phi^*(\zeta)=\zeta+\frac{1}{2\phi''(1)}\zeta^2+o(\zeta^2)$, the dual turns out to be
\begin{align*}
\begin{array}{|ll}
\underset{\gamma,\nu,\mathsf{\lambda},\kappa}{\mbox{minimize}}~ & \displaystyle
\gamma
+\frac{1}{1-\beta}\sum_{i=1}^np_i\lambda_i
+\varepsilon\kappa
+\frac{1}{2\phi''(1)\kappa}\sum_{i=1}^np_i\big(\frac{\lambda_i}{1-\beta}-\nu\big)^2
+\kappa\sum_{i=1}^np_io\big(\frac{\big(\frac{\lambda_i}{1-\beta}-\nu\big)^2}{\kappa^2}\big)
,\\
\mbox{subject to}
~&\mathsf{\lambda}\geq\mathsf{f}-\gamma\mathsf{1},~\mathsf{\lambda}\geq\mathsf{0},~\kappa>0 
\end{array}
\end{align*}
The arithmetic mean--geometric mean inequality implies that
\begin{align*}
\varepsilon\kappa+
\frac{1}{2\phi''(1)\kappa}\sum_{i=1}^np_i\big(\frac{\lambda_i}{1-\beta}-\nu\big)^2
\geq
\sqrt{\frac{2\varepsilon}{\phi''(1)}\sum_{i=1}^np_i\big(\frac{\lambda_i}{1-\beta}-\nu\big)^2}
\end{align*}
where the equality holds when $\kappa=\sqrt{\frac{\sum_{i=1}^np_i\big(\frac{\lambda_i}{1-\beta}-\nu\big)^2}{2\phi''(1)\varepsilon}}$. 
Using some perturbation arguments such as $\gamma=\gamma_0+\varepsilon\gamma_1+o(\varepsilon),\lambda_i=\lambda_{i,0}+\varepsilon\lambda_{i,1}+o(\varepsilon),\nu=\nu_0+\varepsilon\nu_1+o(\varepsilon)$ with $\gamma_0,\lambda_{i,0},\nu_0$ being the elements of the optimal solution for $\varepsilon=0$, i.e., $\gamma_0=\mathrm{VaR}_{\mathsf{p},\beta}(\mathsf{f})$, $\lambda_{i0}=\max\{f_i-\gamma_0,0\}$, and $\nu_0=\frac{1}{1-\beta}\sum_{i=1}^np_i\lambda_{i,0}=\mathrm{CVaR}_{\mathsf{p},\beta}(\mathsf{f})-\mathrm{VaR}_{\mathsf{p},\beta}(\mathsf{f})$
, we have $\kappa=O(\sqrt{\varepsilon})$ and 
$\kappa\sum_{i=1}^np_io\big(\frac{\big(\frac{\lambda_i}{1-\beta}-\nu\big)^2}{\kappa^2}\big)=
O(\varepsilon^{1/2})o(\varepsilon)$ and  
\begin{align*}
&
\begin{array}[t]{|ll}
\underset{\gamma,\nu,\mathsf{\lambda}}{\mbox{minimize}}~& \displaystyle
\underbrace{\gamma+
\frac{1}{1-\beta}\sum_{i=1}^np_i\lambda_i}_{
\mathrm{CVaR}_{\mathsf{p},\beta}(\mathsf{f})+O(\varepsilon)}
+\sqrt{\frac{2\varepsilon}{\phi''(1)}}
\hspace{-1em}
\underbrace{\displaystyle\sqrt{
\sum\limits_{i=1}^np_i\big(\frac{\lambda_i}{1-\beta}-\nu\big)^2}}_{
\sqrt{\sum\limits_{i=1}^np_i\big(\frac{\lambda_{i,0}}{1-\beta}-\nu_{0}\big)^2+O(\varepsilon)+O(\varepsilon^2)}
}
+o(\sqrt{\varepsilon})
,\\
\mbox{subject to}
~&\mathsf{\lambda}\geq\mathsf{f}-\gamma\mathsf{1},~\mathsf{\lambda}\geq\mathsf{0}
\end{array}\\
=&
\begin{array}[t]{ll}
\mathrm{CVaR}_{\mathsf{p},\beta}(\mathsf{f})
+\sqrt{\frac{2\varepsilon}{\phi''(1)}}\sqrt{\sum\limits_{i=1}^np_i\big\{\frac{|f_i-\mathrm{VaR}_{\mathsf{p},\beta}(\mathsf{f})|_+}{1-\beta}-\big(\mathrm{CVaR}_{\mathsf{p},\beta}(\mathsf{f})-\mathrm{VaR}_{\mathsf{p},\beta}(\mathsf{f})\big)\big\}^2
+O(\varepsilon^2)}
+o(\sqrt{\varepsilon})
,\end{array}
\end{align*}
The desired result follows. \hfill$\Box$

\noindent
(b)
With $\mathcal{Q}(\varepsilon)=\mathcal{Q}_{\rm TV}(\varepsilon)$, the robust CVaR objective turns out to be
\begin{align*}
\mathrm{RCVaR}_{\mathsf{p},\beta}^{\mathrm{TV},\varepsilon}(\mathsf{f})&\equiv
\begin{array}[t]{|lll}
\underset{\mathsf{q},\mathsf{Q},\mathsf{u},\mathsf{v}}{\mbox{maximize}}&\mathsf{f}^\top\mathsf{Q}\\
\mbox{subject to}
&\mathsf{1}^\top\mathsf{Q}=1,\\
&\mathsf{1}^\top\mathsf{q}=1,\\
&\mathsf{Q}-\frac{1}{1-\beta}\mathsf{q}\leq\mathsf{0},\\
&\mathsf{1}^\top(\mathsf{u}+\mathsf{v})\leq\varepsilon,\\
&\mathsf{p}-\mathsf{q}=\mathsf{u}-\mathsf{v},\\
&\mathsf{Q},\mathsf{q},\mathsf{u},\mathsf{v}\geq\mathsf{0}
\end{array}\\
&=
\begin{array}[t]{|lll}
\underset{\mathsf{Q},\mathsf{u},\mathsf{v}}{\mbox{maximize}}&\mathsf{f}^\top\mathsf{Q}\\
\mbox{subject to}
&\mathsf{1}^\top\mathsf{Q}=1,&\gets\gamma\\
&\mathsf{1}^\top\mathsf{u}-\mathsf{1}^\top\mathsf{v}=0,&\gets\nu\\
&\mathsf{Q}+\frac{1}{1-\beta}\mathsf{u}-\frac{1}{1-\beta}\mathsf{v}\leq\frac{1}{1-\beta}\mathsf{p},&\gets\mathsf{\lambda}\geq\mathsf{0}\\
&\mathsf{1}^\top\mathsf{u}+\mathsf{1}^\top\mathsf{v}\leq\varepsilon,&\gets\kappa\geq 0\\
&\mathsf{Q},\mathsf{u},\mathsf{v}\geq\mathsf{0}
\end{array}\\
&=
\begin{array}[t]{|lll}
\underset{\gamma,\nu,\mathsf{\lambda},\kappa}{\mbox{minimize}}&
\gamma+\frac{1}{1-\beta}\mathsf{p}^\top\mathsf{\lambda}+\varepsilon\kappa\\
\mbox{subject to}
&\mathsf{1}\gamma+\mathsf{\lambda}\geq\mathsf{f},\\
&\mathsf{1}\nu+\frac{1}{1-\beta}\mathsf{\lambda}+\kappa\mathsf{1}\geq\mathsf{0},\\
&-\mathsf{1}\nu-\frac{1}{1-\beta}\mathsf{\lambda}+\kappa\mathsf{1}\geq\mathsf{0},\\
&\mathsf{\lambda}\geq\mathsf{0},~\kappa\geq 0
\end{array}\\
&=
\begin{array}[t]{lll}
\underset{\gamma,\nu}{\mbox{minimize}}&
\gamma+\frac{1}{1-\beta}\mathsf{p}^\top|\mathsf{f}-\mathsf{1}\gamma|_++\varepsilon\max(\left|\mathsf{1}\nu-\frac{1}{1-\beta}|\mathsf{f}-\gamma\mathsf{1}|_+\right|),
\end{array}\\
&=
\begin{array}[t]{lll}
\underset{\gamma,\nu'}{\mbox{minimize}}&
\gamma+\frac{1}{1-\beta}\mathsf{p}^\top|\mathsf{f}-\gamma\mathsf{1}|_++\frac{\varepsilon}{1-\beta}\max(\big|\mathsf{1}\nu'-|\mathsf{f}-\mathsf{1}\gamma|_+\big|),
\end{array}\\
&=
\begin{array}[t]{lll}
&
\mathrm{VaR}_{\mathsf{p},\beta}(\mathsf{f})+\frac{1}{1-\beta}\mathsf{p}^\top|\mathsf{f}-\mathrm{VaR}_{\mathsf{p},\beta}(\mathsf{f})\mathsf{1}|_++\frac{\varepsilon}{2(1-\beta)}
\big(\max(\mathsf{f})-\mathrm{VaR}_{\mathsf{p},\beta}(\mathsf{f})\big)+O(\varepsilon),
\end{array}
\end{align*}
where the last equality comes from the sensitivity analysis of LP under Assumption \ref{ass:nondegeneracy_for_cvar}, and from the fact that
\begin{center}
$\begin{array}{l}
\displaystyle\arg\min_{\nu'}\{\big|\nu'-|f_{(1)}-f_{(k+1)}|\big|,...,\big|\nu'-|f_{(k)}-f_{(k+1)}|\big|,\big|\nu'-0\big|,...,\big|\nu'-0\big|\}\\
\displaystyle\quad=\arg\min_{\nu'}\{\big|\nu'-|f_{(1)}-f_{(k+1)}|\big|,\big|\nu'-0\big|\}\\
\displaystyle\quad=\frac{|f_{(1)}-f_{(k+1)}|+0}{2}=\frac{\max(\mathsf{f})-\mathrm{VaR}_{\mathsf{p},\beta}(\mathsf{f})}{2}.
\end{array}$
\end{center}
Consequently, we reach the desired result.\hfill$\Box$

\noindent
(c)
With $\mathcal{Q}(\varepsilon)=\mathcal{Q}_{\rm b}(\varepsilon)$, the robust CVaR objective turns out to be
\begin{align}
\mathrm{RCVaR}_{\mathsf{p},\beta}^{{\rm b},\varepsilon}(\mathsf{f})\equiv
\begin{array}[t]{|ll}
\underset{\mathsf{q},\mathsf{Q}}{\mbox{maximize}} & \mathsf{f}^\top\mathsf{Q}\\
\mbox{subject to}                              & \mathsf{1}^\top\mathsf{Q}=1,\\
                                               & \mathsf{1}^\top\mathsf{q}=1,\\
                                               & \mathsf{0}\leq\mathsf{Q}\leq\frac{1}{1-\beta}\mathsf{q},\\
                                               & \mathsf{0}\leq\mathsf{q}\leq(1+\varepsilon)\mathsf{p}.\\
\end{array}
\label{eq:rcvar_lp}
\end{align}
Observe first that when $\varepsilon=0$, the optimal solution to \eqref{eq:rcvar_lp} is given by 
\[
(\mathsf{q},\mathsf{Q})=(\mathsf{p},\mathsf{q}^*),
\]
where $\mathsf{q}^*$ is given by \eqref{eq:dual_cvar_sol}. 

Next consider a small perturbation $\varepsilon>0$. 
 From the theory of LP (or simplex method), we see that if $\varepsilon$ is sufficiently small, the basis variables can remain positive, and due to the greedy algorithm, the solution will be of the form
\[
(\mathsf{q},\mathsf{Q})=(\hat{\mathsf{q}}(\varepsilon),\hat{\mathsf{Q}}(\varepsilon)),
\]
where
\begin{align*}
(\hat{q}_{(i)}(\varepsilon),\hat{Q}_{(i)}({\varepsilon}))&:=
\left\{
\begin{array}{ll}
\displaystyle((1+\varepsilon)p_{(i)},\frac{1+\varepsilon}{1-\beta}p_{(i)}),&i=1,...,k,\\
\displaystyle (1-(1+\varepsilon)\sum_{i=1}^{k}p_{(i)},1-\frac{1+\varepsilon}{1-\beta}\sum_{i=1}^{k}p_{(i)}),&i=k+1,\\
\displaystyle (0,0)&i=k+2,...,n.
\end{array}
\right.
\end{align*}
%
Then we have (for small $\varepsilon$)
\begin{align*}
\mathrm{RCVaR}_{\mathsf{p},\beta}^{{\rm b},\varepsilon}(\mathsf{f})
 &=\frac{1+\varepsilon}{1-\beta}\sum_{i=1}^{k}p_{(i)}f_{(i)}+\big(1-\frac{1+\varepsilon}{1-\beta}\sum_{i=1}^{k}p_{(i)}\big)f_{(k+1)}\\
 &=\frac{1}{1-\beta}\sum_{i=1}^{k}p_{(i)}f_{(i)}+\big(1-\frac{1}{1-\beta}\sum_{i=1}^{k}p_{(i)}\big)f_{(k+1)}\\
 &\qquad+\varepsilon\Big[
   \frac{1}{1-\beta}\sum_{i=1}^{k}p_{(i)}f_{(i)}-\frac{1}{1-\beta}\sum_{i=1}^{k}p_{(i)}f_{(k+1)}
  \Big]\\
 &=\underbrace{\frac{1}{1-\beta}\sum_{i=1}^{k}p_{(i)}f_{(i)}+\big(1-\frac{1}{1-\beta}\sum_{i=1}^{k}p_{(i)}\big)f_{(k+1)}}_{\mathrm{CVaR}_{\mathsf{p},\beta}(\mathsf{f})}\\
 &\qquad+\varepsilon\Big[
   \underbrace{\frac{1}{1-\beta}\sum_{i=1}^{k}p_{(i)}f_{(i)}+\big(1-\frac{1}{1-\beta}\sum_{i=1}^{k}p_{(i)}\big)f_{(k+1)}}_{\mathrm{CVaR}_{\mathsf{p},\beta}(\mathsf{f})} - \underbrace{f_{(k+1)}}_{\mathrm{VaR}_{\mathsf{p},\beta}(\mathsf{f})}
  \Big]
\end{align*}
Therefore, the desired result follows. \hfill$\Box$

\noindent
(d)
With $\mathcal{Q}(\varepsilon)=\mathcal{Q}_{\rm c}(\varepsilon)$, the robust CVaR objective turns out to be
\begin{align*}
\lefteqn{\mathrm{RCVaR}_{\mathsf{p},\beta}^{\mathrm{c},\alpha,\varepsilon}(\mathsf{f})}\\
&\equiv
\begin{array}[t]{|lll}
\underset{\mathsf{q},\mathsf{Q}}{\mbox{maximize}}&\mathsf{f}^\top\mathsf{Q}\\
\mbox{subject to}
&\mathsf{1}^\top\mathsf{Q}=1,&\gets \gamma\\
&\mathsf{1}^\top\mathsf{q}=1,&\gets\nu\\
&\mathsf{Q}-\frac{1}{1-\beta}\mathsf{q}\leq\mathsf{0},&\gets\mathsf{\lambda}\geq\mathsf{0}\\
&-\mathsf{q}\leq-(1-\varepsilon)\mathsf{p},&\gets\mathsf{\pi}\geq\mathsf{0}\\
&\mathsf{q}\leq(1-\varepsilon+\frac{\varepsilon}{1-\alpha})\mathsf{p}&\gets\mathsf{\theta}\geq\mathsf{0}\\
&\mathsf{Q}\geq\mathsf{0},
\end{array}\\
&=
\begin{array}[t]{|lll}
\underset{\gamma,\nu,\mathsf{\lambda},\mathsf{\pi},\mathsf{\theta}}{\mbox{minimize}}&
\multicolumn{2}{l}{\gamma+\nu-(1-\beta)\mathsf{p}^\top\mathsf{\pi}+(1-\varepsilon+\frac{\varepsilon}{1-\alpha})\mathsf{p}^\top\mathsf{\theta}}\\
\mbox{subject to}
&\mathsf{1}\gamma+\mathsf{\lambda}\geq\mathsf{f},&\gets\mathsf{Q}\geq\mathsf{0}\\
&\mathsf{1}\nu-\frac{1}{1-\beta}\mathsf{\lambda}-\mathsf{\pi}+\mathsf{\theta}=\mathsf{0},&\gets\mathsf{q}\\
&\mathsf{\lambda},\mathsf{\pi},\mathsf{\theta}\geq\mathsf{0}
\end{array}\\
&=
\begin{array}[t]{|lll}
\underset{\gamma,\nu,\mathsf{\lambda},\mathsf{\theta}}{\mbox{minimize}}&
\gamma+\varepsilon\nu+\frac{1-\varepsilon}{1-\beta}\mathsf{p}^\top\mathsf{\lambda}+\frac{\varepsilon}{1-\alpha}\mathsf{p}^\top\mathsf{\theta}\\
\mbox{subject to}
&\mathsf{\lambda}\geq\mathsf{f}-\mathsf{1}\gamma,\\
&\mathsf{\theta}\geq-\mathsf{1}\nu+\frac{1}{1-\beta}\mathsf{\lambda},\\
&\mathsf{\lambda},\mathsf{\theta}\geq\mathsf{0}
\end{array}\\
&=
\begin{array}[t]{lll}
\underset{\gamma,\nu}{\min}&
\underbrace{\gamma+\frac{1}{1-\beta}\mathsf{p}^\top\big|\mathsf{f}-\mathsf{1}\gamma\big|_+}_{\mathrm{CVaR}_{\mathsf{p},\beta}(\mathsf{f})+O(\varepsilon)}
+
\varepsilon
\underbrace{
\Big\{
\nu+\frac{1}{1-\alpha}\mathsf{p}^\top\big|\frac{1}{1-\beta}|\mathsf{f}-\mathsf{1}\gamma|_+-\mathsf{1}\nu\big|_+
-
\frac{1}{1-\beta}\mathsf{p}^\top\big|\mathsf{f}-\mathsf{1}\gamma\big|_+
\Big\}}_{\mathrm{CVaR}_{\mathsf{p},\alpha}(\frac{1}{1-\beta}|\mathsf{f}-\mathsf{1}\gamma|_+)-\mathbb{E}_{\mathsf{p}}(\frac{1}{1-\beta}|\mathsf{f}-\mathsf{1}\gamma|_+)+O(\varepsilon)\mbox{ where }\gamma=\mathrm{VaR}_{\mathsf{p},\beta}(\mathsf{f})+O(\varepsilon)},
\end{array}
\end{align*}
so we reach the desired result. \hfill$\Box$


\subsection{Wasserstein uncertainty sets}\label{sec:proof:wcs_wasserstein}

\subsubsection*{Preliminaries: Some useful results from convex duality}

We summarize general properties of the dual problem and the relationship between optimal dual variables and super-gradients of the value function for the optimization problem
\begin{eqnarray}
V(\varepsilon) := \max_{x\in{\Omega}} F(x) \; \mbox{subject to:}\; G(x) \leq \varepsilon
\label{eq:convex-general}
\end{eqnarray}
which we apply to \eqref{eq:W2}.
Here, we assume that $F: X \rightarrow \mathbb R$ and $G:X \rightarrow {\mathbb R}$, where $X$ is a vector space and $\Omega$ is a convex subset of $\mathcal X$.
The associated dual problem is
\begin{eqnarray}
D(\varepsilon) := \min_{\lambda\geq 0} \max_{x\in{\mathcal X}} F(x) + \lambda\big\{\varepsilon-G(x)\big\}\equiv 
\min_{\lambda \geq 0} \Big\{ H(\lambda)+\lambda \varepsilon \Big\}
\label{eq:dual-general}
\end{eqnarray}
where $\lambda \in {\mathbb R}$ is the Lagrange multiplier and
\begin{eqnarray*}
H(\lambda) := \max_{x\in{\Omega}} \Big\{F(x)-\lambda G(x)\Big\}.
\end{eqnarray*}
Note that $H(\lambda)$ is convex in $\lambda$. 
We denote the derivative of $H$ at $\lambda$ by $H'(\lambda)$, and the directional derivaties
\begin{align*}
H'(\lambda^+) & := \lim_{\delta\downarrow 0}\frac{H(\lambda+\delta)-H(\lambda)}{\delta},\\
H'(\lambda^-) & := \lim_{\delta\downarrow 0}\frac{H(\lambda-\delta)-H(\lambda)}{\delta}.
\end{align*}
(We refer to $H'(\lambda^+)$ and $H'(\lambda^-)$ as the right and left derivative of $H$ at $\lambda$, respectively). $H(\lambda)$ is convex in $\lambda$ so $H'(\lambda^+)$ is increasing in $\lambda$ and $H'(\lambda^-)$ is decreasing in $\lambda$.

\begin{lemma} \label{lemma:duality_prop}
Consider the optimization problem \eqref{eq:convex-general}. Assume that $F(x)$ is concave, that $G(x)$ is convex and non-negative, that there is an $x\in\Omega$ such that $G(x) = 0$, and that $V(\varepsilon)$ is finite for every $\varepsilon$.
Then $V(\varepsilon)$ is concave and increasing in $\varepsilon \geq 0$, and differentiable at almost every $\varepsilon > 0$. When $\varepsilon>0$, strong duality holds $(V(\varepsilon)=D(\varepsilon))$  and the maximum of the dual problem is achieved. If $\lambda(\varepsilon)$ is a solution of the dual problem corresponding to $\varepsilon>0$, then $\lambda(\varepsilon)$ is a super-gradient of $V$ at $\varepsilon$. If $V$ is differentiable at $\varepsilon>0$, then the dual problem has a unique solution and  $V'(\varepsilon) = \lambda(\varepsilon)$.
\end{lemma}
\paragraph{Proof of Lemma \ref{lemma:duality_prop}}
When $\varepsilon>0$ the Lagrange Duality Theorem \cite{luenberger1997optimization} implies that strong duality holds and that there exists a solution $\lambda(\varepsilon)$ of the dual problem.
\hfill$\Box$


\begin{proposition} \label{prop:LM_limit}
Assume that $F(x)$ is concave, that $G(x)$ is convex and non-negative, that there is an $x\in\Omega$ such that $G(x) = 0$, and that $V(\varepsilon)$ is finite for every $\varepsilon$.
Let $\{\varepsilon_i\}$ be a sequence of positive numbers such that $\varepsilon_i\downarrow \varepsilon$. Suppose that $V$ is differentiable at $\varepsilon_i$, strong duality holds at $\varepsilon_i$, and there is a solution $\lambda_i\equiv\lambda(\varepsilon_i)$ of the dual problem at $\varepsilon_i$, for every $i$. Then
$\lambda_i \equiv \lambda(\varepsilon_i)$ is increasing as $\varepsilon_i\to\varepsilon$. If $\lambda_i \uparrow \lambda^* < \infty$ when $\varepsilon_i\downarrow\varepsilon$, then $\lambda^*$ is a solution of the dual problem at $\varepsilon$.
\end{proposition}
\paragraph{Proof of Proposition \ref{prop:LM_limit}}
Since strong duality holds at $\varepsilon_i$
\begin{eqnarray*}
V(\varepsilon_i) =
\min_{\lambda\geq 0} \Big\{H(\lambda)+\lambda \varepsilon_i\Big\}  = H(\lambda(\varepsilon_i))+\lambda(\varepsilon_i) \varepsilon_i,
\end{eqnarray*}
where $H(\lambda)$ is a convex function of $\lambda\geq 0$, and hence is differentiable at almost every $\lambda\geq 0$.

Consider first the case that $\varepsilon_i\downarrow \varepsilon$ where $\varepsilon > 0$. Since $\lambda_i\equiv \lambda(\varepsilon_i)$ is a solution of the dual problem
\begin{align}
H'(\lambda_i^+) + \varepsilon_i  &\geq 0, \nonumber \\
H'(\lambda_i^-) - \varepsilon_i  &\geq 0.
\label{eq:dd}
\end{align}
Since $\varepsilon_i\downarrow \varepsilon$ and $\lambda_i \uparrow \lambda^*$ as $i\rightarrow\infty$, and the right derivative  is increasing in $\lambda$, it follows that
\begin{eqnarray*}
H'({\lambda^*}^+) + \varepsilon \geq \lim_{i\rightarrow\infty} \Big\{ H'({\lambda_i^+}) + \varepsilon_i \Big\} \geq 0.
\end{eqnarray*}
On the other hand, the left derivative $H'(\lambda^-)$ is left-continuous in $\lambda$ so
\begin{eqnarray*}
H'({\lambda^*}^-) - \varepsilon = \lim_{i\rightarrow\infty} \big\{ H'_-(\lambda(\varepsilon_i)) - \varepsilon_i \Big\} \geq 0.
\end{eqnarray*}
It follows that $\lambda^*$ is optimal for the dual problem at $\varepsilon$.

When $\varepsilon = 0$, we need only consider the derivative from the right at $\lambda(0)$. In particular, it follows from \eqref{eq:dd} that
\begin{center}
\hfill
$\displaystyle
H'({{\lambda^*}^+}) \geq \lim_{\varepsilon_i \downarrow 0} \Big\{H'_+(\lambda(\varepsilon_i)) - \varepsilon_i\Big\} \geq 0.
$
\hfill
$\Box$
\end{center}

\subsubsection{Proof of Proposition \ref{prop:eps0}}
For every $\lambda\geq 0$ we have
\begin{eqnarray*}
\sum_{i=1}^
n p_i\max_{z_i} \Big\{f(z_i)-f(Y_i)- \lambda\|z_i-Y_i\|\Big\} \geq 0.
\end{eqnarray*}
We characterize the optimal dual variables by finding the values of $\lambda$ such that the lower bound is attained.

If
\begin{eqnarray*}
\lambda \geq \max_{i=1,\cdots,\,n} \max_{z_i}\frac{f(z_i)-f(Y_i)}{\|z_i-Y_i\|}
\end{eqnarray*}
it follows that
\begin{eqnarray*}
f(z_i)-f(Y_i) - \lambda \|z_i-Y_i\|\geq 0
\end{eqnarray*}
for every $i$, with equality when $z_i=Y_i$. It follows that
\begin{eqnarray*}
\sum_{i=1}^
n p_i\max_{z_i} \Big\{f(z_i)-f(Y_i)- \lambda\|z_i-Y_i\|\Big\}=0
\end{eqnarray*}
so $\lambda$ is a solution of the dual problem.
If
\begin{eqnarray*}
\lambda < \max_{i=1,\cdots,\,n} \max_{z_i}\frac{f(z_i)-f(Y_i)}{\|z_i-Y_i\|}
\end{eqnarray*}
then 
there exists
 $i$ such that
\begin{eqnarray*}
\lambda <  \max_{z_i}\frac{f(z_i)-f(Y_i)}{\|z_i-Y_i\|}
\end{eqnarray*}
we can find $z_i$ such that $f(z_i)-f(Y_i) - \lambda \|z_i-Y_i\| > 0$. It follows that
\begin{eqnarray*}
\sum_{i=1}^
n p_i\max_{z_i} \Big\{f(z_i)-f(Y_i)- \lambda\|z_i-Y_i\|\Big\}>0
\end{eqnarray*}
so $\lambda$ is not a solution of the dual problem.
\hfill$\Box$

\subsubsection{Dual variables}
Note that the cost and constraint functionals
\begin{align*}
F(\gamma) & = \max_\gamma \int_z f(x,\,z) \Big(\sum_{i=1}^{n} \gamma_i({\rm d}z) \Big)\\
G(\gamma) & := \sum_{i=1}^{n}\int_z \|z-Y_i\|\gamma_i({\rm d}z)
\end{align*}
are linear in $\gamma$, so \eqref{eq:W2} is a convex optimization problem, and $V_{\rm w}(\varepsilon)$ is concave, increasing and differentiable in $\varepsilon$ almost everywhere \cite{
luenberger1997optimization}.

Since solutions of the dual problem of \eqref{eq:W2} are supergradients of the value function $V_{\rm w}(\varepsilon)$ \cite{luenberger1997optimization}, we study worst-case sensitivity $V_{\rm w}'(0^+)\equiv\lim_{\varepsilon\downarrow 0}(V_{\rm w}(\varepsilon)-V_{\rm w}(0))/\varepsilon$ by studying dual solutions when $\varepsilon\downarrow 0$.

Let $\lambda\geq 0$ be the Lagrange multiplier for the Wasserstein constraint. The dual problem is
\begin{align}
\lefteqn{\min_{\lambda\geq 0} \max_{\gamma \in{\mathcal X}} \Big\{ \sum_{i=1}^{n} \int_z f(z) \gamma_i({\rm d}z)
+ \lambda\Big(\varepsilon- \sum_{i=1}^n\int_z \|z_i - Y_i\|\gamma_i({\rm d}z)\Big)} \nonumber  \\[5pt]
& = \min_{\lambda\geq 0} \max_{\gamma \in {\mathcal X}} \sum_{i=1}^{n}\int_{z_i}\Big[f(z_i)-\lambda \|z_i - Y_i\| \Big]\gamma_i({\rm d}z_i) + \lambda \varepsilon
\label{eq:D1}
\end{align}
where to ease notation, we drop the decision variable from the notation and write $f(z) \equiv f(x,\,z)$. This can be written
\begin{align*}
\lefteqn{\min_{\lambda\geq 0} \Big\{\sum_{i=1}^{n} p_i \max_{z_i}\Big\{f(z_i) -\lambda \|z_i - Y_i\|\Big\} + \lambda\varepsilon\Big\}} \\
& = \sum_{i=1}^{n} p_i f(Y_i) + \min_{\lambda\geq 0} \Big\{\sum_{i=1}^{n} p_i \max_{z_i}\Big\{f(z_i) - f(Y_i) -\lambda \|z_i - Y_i\|\Big\} + \lambda\varepsilon \Big\}.
\end{align*}
Intuitively, for every given transportation cost $\lambda$, the inner maximization in \eqref{eq:D1} defines a worst-case measure that moves probability mass $p_i$ from $Y_i$ to
\begin{align*}
z_i^* &= \argmax_{z_i}\Big\{f(z_i) -\lambda \|z_i - Y_i\|\Big\}.
\end{align*}
If $\lambda(\varepsilon)$ is any solution of the dual problem at $\varepsilon$,  $\lambda(\varepsilon)$ is a super-gradient of $V_{\rm w}(\varepsilon)$ at $\varepsilon$ \cite{
luenberger1997optimization}, so concavity of $V_{\rm w}(\varepsilon)$ means that worst-case sensitivity $V_{\rm w}'(0^+)\leq \lambda(0)$. We compute the sensitivity at $\varepsilon=0$ by characterizing the solutions of the dual problem of \eqref{eq:W2} when $\varepsilon=0$, and showing that one of these actually equals the right derivative $V_{\rm w}'(0^+)$.

By Lemma \ref{lemma:duality_prop}, strong duality holds when $\varepsilon>0$. The following result shows that  strong duality also holds when $\varepsilon=0$, and characterizes the set of optimal dual variables.
\begin{proposition} \label{prop:eps0}
Assume that there exists constant $L$ such that $|f(z)-f(Y_i)|\leq L\|z-Y_i\|$ for every $z$ and $i=1,\cdots,\,n$. Then
\begin{align*}
\min_{\lambda\geq 0}\sum_{i=1}^n
 p_i\max_{z_i} \Big\{f(z_i)-f(Y_i)- \lambda\|z_i-Y_i\|\Big\}=0
\end{align*}
and strong duality holds when $\varepsilon=0$
\begin{align*}
V_{\rm w}(0) & = \sum_{i=1}^{n}p_i f(Y_i) + \min_{\lambda\geq 0}\sum_{i=1}^n
p_i\max_{z_i} \Big\{f(z_i)-f(Y_i)- \lambda\|z_i-Y_i\|\Big\}  \\
             & = \sum_{i=1}^{n}p_i f(Y_i),
\end{align*}
and hence for all $\varepsilon\geq 0$.
The set of optimal solutions of the dual problem when $\varepsilon=0$  is
\begin{align*}
\nonumber
\Big\{\lambda\,\Big|\, \lambda \geq \max_{i=1,\cdots,\,n} \max_{z_i}\frac{f(z_i)-f(Y_i)}{\|z_i-Y_i\|}\Big\}
= \argmin_{\lambda\geq 0} \sum_{i=1}^n
 p_i\max_{z_i} \Big\{f(z_i)-f(Y_i)- \lambda\|z_i-Y_i\|\Big\}.
\end{align*}
\end{proposition}
Proposition \ref{prop:eps0} implies that if
\begin{align*}
\lambda \geq \max_{i=1,\cdots,\,n} \max_{z_i}\frac{f(z_i)-f(Y_i)}{\|z_i-Y_i\|},
\end{align*}
then $\lambda$ is a supergradient of $V_{\rm w}$ at $\varepsilon=0$, and hence is an upper bound of the right derivative
\begin{align}
V_{\rm w}'(0^+) \leq \max_{i=1,\cdots,\,n} \max_{z_i}\frac{f(z_i)-f(Y_i)}{\|z_i-Y_i\|}.
\label{eq:Wass_temp1}
\end{align}
The following result shows that this inequality is actually an equality, so $V_{\rm w}'(0^+)$ is also a solution of the dual problem at $\varepsilon=0$. This allows us to identify the identify worst-case sensitivity with the lower bound of the set of dual solutions at $\varepsilon=0$.

\subsubsection{Proof of Proposition \ref{prop:Wass_sensitivity}}
Suppose that the right derivative of $V$ at $\varepsilon=0$ is finite (i.e., $|V'(0^+)|<\infty$). Let $\{\varepsilon_i\}$ be any sequence such that $\varepsilon_i>0$, $\varepsilon_i\downarrow 0$ and $V$ is differentiable at $\varepsilon_i$. The Lagrange Duality Theorem \cite{luenberger1997optimization} implies that strong duality holds and that there exists a unique solution $\lambda_i\equiv\lambda(\varepsilon_i)$ of the dual problem for each $\varepsilon_i$. Since $\{\lambda(\varepsilon_i)\}$ is an increasing sequence that is bounded above by $0\leq V'(0^+) <\infty$ (we have selected $\{\varepsilon_i\}$ such that $V$ is differentiable at $\varepsilon_i$, so $\lambda(\varepsilon_i)= V'(\varepsilon_i)$; in addition, $0 \leq V'(\varepsilon_i) \leq V'(0^+)$ by the concavity of $V(\varepsilon)$), it converges to a limit that is also bounded by $V'_+(0)$
\begin{eqnarray*}
\lambda^* = \lim_{i\rightarrow\infty} \lambda(\varepsilon_i) =   \lim_{i\rightarrow\infty} V'(\varepsilon_i)  \leq V'(0^+) < \infty.
\end{eqnarray*}
By Proposition \ref{prop:LM_limit}, the limit $\lambda^*$ is a solution of the dual problem at $\varepsilon=0$. Since \begin{eqnarray*}
\left[\max_{i=1,\cdots,\,n} \max_{z_i}\frac{f(z_i)-f(Y_i)}{\|z_i-Y_i\|},\,\infty\right)
\end{eqnarray*}
is the set of Lagrange multipliers that solve the dual problem when $\varepsilon=0$,
\begin{eqnarray*}
 V'(0^+) \geq \lambda^*  \geq \max_{i=1,\cdots,\,n} \max_{z_i}\frac{f(z_i)-f(Y_i)}{\|z_i-Y_i\|}
\end{eqnarray*}
Together with \eqref{eq:Wass_temp1}, it now follows that worst-case sensitivity satisfies \eqref{eq:Wass-sensitivity}.
\hfill$\Box$

\section{Comparison of sensitivity measures}\label{sec:comparingWCS}
This section compares the WCS formulas derived in Sections \ref{sec:phi-div} and \ref{sec:nonsmooth_phi} when the nominal distribution is empirical one ($\mathsf{p}=\mathsf{1}/n$). 

It is known (e.g., \cite{shiffler1980}) that for any $\mathsf{f}$ and $\mathsf{p}=\mathsf{1}/n$,
\[
\sqrt{\mathbb{V}_{\mathsf{p}}(\mathsf{f})}\leq \frac{1}{2}\mbox{Range}(\mathsf{f}),
\]
and, accordingly, we have
\begin{eqnarray}
\sqrt{\frac{\phi''(1)}{2}}V'_{\phi}(0^+)\leq V'_{{\rm TV}}(0^+)
\label{eq:phi-TV bound}
\end{eqnarray}
where $V'_\phi(0+)$ and $V'_{\rm TV}(0+)$ denote the worst-case sensitivity \eqref{eq:phi-sensitivity} for smooth $\phi$ and \eqref{eq:sen_tv} for total variation, respectively.

We can associate CVaR deviation with the standard deviation.
\begin{proposition}
\label{propo:M-Stdev_as_UB_of_CVaR}
Let $\mathsf{f}\in\mathbb{R}^n$ and $\alpha\in(0,1)$. 
For $\mathsf{p}=\mathsf{1}/n$,
we have
\begin{equation}
\mbox{``CVaR Deviation of }\mathsf{f}\mbox{''}~\equiv~
\mathrm{CVaR}_{\mathsf{p},\alpha}(\mathsf{f})-\mathbb{E}_{\mathsf{p}}(\mathsf{f})\leq C_{\alpha,n}\sqrt{\mathbb{V}_{\mathsf{p}}(\mathsf{f})},
\label{eq:M-Stdev_as_UB_of_CVaR}
\end{equation}
where 
\[
C_{\alpha,n}:=\frac{\sqrt{n\Big\{\lfloor\kappa\rfloor+\big(\kappa-\lfloor\kappa\rfloor\big)^2\Big\}-\kappa^2}}{\kappa}
\]
with $\kappa:=n(1-\alpha)$.
The inequality \eqref{eq:M-Stdev_as_UB_of_CVaR} is tight, i.e., there is a vector $\mathsf{f}$ which attains the equality.
\end{proposition}
\paragraph{Proof of Proposition \ref{propo:M-Stdev_as_UB_of_CVaR}}
Deriving the inequality
\eqref{eq:M-Stdev_as_UB_of_CVaR}
is equivalent to showing the reciprocal of the minimum of the following optimization problem is equal to 
$1/C_{\alpha,n}$
\begin{equation}
\begin{array}{r|ll}
\frac{1}{C_{\alpha,n}}
=&\underset{\mathsf{f}}{\mbox{minimize}}&\displaystyle\frac{\sqrt{\mathbb{V}_{\mathsf{p}}(\mathsf{f})}}{\mathrm{CVaR}_{\mathsf{p},\alpha}(\mathsf{f})-\mathbb{E}_{\mathsf{p}}(\mathsf{f})}\\
&\mbox{subject to}&\mathsf{f}\neq C\mathsf{1}\mbox{ for any }C,
\end{array}
\label{eq:frac1}
\end{equation}
where $\mathsf{p}=\mathsf{1}/n$.
Noting that both the denominator and numerator are positively homogeneous and that 
CVaR is translation invariant and, thus, $\mathrm{CVaR}_{\mathsf{p},\alpha}(\mathsf{f})-\mathbb{E}_{\mathsf{p}}(\mathsf{f})=\mathrm{CVaR}_{\mathsf{p},\alpha}(\mathsf{f}-\mathbb{E}_{\mathsf{p}}(\mathsf{f})\mathsf{1})$, the fractional program \eqref{eq:frac1} is equivalently rewritten as
\[
\begin{array}{r|ll}
\frac{1}{C_{\alpha,n}^2}
=&\underset{\mathsf{f},\mathsf{z}}{\mbox{minimize}}&\displaystyle\frac{1}{n}\mathsf{z}^\top\mathsf{z}\\
&\mbox{subject to}&\mathrm{CVaR}_{\mathsf{p},\alpha}(\mathsf{z})=1,\\
&&\mathsf{z}=\mathsf{f}-\frac{1}{n}\mathsf{1}\mathsf{1}^\top\mathsf{f},~\mathsf{f}\neq C\mathsf{1}\mbox{ for any }C,
\end{array}
\]
where the objective function is squared for the convenience. 
Noting that the final constraint implies $
\mathsf{1}^\top\mathsf{z}=0$,
we consider the following relaxed optimization problem:
\begin{equation}
\begin{array}{r|ll}
&\underset{\mathsf{z}}{\mbox{minimize}}&\frac{1}{n}(z_1^2+\cdots+z_n^2)\\
&\mbox{subject to}&\displaystyle\frac{1}{\kappa}\{z_1+\cdots+z_{k}+(\kappa-k)z_{k+1}\}=1,\\
&&z_1\geq z_2\geq \cdots\geq z_n,\\
&&z_1+\cdots+z_n=0,
\end{array}
\label{eq:relaxed_problem_in_proof}
\end{equation}
where $\kappa:=n(1-\alpha)$ and $k:=\lfloor\kappa\rfloor$.

Let us tentatively consider to remove the inequality constraints from \eqref{eq:relaxed_problem_in_proof}:
\[
\begin{array}{r|ll}
&\underset{\mathsf{z}}{\mbox{minimize}}&\frac{1}{n}(z_1^2+\cdots+z_n^2)\\
&\mbox{subject to}&z_1+\cdots+z_{k}+(\kappa-k)z_{k+1}=\kappa,\\
&&(1-\kappa+k)z_{k+1}+z_{k+2}+\cdots+z_{n}=-\kappa.\\
\end{array}
\]
The Lagrangian of this is given by
\[
L(\mathsf{z},\lambda,\mu)=z_1^2+\cdots+z_n^2
-2\lambda\{z_1+\cdots+z_{k}+(\kappa-k)z_{k+1}-\kappa\}
-2\mu\{(1-\kappa+k)z_{k+1}+z_{k+2}+\cdots+z_{n}+\kappa\}
\]
and the optimality condition is given by
\[\renewcommand{\arraystretch}{1.5}
\left\{
\begin{array}{ll}
\displaystyle\frac{\partial L}{\partial z_j}(\mathsf{z},\lambda,\mu)=0,&j=1,...,n,\\
\displaystyle\frac{\partial L}{\partial \lambda}(\mathsf{z},\lambda,\mu)=0&\leftrightarrow z_1+\cdots+z_{k}+(\kappa-k)z_{k+1}=\kappa,\\
\displaystyle\frac{\partial L}{\partial \mu}(\mathsf{z},\lambda,\mu)=0&\leftrightarrow (1-\kappa+k)z_{k+1}+z_{k+2}+\cdots+z_{n}=-\kappa.\\
\end{array}
\right.
\]
Here the first condition is
\[
\frac{\partial L}{\partial z_j}(\mathsf{z},\lambda,\mu)=
\left\{
\begin{array}{ll}
2z_j-2\lambda=0,&1\leq j\leq k,\\
2z_{k+1}-2(\kappa-k)\lambda-2(1-\kappa+k)\mu=0,&j=k+1,\\
2z_j-2\mu=0,&k+2\leq j\leq n.
\end{array}
\right.
\]
Substituting $z_j$'s derived from this to the last two conditions, 
we have
\[
\left\{
\begin{array}{rrl}
\displaystyle \{k+(\kappa-k)^2\}\lambda&+(\kappa-k)(1-\kappa+k)\mu&=\kappa,\\
\displaystyle (\kappa-k)(1-\kappa+k)\lambda&+\{n-k-1+(1-\kappa+k)^2\}\mu&=-\kappa,
\end{array}
\right.
\]
and
\[
\left\{
\begin{array}{ll}
\displaystyle \lambda=\frac{\kappa(n-\kappa)}{n\{k+(\kappa-k)^2\}-\kappa^2}\\
\displaystyle \mu=\frac{-\kappa^2}{n\{k+(\kappa-k)^2\}-\kappa^2},\\
\end{array}
\right.
\]
so
\begin{equation}
\left\{
\begin{array}{lll}
\displaystyle z_j&\displaystyle =\frac{\kappa(n-\kappa)}{n\{k+(\kappa-k)^2\}-\kappa^2},&j=1,...,k,\\
\displaystyle z_{k+1}&\displaystyle =\frac{-\kappa^2+n\kappa(\kappa-k)}{n\{k+(\kappa-k)^2\}-\kappa^2},&\\
\displaystyle z_j&\displaystyle =\frac{-\kappa^2}{n\{k+(\kappa-k)^2\}-\kappa^2},&j=k+2,...,n,\\
\end{array}
\right.
\label{eq:ratio_kkt_sol}
\end{equation}
and
the objective value is then
\begin{equation}
\frac{1}{n}\mathsf{z}^\top\mathsf{z}=\frac{\kappa^2}{n\{k+(\kappa-k)^2-\kappa^2\}}.
\label{eq:squared_obj_val}
\end{equation}
The vector
$\mathsf{z}=(z_1,...,z_n)^\top$ given by \eqref{eq:ratio_kkt_sol} satisfies 
$z_1\geq\cdots\geq z_n$, and turns out to be optimal to \eqref{eq:relaxed_problem_in_proof} as well.
Furthermore,
for the vector $\mathsf{z}$ given by \eqref{eq:ratio_kkt_sol} and any constant $C$,
we can reproduce a nonconstant vector $\mathsf{f}$ as $\mathsf{f}=\mathsf{z}+C\mathsf{1}$, which satisfies $\mathsf{z}=\mathsf{f}-\frac{1}{n}\mathsf{1}\mathsf{1}^\top\mathsf{f}$. 
Consequently, 
the square root of 
\eqref{eq:squared_obj_val} is the optimal value of
\eqref{eq:frac1}.
\hfill$\Box$

\begin{remark}
Note that $C_{\alpha,n}\leq\sqrt{\frac{\alpha}{1-\alpha}}$ for all $\alpha\in[0,1)$, and especially when $n(1-\alpha)\in\mathbb{Z}$, the equality holds. Accordingly, \eqref{eq:M-Stdev_as_UB_of_CVaR} suggests a tighter version of nonparametric bound of  
$\alpha$-CVaR with Mean--Standard Deviation than \cite{rockafellar2014superquantile}:
\begin{align*}
\mathrm{CVaR}_{\mathsf{p},\alpha}(\mathsf{f})&\leq \mathbb{E}_{\mathsf{p}}(\mathsf{f})+\sqrt{\frac{\alpha}{1-\alpha}}\sqrt{\mathbb{V}_{\mathsf{p}}(\mathsf{f})},
\end{align*}
While Proposition 1 of \cite{rockafellar2014superquantile} shows a similar bound for random variables in the $L^2$-space, their coefficient is $1/\sqrt{1-\alpha}$, which is larger than $C_{\alpha,n}$.
\begin{figure}[h]\centering
\includegraphics[scale=0.75]{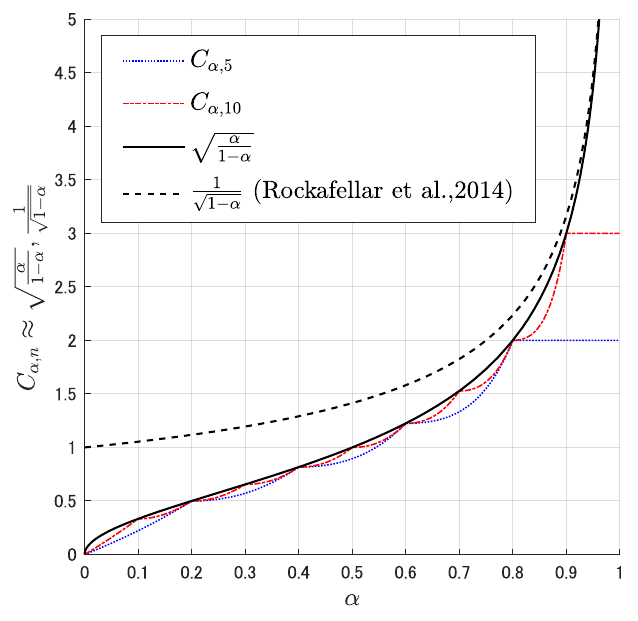}
\caption{$C_{\alpha,n}$, $\sqrt{\frac{\alpha}{1-\alpha}}$, and $\frac{1}{\sqrt{1-\alpha}}$}
\end{figure}
\end{remark}

The inequality \eqref{eq:M-Stdev_as_UB_of_CVaR} is applicable  to the worst-case sensitivity results. 
Let $V_{\rm b}'(0^+)$ and $V_{\rm c}'(0^+)$ denote the worst-case sensitivities \eqref{eq:sen_cvar} and \eqref{eq:sen-CVaRdev}, respectively. 
First of all,
\begin{eqnarray}
V'_{{\rm c}}(0^+)\leq C_{\alpha,n}\sqrt{\frac{\phi''(1)}{2}}V'_{\phi}(0^+) \leq \sqrt{\frac{\phi''(1)\alpha}{2(1-\alpha)}}V'_{\phi}(0^+).
\label{eq:convex-phi bound}
\end{eqnarray}
This inequality suggests that the sensitivity of the DRO with the convex combination of mean and CVaR is bounded above by that with any smooth $\phi$ divergence.

Recalling Proposition \ref{cor:sen_cvar}, we have a tight bound:
\begin{align*}
V'_{\rm b}(0^+)=
\mathbb{E}_{\mathsf{p}}(\mathsf{f})-\min(\mathsf{f})
&=\frac{\mathrm{CVaR}_{\mathsf{p},\frac{\varepsilon}{1+\varepsilon}}(\mathsf{f})-\mathbb{E}_\mathsf{p}(\mathsf{f})}{\varepsilon} 
\\
&\leq \varepsilon^{-1/2}\sqrt{\mathbb{V}_{\mathsf{p}}(\mathsf{f})}
=\sqrt{\frac{\phi''(1)}{2\varepsilon}}V'_{\phi}(0^+)
\end{align*}
for $\mathsf{p}=\mathsf{1}/n$.
Since \eqref{eq:sen_cvar} holds for any $\varepsilon\in(0,\frac{p_{(n)}}{1-p_{(n)}})=(0,\frac{1}{n-1})$, taking $\varepsilon=\frac{1}{n-1}$, we have a (loose) bound
\[
V'_{\rm b}(0^+)<\sqrt{\frac{(n-1)\phi''(1)}{2}}V'_{\phi}(0^+).
\]
From this we see that when $n$ is large, the difference between the smooth $\phi$ and the budgeted uncertainty $\phi=\delta_{[0,1+\varepsilon]}$ can be large for small uncertainty sets.
In contrast, the bounds \eqref{eq:phi-TV bound} and \eqref{eq:convex-phi bound} relating the sensitivities $V'_{\rm TV}(0^+)$, $V'_{\rm c}(0^+)$ and $V'_{\phi}(0^+)$ are tight and independent of $n$.
The potentially large difference likely reflects the fact that the sensitivity of budgeted uncertainty depends only on the lower part of the cost distribution whereas  $V'_{\phi}(0^+)$ depends on the entire distribution.
More generally, this suggests the possibility that solutions of DRO problems with budgeted uncertainty may differ quite substantially from those for other uncertainty sets.

\section{Logistic regression}
\label{sec:LR}
%
The ordinary logistic regression is SAA where the cost of the $i$-th sample is defined as
\begin{align*}
f(\mathsf{w},(y_i,\mathsf{x}_i))&:=\ln\big(1+\exp(-y_i\mathsf{w}^\top\mathsf{x}_i)\big),
\end{align*}
where $\mathsf{x}_i\in\mathbb{R}^{d}$ and $y_i\in\{\pm 1\}$ denote the input vector and binary label of the $i$-th sample, respectively, and $\mathsf{w}\in\mathbb{R}^{d}$ is the vector of decision variables. 

 For the Wasserstein DRO, \cite{abadeh2015distributionally} shows that when the transportation cost is given by $c((\mathsf{x}_i,y_i),(\mathsf{x}_j,y_j))=\|\mathsf{x}_i-\mathsf{x}_j\|_2+\kappa|y_i-y_j|$ for a constant $\kappa\geq 0$, the robust counterpart can be reduced to the following convex optimization problem:
\[
\begin{array}{lll}
\underset{\lambda\geq 0,\mathsf{s},\mathsf{w}}{\mbox{minimize}}
 & \varepsilon\lambda + \frac{1}{n}\sum\limits_{i=1}^ns_i\\
\mbox{subject to}& s_i\geq \ln\Big(1+\exp\big(-y_i(\mathsf{x}_i^\top\mathsf{w})\big)\Big),             & i=1,...,n,\\
                 & s_i\geq \ln\Big(1+\exp\big(y_i(\mathsf{x}_i^\top\mathsf{w})\big)\Big)-\kappa\lambda,& i=1,...,n,\\
                 & \lambda \geq \|\mathsf{w}\|_2.
\end{array}
\]
If the decision maker is not interested in the mislabeling, s/he can set as $\kappa=\infty$, so that the formulation results in the regularized logistic regression:
\begin{align*}
\underset{\lambda\geq 0,\mathsf{s},\mathsf{w}}{\mbox{minimize}}\quad
 & \varepsilon\|\mathsf{w}\|_2 + \frac{1}{n}\sum\limits_{i=1}^n\ln\Big(1+\exp\big(-y_i(\mathsf{x}_i^\top\mathsf{w})\big)\Big).
\end{align*}
Then the worst-case sensitivity is given by $\mathcal{S}_{\mathsf{p}}(\mathsf{f})=\|\mathsf{w}_{\rm SAA}\|_2
$ where 
$\mathsf{w}_{\rm SAA}$ is a solution to the ordinary logistic regression. 

Figure \ref{fig:lr:hearfail:frontiers} shows the 
mean--smooth $\phi$-sensitivity frontiers of the seven DROs for the logistic regression using the heart failure clinical records dataset \cite{chicco2020machine}, which consists of 299 samples having 12 covariates.
For this example, all frontiers 
are similar, and the choice of uncertainty set is not expected to make a significant difference. 
%
%
%
%
\begin{figure}[h]
\includegraphics[scale=0.5]{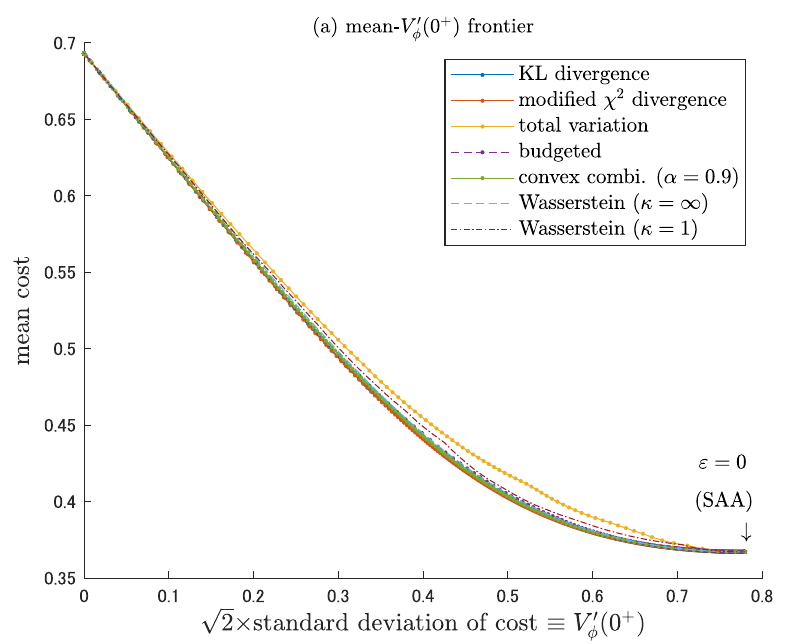}
\caption{Mean--sensitivity frontiers of solutions of the robust logistic regression problem with the Heart Failure data set}
\label{fig:lr:hearfail:frontiers}
{\raggedright\footnotesize
Each frontier corresponds to the family of DRO solutions for one of the seven uncertainty sets. Worst-case sensitivity in all the frontier plots is measured by $\sqrt{2}\times{\mbox{standard deviation}}$ (worst-case sensitivity of the modified $\chi^2$ uncertainty set) of the associated in-sample cost distribution. 
\par}
\end{figure}

\end{document}